\documentclass[12pt]{article}

\usepackage{arxiv}
\usepackage[numbers]{natbib}
\usepackage[ruled,linesnumbered]{algorithm2e}
\usepackage[utf8]{inputenc} 
\usepackage[T1]{fontenc}    
\usepackage{hyperref}       
\usepackage{url}            
\usepackage{booktabs}       
\usepackage{amsfonts}       
\usepackage{nicefrac}       
\usepackage{microtype}      
\usepackage{lipsum}		
\usepackage{graphicx}
\usepackage{natbib}
\usepackage{doi}
\usepackage{color}
\usepackage{amsmath}
\usepackage{amssymb}
\usepackage{amsthm}
\usepackage{algorithmic}
\newtheorem{assumption}{Assumption}

\newtheorem{theorem}{Theorem}
\newtheorem{definition}{Definition}
\newtheorem{lemma}{Lemma}

\usepackage{lyu}
\DeclareMathOperator{\Diag}{Diag}
\author{
   Liyao Lyu\\
   Department of Computational Mathematics, Science \& Engineering\\ Michigan State University, MI 48824, USA
    \And
    Jingrun Chen\\
   School of Mathematical Sciences and Suzhou Institute for Advanced Research\\
   University of Science and Technology of China, Suzhou 215127, China 
    }
\title{Consensus Based Stochastic Control}

\begin{document}
\onecolumn
\maketitle
We propose a \emph{gradient-free} deep reinforcement learning algorithm to solve \emph{high-dimensional}, finite-horizon stochastic control problems. 
Although the recently developed deep reinforcement learning framework has achieved great success in solving these problems, direct estimation of policy gradients from Monte Carlo sampling often suffers from high variance.  To address this, we introduce the Momentum Consensus-Based Optimization (M-CBO) and Adaptive Momentum Consensus-Based Optimization (Adam-CBO) frameworks. These methods optimize policies using Monte Carlo estimates of the value function, rather than its gradients. Adjustable Gaussian noise supports efficient exploration, helping the algorithm converge to optimal policies in complex, nonconvex environments. Numerical results confirm the accuracy and scalability of our approach across various problem dimensions and show the potential for extension to mean-field control problems. Theoretically, we prove that M-CBO can converge to the optimal policy under some assumptions.
\section{Introduction}

Stochastic optimal control (SOC) problems~\cite{stengel1986stochastic,fleming2012deterministic}, along with their mean-field variants, have been extensively studied throughout the twentieth century and have had a wide range of applications in various areas, such as  finance~\cite{pham2009continuous,fleming2004stochastic,carmona2003pricing,cousin2011mean,lachapelle2016efficiency,cardaliaguet2018mean}, economics~\cite{gueant2009mean,gomes2015economic,gueant2010mean,achdou2014partial,achdou2022income}, chemistry~\cite{welch2019describing,holdijk2024stochastic}, and biology~\cite{lachapelle2011mean,aurell2018mean,achdou2019mean}. Readers seeking an overview of these developments may refer to the recent review~\cite{hu2024recent}. Traditional methods for solving the SOC problem, such as the finite-volume method~\cite{richardson2006numerical,wang2003numerical}, the Galerkin method~\cite{beard1997galerkin,bea1998successive}, and the monotone approximation method~\cite{forsyth2007numerical}, aim to solve the corresponding Hamilton-Jacobi-Bellman (HJB) equations. However, these methods struggle to scale in high-dimensional spaces due to the curse of dimensionality, where the computational complexity grows exponentially with the dimension of state and action variables. This limitation hinders their application in large-scale systems where efficiency is critical.

Significant advances have been made in addressing the high-dimensional SOC problem by modeling control strategies using deep neural networks, leveraging their capability to approximate functions in high-dimensional spaces. One prominent approach is the value-based method
~\cite{Li2024value,Lien2024Enhancing,Obando2024value,zhang2024model,mou2024bellman}, such as the deep-backward stochastic differential equation (BSDE) method~\cite{e2017deep,han2018solving,nusken2021solving,pham2021neural}. Based on the Bellman principle, the optimal control can be modeled as a function of the value function and its gradient. Therefore, solving the value function from the BSDE that it satisfies can automatically give the optimal control of the SOC. These methodologies are commonly referred to as \emph{model-based} methods because they need an explicit connection between the optimal control and the value function. This kind of connection usually depends on accurate modeling of the transition kernel between different states. However, modeling the (mean-field) transition kernel for a real-world process in practical applications can be extremely challenging~\cite{lyu2023construction,lu2019nonparametric}.

Recently, \emph{model-free} methods have gained attention in control and reinforcement learning~\cite{agrawalpolicy,chenaccelerated,chenE32024,Daisafe2024,Hisaki2024Average,Hong2024model,Hu2024Q,Park2024Max,tang2024solving}, such as Deep Q Networks~\cite{mnih2015human}, Proximal Policy Optimization~\cite{heess2017emergence,schulman2015high,schulman2017proximal}, Trust Region Policy Optimization~\cite{schulman2015trust}, Deep Deterministic Policy Gradient~\cite{silver2014deterministic,lillicrap2015continuous} and Soft Actor-Critic~\cite{haarnoja2018soft,haarnoja2018learning}.  These approaches address this issue by directly optimizing the policy without explicit transition kernel modeling. 
Nevertheless, these methods rely on the evaluation of policy gradients~\cite{jia2022actor-critic,jia2022temporal-difference} or depend on the action and state space discretization~\cite{gu2021mean,carmona2023model}. The evaluation of policy gradients often has high variance and is computationally intensive~\cite{hua2024simulation}, and the discretization of action and state space reintroduces dimensionality constraints.
Consequently, current approaches face a trade-off between model fidelity and scalability, motivating the need for a method that can achieve robust performance without gradient estimation and state-action discretization.

In this work, we introduce a novel approach to overcome the limitations of both model-based and model-free reinforcement learning methods by applying the Adam-CBO~\cite{chen2022consensus} framework to high-dimensional SOC problems. Unlike value-based methods, our approach is entirely model-free, directly optimizing the policy without requiring an explicit formulation of the transition kernel. In addition, it is gradient-free, avoiding the high-variance issue associated with policy gradients, and mesh-free, eliminating the need to discretize state and action spaces. These features allow our method to scale efficiently in high-dimensional environments, making it particularly suited for finite-horizon problems where the optimal control is time-dependent. Contrary to concerns that direct policy optimization may lead to local optima, our method demonstrates superior accuracy in handling nonconvex issues, as evidenced by extensive numerical results.

Beyond numerical validation, our study contributes a rigorous theoretical foundation by providing the convergence analysis for the M-CBO method, a simplified version of Adam-CBO without adaptive timestep. This proof establishes that, under certain assumptions, our algorithm reliably converges to the optimal policy, addressing a crucial gap in reinforcement learning for the SOC problem, where theoretical guarantees are often challenging to obtain.

\section{Problem Formulation}
Consider a control problem over a finite time horizon $t\in [0, T]$ for some $T<\infty$.
The state space is denoted by $\mathcal S \subset \mathbb R ^d$, and the action space by $\mathcal A \subset \mathbb R^m $.
An agent governs its state process $\bx_t$ through an action process $\ctrl_t$ with a transition kernel $p(\bx'|t,\bx,\ctrl)$ that describes the evolution from state $\bx$ to state $\bx'$ under the action $\ctrl$ at time $t$.  
The agent's goal is to minimize the combined terminal cost $g(\bx_T)$ and the running cost $f(t,\bx,\ctrl)$ incurred during the process.   
The total cost function is generally represented by \[J[\boldsymbol\alpha] = \mathbb{E}\l [\int_0^T f(t,\bx_t,\ctrl_t)\intd t + g(\bx_T)\r  ].\] 
In this work, we model the policy $\ctrl(t,\bx;\theta)$ as a fully connected neural network parameterized by $\theta$. The rest of the paper will focus on finding the optimal $\theta \in \mathbb{R}^{D}$ such that it minimizes the cost function $\mathcal J(\theta) = J[\ctrl[(t,\bx;\theta)]$.

\section{Gradient-free Policy Update}
\label{sec:gradient_free_policy_update}
We propose two algorithms to find the optimal policy: M-CBO and Adam-CBO. The Adam-CBO algorithm improves on M-CBO by adaptively adjusting the timestep, resulting in better numerical performance. 

\subsection{Momentum Consensus-Based Optimization}
In M-CBO,  we begin by initializing a population of $N$ agents represented by $\l (\boldsymbol{\Theta},\boldsymbol{\Omega}\r )= (\Theta^1,\Omega^1,\cdots,\Theta^N,\Omega^N) \in \mathbb{R}^{2ND}$. Here $\Theta^i \in \mbRd$ denotes the policy parameterization of the $i\mhyphen$th agent, and $\Omega^i\in \mbRd$ represents its momentum. To exploit the current group of policies, we estimate a consensus policy as 
\begin{equation*}
\mathcal{M}_\beta\l(\boldsymbol{\Theta}\r) = 
\sum_{i=1}^N\frac{\Theta^i w_\beta\l(\Theta^i\r)}{\sum_{j=1}^N w_\beta\l(\Theta^j\r)},
\end{equation*}
where $ w_\beta(\Theta) = \exp\l(-\beta \mathcal J(\Theta)\r)$.
Here $\beta\geq 0 $ is an inverse temperature parameter, controlling how strongly each agent's performance (determined by the objective function $\mathcal J (\Theta)$) influences the consensus. 
Using the consensus policy, we define the following dynamics to guide each policy toward consensus: 
\begin{equation}
\label{equ:CBO_SDE}
\begin{aligned}
        \intd \Theta^i_t =&   \Omega_t^i  \intd t  -\gamma_1 \l(\Theta^i_t - \mathcal{M}_\beta(\boldsymbol{\Theta})\r) + \sigma(t)\intd W^i_{\theta,t} ,\\ 
    \intd \Omega^i_t =&  - m \l(\Theta^i_t - \mathcal{M}_\beta(\boldsymbol{\Theta}) \r) \intd t 
    \\
    & -\gamma_2 \Omega^i_t\intd t + \sqrt{m} \sigma(t)\intd W^i_{\omega,t},
\end{aligned}
\end{equation}
where $m,\gamma_1$, and $\gamma_2$ are positive constants and $ W^i_{\theta,t}, W^i_{\omega,t}$ are $D$ dimensional Wiener processes that introduce stochasticity into the dynamics. This facilitates the exploration of unknown regions, with a parameter $\sigma(t)$ regulating the exploration strength. Using the Euler-Maruyama (EM) scheme for Equation \eqref{equ:CBO_SDE}, we get the M-CBO algorithm, as detailed in Algorithm \ref{alg:CBO-moment}.

\begin{algorithm}[tb]
\caption{ Consensus Based Optimization with Momentum}
\label{alg:CBO-moment}
\begin{algorithmic}
    \STATE {\bfseries Input:} {time step $\lambda$, Number of player $N$, Batch size $M$, total time $t_N$, parameters  $\beta$, $\gamma_1$, $\gamma_2$, $m$}
    \STATE Initialize  $\Theta^i_0\sim \mathcal N(0,\mathbb I_D)$, $i = 1, \dots, N$
    \STATE  Initialize $\Omega^i_0=0$, $i = 1, \dots, N$;\\
    \FOR{$t = 0$ {\bfseries to} $t_N$ }
    
    \STATE Partition the indices $\{1, 2, \dots, N\}$ into batches $B^1, \dots, B^{\frac{N}{M}}$, each containing $M$ particles
    
    \FOR{$j = 1$ {\bfseries to} $\frac{N}{M}$}
        \STATE  $\mathcal J^i = \mathcal J(\Theta^i_t)$, where $i\in B^j$
        
        \STATE  $M =\sum\limits_{k\in B^j} \frac{\Theta_t^k w^k}{\sum\limits_{i\in B^j} w^i}$, where $w^i = \exp\l ( -\beta \mathcal J^i\r ) $
        \STATE  Update the policies and their momentum:
        \[
        \begin{aligned}
        \Theta^i_{t+1} &= \Theta^i_t +\lambda \Omega^i_t - \gamma_1\lambda  (\Theta^i_t - M)+\sqrt{\lambda}\xi_{\theta}^i, \\
        \Omega^i_{t+1} &= \Omega^i_t - \lambda m (\Theta^i_t - M) -\lambda \gamma_2 \Omega^i_t  + \sigma(t)\sqrt{\lambda m}  \xi_{\omega}^i,
        \end{aligned}
        \]    where  $\xi_{\theta}^i,\xi_{\omega}^i \sim \mathcal N(0,\mathbb I_D) $
    \ENDFOR
    \ENDFOR
\STATE {\bfseries Output:} {$\Theta_{t_N}^i$, $i = 1, \dots, N$}
\end{algorithmic}
\end{algorithm}

The original CBO method~\cite{fornasier2024consensus2} aims to achieve a monotonic reduction in the distance between the optimal policy $\tilde\theta$ and the policies of agents. Specifically, this is represented as: $\frac{1}{N}\sum_{i=1}^N\|\Theta^i_t-\tilde\theta\|^2 \simeq \int\|\theta-\tilde\theta\|^2\intd \mu_t(\theta)$, where $\mu_t$ represents the law of agents $\boldsymbol\Theta_t$. Our method minimizes a combined expression $\frac{1}{N}\sum_{i=1}^N\l(\|\Theta^i_t-\tilde\theta\|^2  + m^{-1}\|\Omega^i_t\| \r)\simeq \int\|\theta-\tilde\theta\|^2 + m^{-1} \|\omega\|^2\intd \rho_t(\theta,\omega)$, where $\rho_t$ represents the joint distribution of policies $\boldsymbol\Theta$ and $\boldsymbol\Omega$ at time $t$. In particular, the M-CBO method does not force the monotonic reduction of $\frac{1}{N}\sum_{i=1}^N\|\Theta^i_t-\tilde\theta\|^2$,  allowing for the additional momentum term $\omega$  to enhance the exploration capability. It provides greater flexibility and reduces the risk of becoming trapped in local minima; see Section \ref{sec:convergence} for a more detailed analysis.

\subsection{Adaptive  Momentum Consensus-Based Optimization}
In the Adam-CBO method, we extend M-CBO by replacing the constant momentum term $m$ with an adaptive term based on the inverse of the second moment of the agents' policies. Specifically, we replace $m$ with $(\mathcal{V_\beta[\boldsymbol{\Theta}]}+\epsilon \mathbf{I})^{-1}$ , where $V_\beta[\boldsymbol\Theta]$ is the second moment defined as:
\[V_\beta[\boldsymbol{\Theta}]= \sum_{i=1}^N\frac{(\Theta^i -\mathcal{M}_\beta[\boldsymbol{\Theta}])^2w_\beta(\Theta^i)}{\sum_{j=1}^N w_\beta(\Theta^j)}.\]
In particular, $\epsilon=1\times10^{-8}$ is used to keep the positivity. This adaptive adjustment introduces a mechanism similar to the Adam optimizer, where updates are scaled by a normalized second moment, allowing for faster convergence and improved numerical performance. The detailed algorithm is shown in Algorithm \ref{alg:CBO-adam}.

\begin{algorithm}[tb]
\caption{Consensus-based Optimization with Adaptive Momentum}
\label{alg:CBO-adam}
\begin{algorithmic}
    \STATE {\bfseries Input:} time step $\lambda$, Number of player $N$, Batch size $M$, total time $t_N$, parameters  $\beta$, $\beta_1$, $\beta_2$
    \STATE Initialize $\Theta^i_0\sim \mathcal N(0,\mathbb I_D)$, $i = 1, \dots, N$ 
    \STATE  Initialize $\Omega^i_0=0$, $i = 1, \dots, N$
    \STATE Initialize $M_0, V_0 = 0$
    \FOR{$t = 0$ {\bfseries to} $t_N$  }
    
    \STATE Partition the indices $\{1, 2, \dots, N\}$  into batches $B^1, \dots, B^{\frac{N}{M}}$, each containing $M$ particles
    \FOR{$j = 1$ {\bfseries to} $\frac{N}{M}$}
    \STATE  $\mathcal J^i := \mathcal J(\Theta^i_t)$, where $i\in B^j$
        
    \STATE  $M =\sum\limits_{k\in B^j} \frac{\Theta_t^k w^k}{\sum\limits_{i\in B^j} w^i}$ 
        
    \STATE 
        $V = \sum\limits_{k\in B^j} \frac{(\Theta_t^k - M )^2w^k}{\sum\limits_{i\in B^j} w^i}$
    \STATE Update the moving average moment estimate:
        \[
        \begin{aligned}
        M_{t+1} = \beta_1 M_{t} +(1-\beta_1) M, \quad \hat{M}_{t+1} = \frac{M_{t+1}}{1-\beta_1^t},\\
        V_{t+1} = \beta_2 V_{t} +(1-\beta_2) V, \quad \hat{V}_{t+1} = \frac{V_{t+1}}{1-\beta_2^t}
        \end{aligned}
        \]
        
    \STATE  Update the policies and their momentum:
        \[
        \begin{aligned}
        \Theta^i_{t+1} =& \Theta^i_t +\lambda V^i_t, \\
        \Omega^i_{t+1} =& \Omega^i_t - \lambda \Diag(\hat{V_t^i} + \epsilon)^{-1} (\Theta^i_t - \hat{M_t}) 
        \\
        &+ \gamma\lambda \Omega^i_t + \sigma(t) \sqrt{\lambda}  \xi_i,
        \end{aligned}
        \] where  $\xi_{i} \sim \mathcal N(0,\mathbb I_D) $
    \ENDFOR
    \ENDFOR
\STATE {\bfseries Output: }{$\Theta_{t_N}^i$, $i = 1, \dots, N$}
    
\end{algorithmic}

\end{algorithm}

\section{Convergence Analysis}\label{sec:convergence}
In Section \ref{sec:gradient_free_policy_update}, we propose two dynamics that converge to the consensus policies. A natural question we want to answer here is whether policies can converge to the optimal policies. From the theoretical perspective, for simplicity, we focus on proving the convergence of the M-CBO method in this work. We begin by establishing the well-posedness of the M-CBO method, ensuring the uniqueness and existence of solutions under certain regularity conditions on the cost function $\mathcal J$.
\begin{assumption}\label{assum:Lip_J}
 The following assumptions are imposed on the cost function $\mathcal J $
\begin{enumerate}
    \item There exist $\tilde\theta$ such that $\mathcal{J}(\tilde\theta) =\inf_\theta \mathcal {J} (\theta) =: \underline{ J} $. Also, it is bounded from above by $\sup \mathcal{J} \leq \overline{J}$.
    \item The cost function $\mathcal J$ is locally Lipschitz continuous $\|\mathcal J[\theta_1]-\mathcal J [\theta_2]\| \leq L_J (\|\theta_1\|+\|\theta_2\|)\|\theta_1-\theta_2\|$.
    \item There exists a constant $c_{\mathcal J}>0$ such that $\mathcal J(\theta ) - \underline{J} \leq c_{\mathcal J}(1+\|\theta\|^2)$.
    \item 
    There exist $\delta_J,R_0,\eta >0$ such that 
    $
        \|\theta - \tilde\theta\| \leq \frac{\mathcal{J} -\underline{ J}}{\eta},
        $
        for all $ \theta \in B_{\theta,R_0}(\tilde\theta)=\{\theta :  \|\theta -\tilde\theta\|\leq R_0\} ,$
and     $
        \mathcal J (\theta) - \underline{ J } > \delta_J$  for all $ \theta \in \l ( B_{\theta,R_0}(\tilde\theta)\r ) ^c$.
    \item The parameters we choose $\sigma(t)$ has upper and lower bound $\underline{\sigma}\leq \sigma(t)\leq \overline{\sigma}$.
\end{enumerate}
\end{assumption}

\begin{theorem}\label{thm:well-posedness_m_cbo_fix_N}
    Under the Assumption \ref{assum:Lip_J}, for each $N\in \mathbb N$, the stochastic differential equation \eqref{equ:CBO_SDE} has a unique strong solution $\l\{\l (\bTheta\pNt, \bOmega \pNt)\r )|t>0\r\}$ for any initial condition $\l (\bTheta^{(N)}_0,\bOmega^{(N)}_0 \rg$ satisfying $\mathbb{E}\l ( \|\bTheta^{(N)}_0\|+\|\bOmega^{(N)}_0\|\r ) \leq \infty$.  
\end{theorem}
\begin{proof}
    See Appendix \ref{app:well-posedness_m_cbo_fix_N}.
\end{proof}
By letting the number of agents $N\to \infty$ in Equation~\eqref{equ:CBO_SDE}, the mean-field limit of the model is formally given by the following McKean–Vlasov stochastic differential equation
\begin{equation}
\label{equ:mean_M_CBO}
    \begin{aligned}
    \intd \bar{\Theta}_t =&    \bar\Omega_t  \intd t  -\gamma_1  \l (\bar\Theta_t - \mathcal{M}_\beta[\mu_t] \r )\intd t+\sigma(t)\intd W_{\theta,t},\\ 
     \intd 
 \bar{\Omega}_t = &
- m \l (\bar\Theta_t - \mathcal{M}_\beta[\mu_t] \r ) \intd t 
\\ & -\gamma_2 \bar \Omega_t \intd t + \sqrt{m} \sigma(t)\intd W_{\omega,t},
    \end{aligned}
\end{equation}
where 
$
\CM[\mu] = \frac{\int \theta \exp(-\beta \mathcal J(\theta)\mu(\intd \theta)}{\int  \exp(-\beta \mathcal J(\theta)\mu(\intd \theta)}
$, 
$\mu_t(\theta)  = \int \rho_t(\theta,\intd \omega)$, and $\rho_t = Law (\bar{\Theta}_t,\bar\Omega_t)$.
Then the corresponding Fokker-Planck equation is 
\begin{equation}
\label{equ:FP_M_CBO}
\begin{aligned}
    \partial_t \rho_t =&  -\nabla_\theta\cdot \l ( \l ( \omega-\gamma_1(\theta-\CM [\mu_t])\r ) \rho_t\r )   \\
    & + \nabla_\omega\cdot  \l ( \l( m \l (\theta - \CM \l [\mu_t \r] \r ) +\gamma_2 \omega\r)\rho_t\r )  \\
    & +\frac{\sigma(t)^2m}{2} \Delta_\omega \rho_t+\frac{\sigma(t)^2 }{2} \Delta_\theta \rho_t.
\end{aligned}
\end{equation}
Next, we will prove the above equation \eqref{equ:mean_M_CBO} and \eqref{equ:FP_M_CBO} are well-posed.
\begin{theorem}\label{thm:well_possness_mean_field}
    Let $\mathcal J$ satisfy the Assumption \ref{assum:Lip_J} and $\rho_0\in \PP{4}$. Then there exists a unique nonlinear process $(\bar{\Theta},\bar{\Omega}) \in   \MC{[0,T],\mbRd\times\mbRd} , T>0$, satisfying \eqref{equ:mean_M_CBO} with initial distribution $(\bar{\Theta},\bar{\Omega}) \sim \rho_0$ in the strong sense, and $\rho_t =  \text{Law} (\bar{\Theta},\bar{\Omega}) \in \MC{ [0,T],\PP{4}}$ satisfies the corresponding Fokker-Planck equation \eqref{equ:FP_M_CBO} in the weak sense with $\lim_{t\to \infty}\rho_t = \rho_0$ .
\end{theorem}
\begin{proof}
    See Appendix \ref{app:well-posedness_mean_field}.
\end{proof}
Then we present the result showing that \eqref{equ:mean_M_CBO}  and  \eqref{equ:FP_M_CBO}  model the mean-field limit of Equation \eqref{equ:CBO_SDE}.
\begin{theorem}\label{thm:mean_field_limit}
Let $\mathcal J $ satisfy Assumption \ref{assum:Lip_J} and $\rho_0\in \mathcal{P}_4(\mbRd\times\mbRd)$. For any $N\geq 2$, assume that $\{(\Theta_t^{(i,N)},\Omega_t^{(i,N)})_{t\in[0,T]}\}_{i=1}^N$ is the unique solution to the particle system \eqref{equ:CBO_SDE} with $\rho_0^{\otimes N}$-distributed initial data $\{(\Theta_0^{(i,N)},\Omega_0^{(i,N)})\}_{i=1}^N$.   Then the limit (denoted by $\rho$) of the sequence of the empirical measure $\rho^N=\frac{1}{N}\sum_{i=1}^N\delta_{\l (\Theta^{(i,N)},\Omega ^{(i,N)}\rg}$ exists. Moreover, $\rho$ is deterministic and it is the unique weak solution to PDE \eqref{equ:FP_M_CBO}.
\end{theorem}
\begin{proof}
    See in Appendix \ref{app:mean_filed}.
\end{proof}
To prove the global convergence of the M-CBO method, we define the energy functional as 
\begin{equation}\label{equ:energy_functional}
    E[\rho] = \frac{1}{2} \int \|\theta-\tilde\theta\|^2+ m^{-1}\|\omega\|^2 \intd\rho.
\end{equation}
The above definition $E[\rho]$ provides a measure of the distance between the distribution of the agents $\rho$ and the Dirac measure at $(\tilde\theta,0)$, denoted as $\delta_{(\tilde\theta,0)}$. Specifically, we have the relationship $2E[\rho_t] = W^2_2(\rho_t,\delta_{(\tilde\theta,0)} )$.

\begin{theorem}\label{thm:converge}
    Let $\mathcal J$ satisfy the Assumption \ref{assum:Lip_J}.  Moreover, let $\rho_{0} \in \mathcal P_4(\mathbb R^{2D})$ and $(\tilde\theta,0)\in supp(\rho_{0})$. By choosing parameters $\sigma(t)$ is exponentially decaying as $\sigma(t) = \sigma_1 \exp(-\sigma_2 t)$ with $\sigma_1>0$ and $\sigma_2>1$ and $\lambda = \max\{m,\gamma_1\} \geq 2\sigma_2 $ and $\gamma = \min\{\gamma_1,\gamma_2\}>0$. Fix any $\epsilon\in (0,E[\rho_0])$ and $\tau \in (0,1-\frac{2\sigma_2}{\lambda})$,  and define the time horizon 
    \begin{equation}
        T^* := \frac{1}{(1-\tau)\lambda}\log\l ( \frac{ E[\rho_{T_0}])}{\epsilon}\r )
    \end{equation}
    Then there exists $\beta >0$ such that for all $\beta>\beta_0$, if $\rho\in \mathcal C([0, T^*],\mathcal P_4(\mathbb R^{2D }))$ is a weak solution to the Fokker-Planck equation in the time interval $[0, T^*]$ with initial condition $\rho_0$, we have 
    \[
    \min_{t\in [0,T^*]}E[\rho_t] \leq \epsilon.
    \]Furthermore, until $E[\rho_t]$ reaches the prescribed accuracy $\epsilon$, we have the exponential decay 
    \begin{equation}
        E[\rho_t] \leq E[\rho_0]\exp(-(1-\tau)\lambda t)
    \end{equation}
    and, up to a constant, the same behavior for $W^2_2(\rho_t,\delta_{(\tilde\theta,0)})$.
\end{theorem}
\begin{proof}
    See Appendix \ref{app:converge}.
\end{proof}

\section{Numerical Results}
\label{sec:numerical}
We evaluate the performance of the Adam-CBO method across various problem settings, including the linear quadratic control problem in $1, 2, 4, 8,$ and $16$ dimensions, the Ginzburg-Landau model, and the systemic risk mean-field control problem with $50,100,200,400,800$ agents. Even though our method is model-free, which means it does not depend on the known explicit knowledge of the transition kernel as well as the precise dependency of the value function 
$u(t, \bx)$
on the optimal control  $\ctrl(t,\bx,\nabla u, \Hessian u)$. The value function is expressed as:
\[u(t,\bx)=  \inf_{\alpha\in \mathcal A}\mathbb E\l [\int_t^T  f(s,\bx_s,\ctrl_s)\intd s + g(X_T)\l| x(t)= x\r.\r].
\]
To measure the accuracy of our method, we compare $u(t,\bx)$ 
or the $\|\ctrl(t,\bx,\nabla u, \Hessian u)-\ctrl(t,\bx;\theta)\|$ as a metric. Our code is available at \url{https://github.com/Lyuliyao/Adam_CBO_Control}.

\subsection*{Linear Quadratic Control Problem}
We begin by considering a classical linear quadratic Gaussian (LQG) control problem.
The value function is known as
$u(t, \bx) = -\ln \l (   \mathbb{E} \l [ \exp \l (  - g \l (  \bx + \sqrt{2} \mathbf W_{T-t} \r )  \r )  \r  ] \r ) $, which we refer to Appendix \ref{sec:simulation_datail} for details. 
The numeric value of $u(t,\bx)$ can be computed by Monte Carlo (MC) estimation directly as a reference to measure the accuracy.

We investigate the LQG problem in dimension $d= 1, 2, 4, 8,$ and $16$, with a terminal time of $T=1$ and a timestep of $\frac{T}{20}$.
We compare our method with the BSDE method in \cite{han2018solving}. In both methods, the number of SDE to compute the value function is 64 and the learning rate is $1\times 10^{-2}$. In M-CBO and Adam-CBO methods, the number of agents is specified as $N = 5000$, and $M=50$ agents are randomly selected to update in each step.

\begin{figure}[ht]
\vskip 0.2in
\begin{center}
    \includegraphics[width=0.49\linewidth]{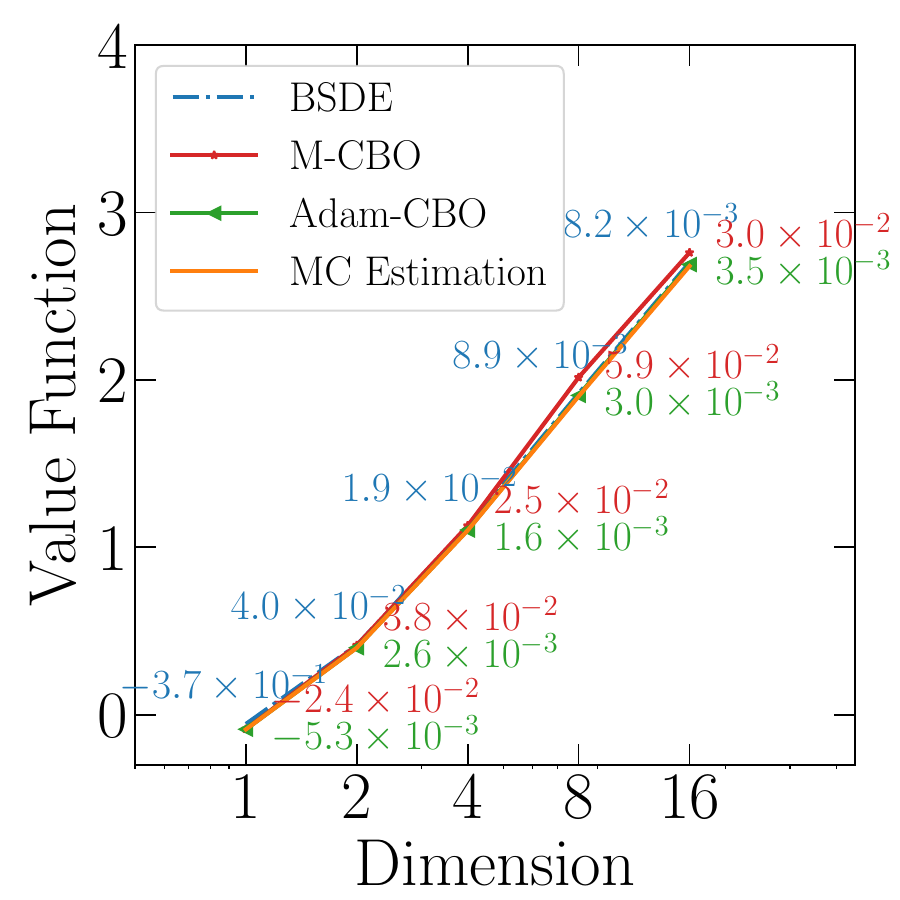}
    \includegraphics[width=0.49\linewidth]{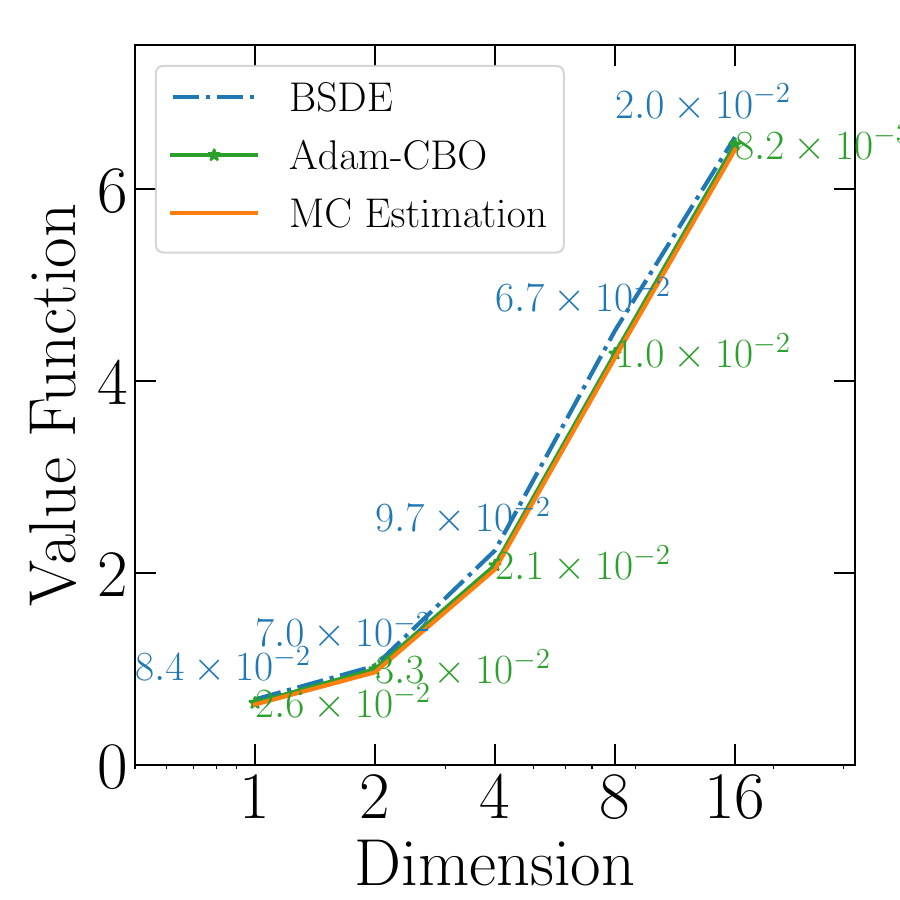}
    \caption{ The value function \( u(t=0, \bx=(0, \ldots, 0)) \) evaluated using BSDE method, M-CBO method, Adam-CBO method (\emph{our method}), and MC estimation (\emph{reference}) for problems in  $1, 2, 4, 8,$ and $16$ dimensions. (a) The terminal cost function \( g(\bx) = \ln \frac{1 + \| \bx \|^2}{2} \). (b) The terminal cost function \( g(\bx) = \ln \frac{1 + (\| \bx \|^2 - 1)^2}{2} \).}
    \label{fig:value_vs_dim}
\end{center}
\vskip -0.2in
\end{figure}

The value function $u(t=0,\bx=(0, \ldots, 0))$  for two different terminal costs - a convex cost: \( g(\bx) = \ln \frac{1 + \|\bx\|^2}{2} \) and a double-well terminal cost: \( g(\bx) = \ln \frac{1 + (\|\bx\|^2 - 1)^2}{2} \) is illustrated in Figure \ref{fig:value_vs_dim}  across varying dimensions. The value function from MC estimation is worked as a reference. The value function of Adam-CBO and M-CBO methods is computed from the expectation of $5000$ controlled dynamics. The value function of the BSDE method is a direct output of the neural network. 

In the convex terminal cost, we can see that both the M-CBO method and the Adam-CBO method outperform the BSDE method in a low-dimensional setting. As the dimensionality increases, Adam-CBO continues to outperform the BSDE method, demonstrating its scalability. Consequently, in the remaining examples, we focus exclusively on the Adam-CBO method because of its superior performance in high dimensions. 

\begin{figure*}[ht]
    \centering
    \includegraphics[width=0.8\linewidth]{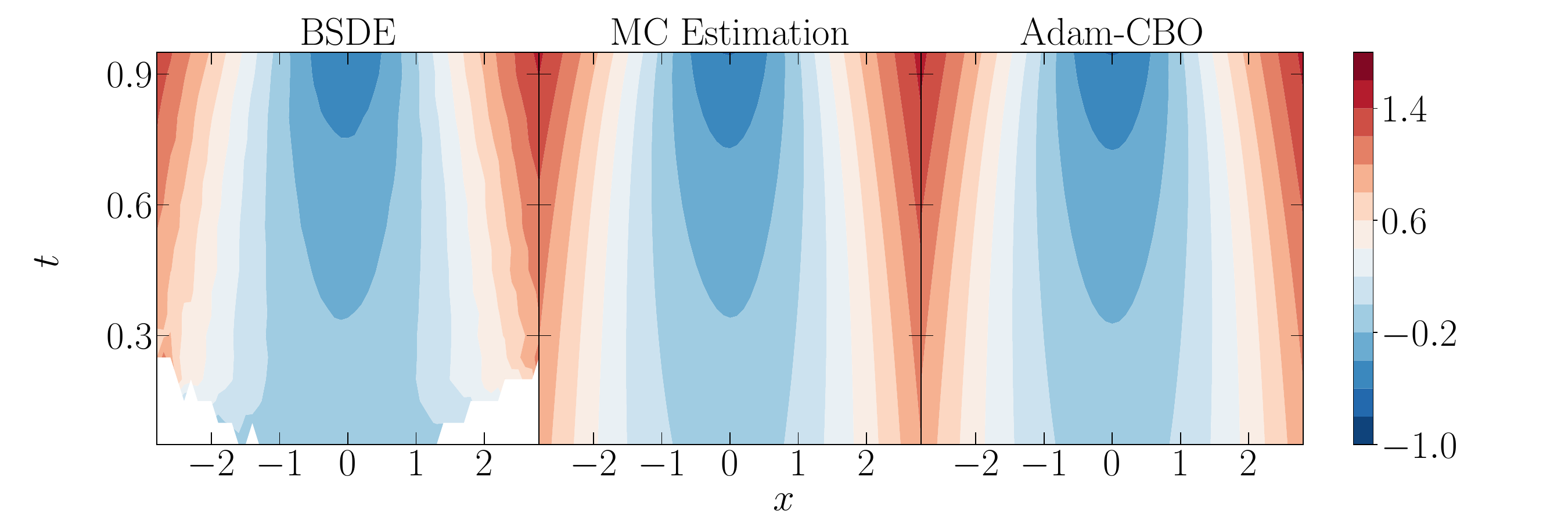}
    \caption{The value function $u(t,x)$ in the one-dimensional case, computed using BSDE method, MC Estimation (\emph{reference}), and Adam-CBO (\emph{our method}), with terminal cost \( g(\bx) = \ln \frac{1 + \| \bx \|^2}{2} \). }
    \label{fig:value_case1}
\end{figure*}

\begin{figure*}[ht]
    \centering
    \includegraphics[width=0.8\linewidth]{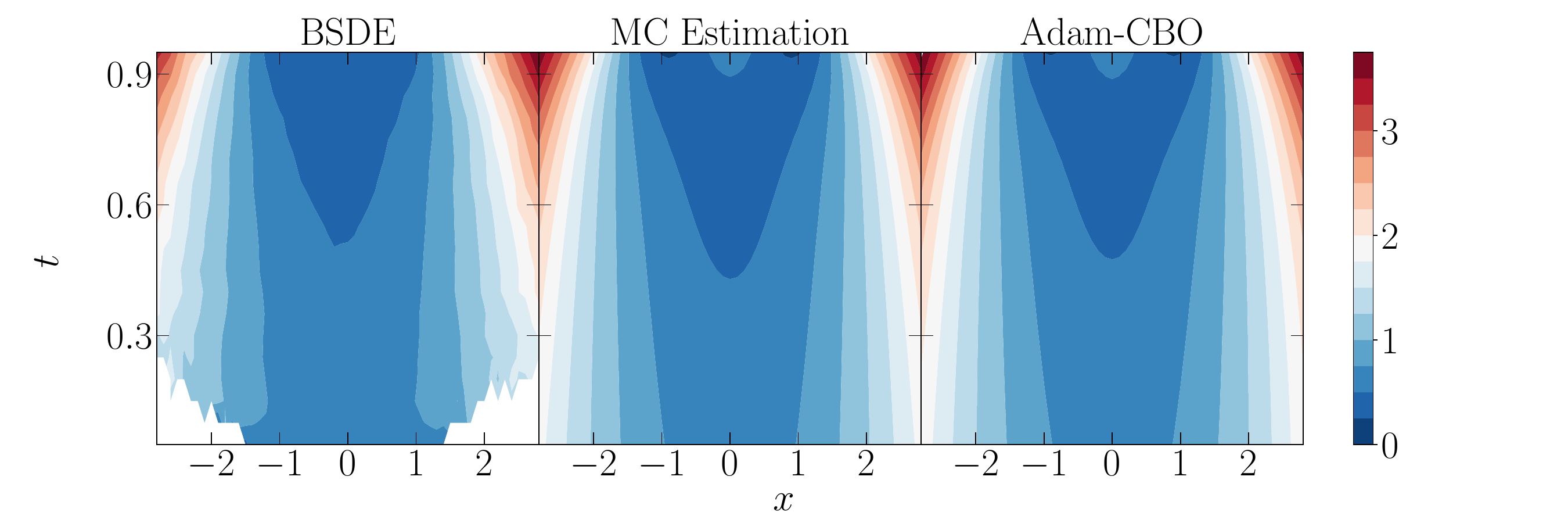}
    \caption{The value function $u(t,x)$ in one-dimensional case, computed using BSDE method, MC Estimation (\emph{reference}), and Adam-CBO (\emph{our method}), with terminal cost \( g(\bx) = \ln \frac{1 + (\| \bx \|^2 - 1)^2}{2} \). }
    \label{fig:value_case2}
\end{figure*}

In the case of the double-well terminal cost, which is nonconvex, our method shows significantly improved accuracy over the BSDE approach. This enhancement can be attributed to several factors. First, CBO-based methods have a higher likelihood of converging to global minima in non-convex settings. Secondly, although both methodologies utilize a discretization of the $20$ time steps during the training phase, our approach facilitates additional refinement of the time steps when assessing the cost function, thereby improving precision. In contrast, the structure of the neural network of the BSDE method is inherently tied to the chosen discretization, necessitating the same time step for both training and evaluation, thus limiting flexibility.
We also visualize the function $u(t,x)$ in the one-dimensional case for both types of terminal costs in Figure \ref{fig:value_case1} and Figure \ref{fig:value_case2}. It is evident that our method aligns more closely with the exact solution than the BSDE-based method.

\begin{figure}[ht]
    \centering
    \includegraphics[width=0.5\linewidth]{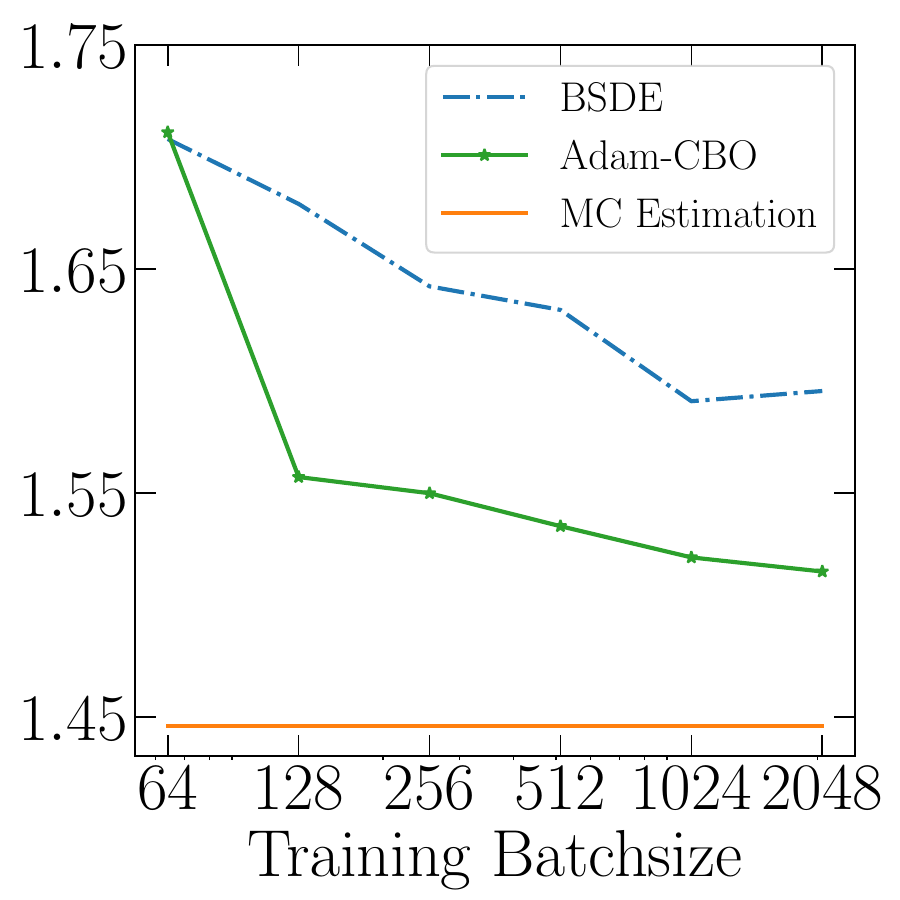}
    \caption{The value function $u(t=0,0,0,0,0)$ of 4D LQC problem, computed using BSDE method, MC Estimation (\emph{reference}), and Adam-CBO (\emph{our method}), 
 with terminal cost \( g(\bx) = \ln \frac{1 + (\|\bx\|^2 - 5)^2}{2} \), evaluated under varying sample sizes per step.".}
    \label{fig:batchsize-value}
\end{figure}

We further investigate the influence of batch size (the number of control processes to compute the cost function) on the problem. We consider a $4$ dimensional problem with a nonconvex terminal cost given by \( g(\bx) = \ln \frac{1 + (\|\bx\|^2 - 5)^2}{2} \). Figure \ref{fig:batchsize-value} illustrates the value function $u(t=0,0,0,0,0)$ evaluated under varying batch sizes during training and compared with a precise estimation of the MC that uses a sufficiently large sample size. The first insight is that the accuracy of training is sensitive to batch size in the training process, which inspired us to develop an improved sampling method to enhance the efficiency of the sampling process in the future. Additionally, our method consistently demonstrates greater accuracy than the BSDE-based approach, confirming its robustness in higher-dimensional and nonconvex settings.

\subsection*{Ginzburg-Landau Model}
We also consider the problem of controlling superconductors in an external electromagnetic field, modeled using the stochastic Ginzburg-Landau theory. 
The dynamics are given by 
\begin{equation}
    \intd \bx_t = \mathbf{b}(\bx_t,\alpha_t) \intd t + \sqrt{2} \intd \mathbf W_t,
\end{equation}
where the drift term is defined as
$$
b(\bx,a) = -\nabla_\bx U(\bx) + 2 \alpha \boldsymbol\omega.
$$
$U$ here is the Ginzburg-Landau free energy, while $\boldsymbol{\omega} \in \mathbb R^d$ specifies the spatial domain of the external field. For further implementation specifics, see Appendix \ref{sec:simulation_datail}.

Since this problem lacks an exact analytical solution, we assess the performance of our trained control $\alpha(t,\bx;\theta)$ by comparing it to the theoretically optimal control 
$- \boldsymbol{\omega} \cdot \nabla_{\bx} u(t,\bx,)$, where $u(\bx,t)$ is the value function. Notably, this value function is different from the last case with an analytical solution; it was computed by taking the expectation of running controlled dynamics and its gradient is computed by taking the finite difference of two starting states. Therefore, this comparison is not intended as a true error metric. Instead, it serves to evaluate the consistency between our trained control and the theoretically optimal control, which many value-based methods use to define the loss.

\begin{figure}[ht]
    \centering
    \includegraphics[width=0.5\linewidth]{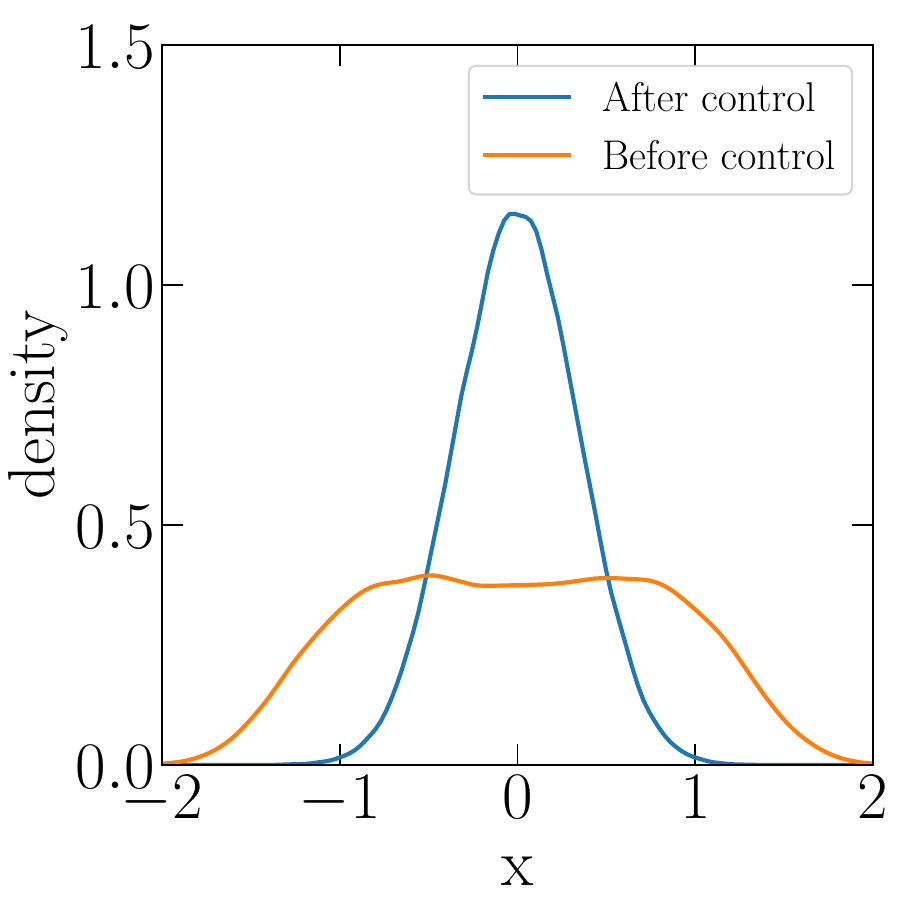}
    \caption{Distribution of $x_1$ before and after control in the 1D Ginzburg-Landau model.}
    \label{fig:1D_2_density}
\end{figure}

\begin{figure}[ht]
    \centering
    \includegraphics[width=\linewidth]{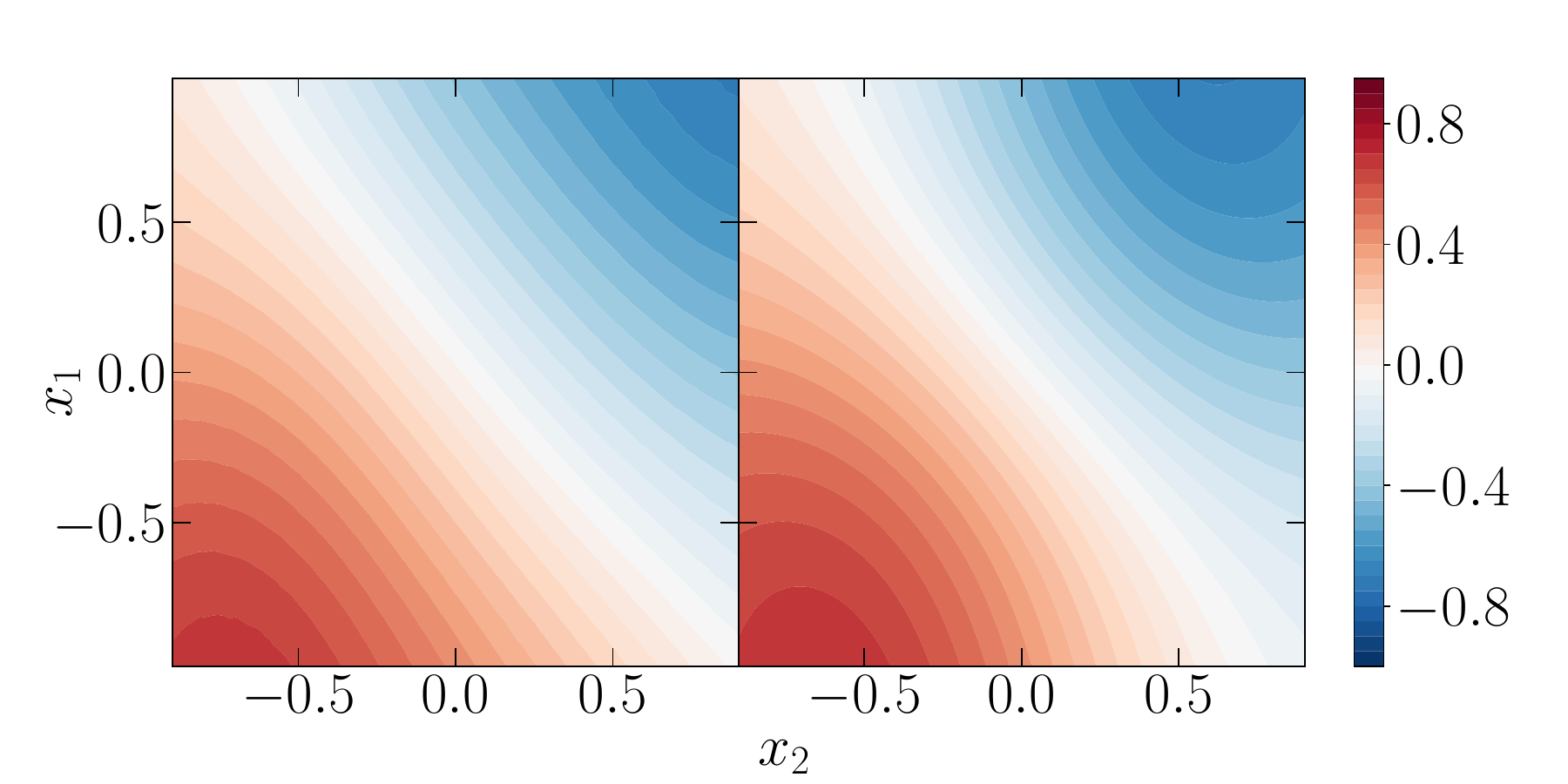}
    \caption{The left figure shows the $\alpha(0.5,\bx;\theta)$ and the right figure shows the $-\boldsymbol{\omega}\cdot\nabla_\bx u(0.5,\bx)$ computed by our method for the 2D Ginzburg-Landau model.}
    \label{fig:compare_control_GL}
\end{figure}

\begin{figure}[ht]
    \centering
    \includegraphics[width=0.45\linewidth]{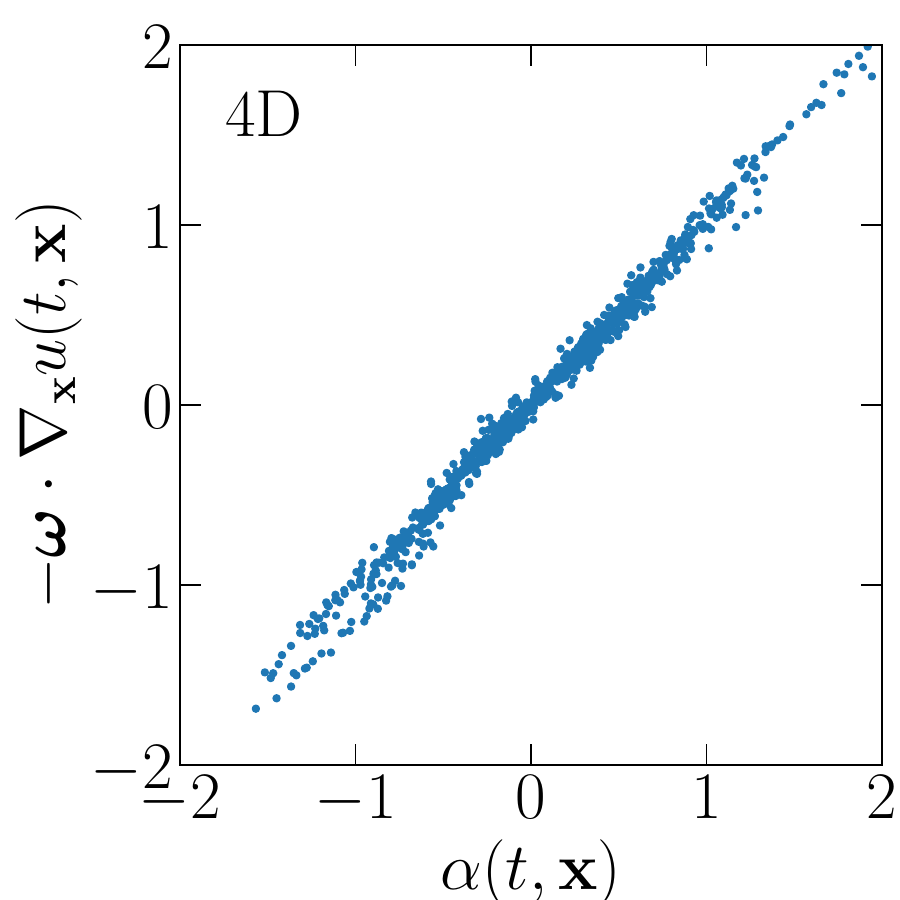}
    \includegraphics[width=0.45\linewidth]{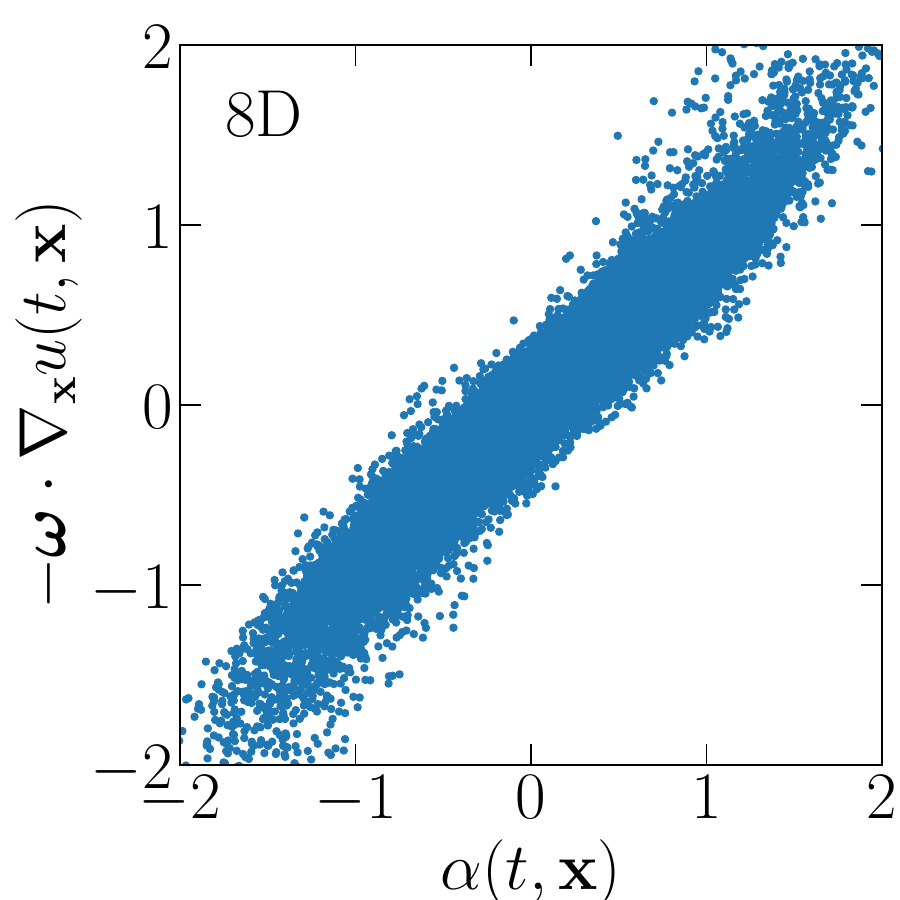}
    \includegraphics[width=0.45\linewidth]{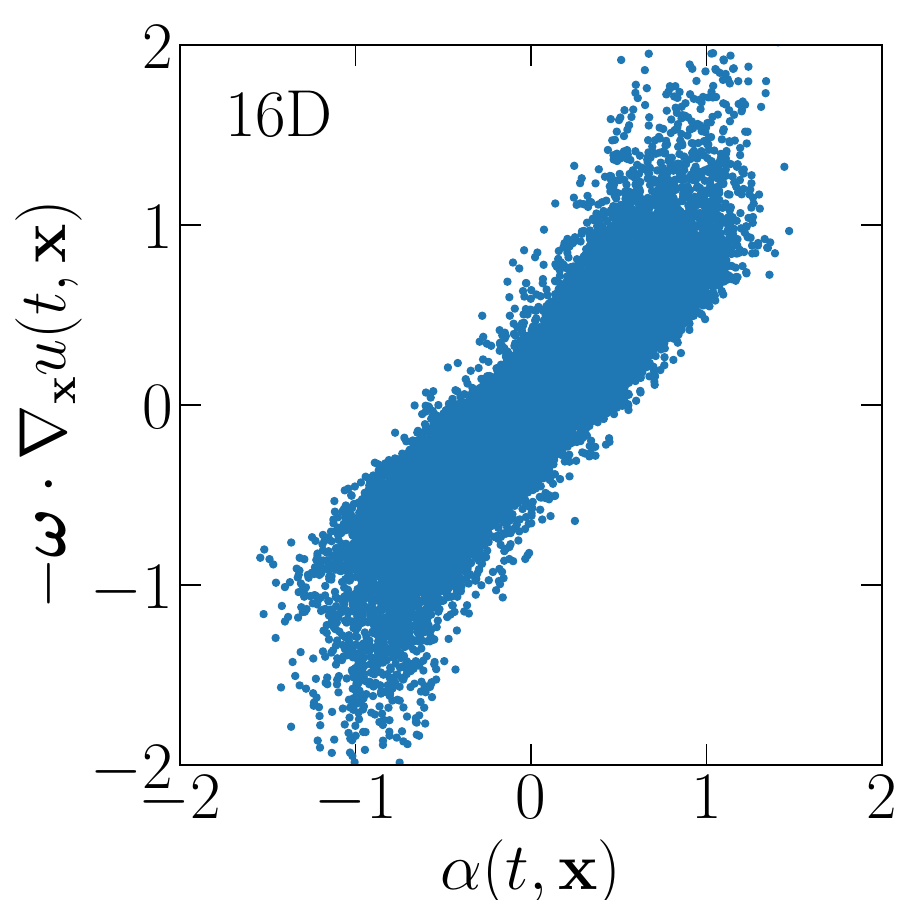}
    \includegraphics[width=0.45\linewidth]{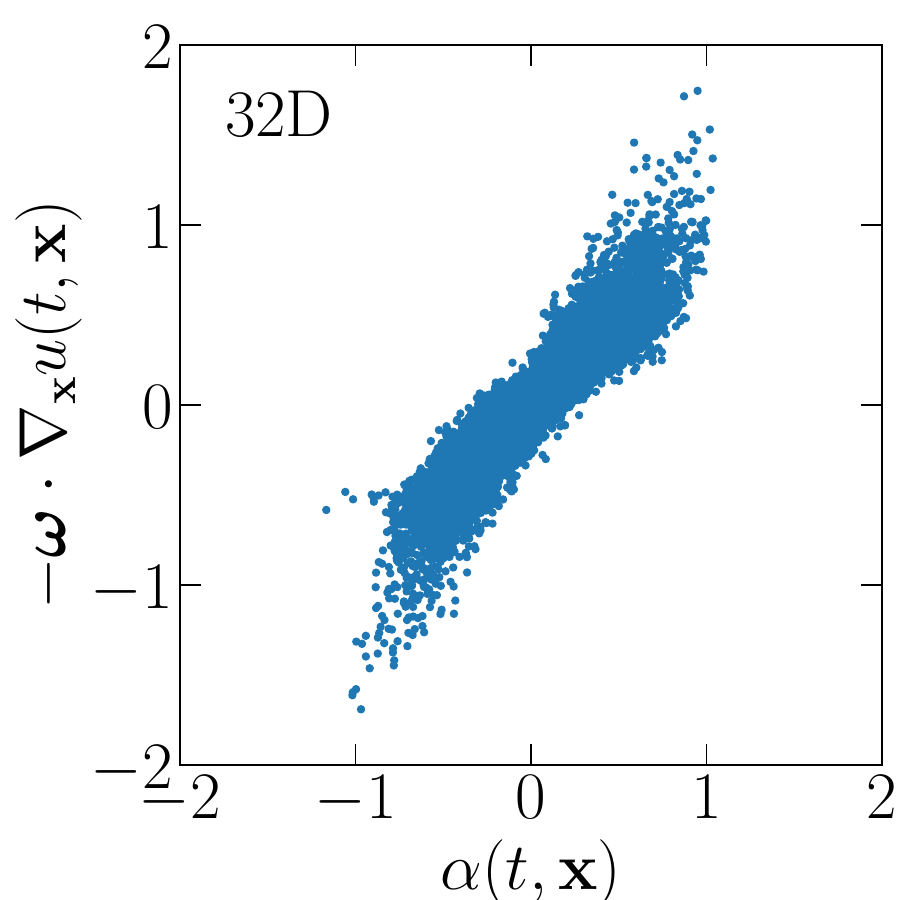}
    \caption{Ginzburg-Landau model with $ d=4,8,16,32$ respectively. The $x$-axis shows the $\alpha(0.5,\bx)$ and $y$-axis  shows the $-\boldsymbol{\omega}\cdot\nabla_\bx u(0.5,\bx)$ }
    \label{fig:compare_control_GL_high}
\end{figure}

We start with a simple case with $d= 2$, $\mu=10$, $\lambda=0.2$. We compare the distribution of $x_1$ before and after the control in Figure \ref{fig:1D_2_density}. One can find that before the control the particles will stay in a stable state $-1,1$, while after control the particles will stay near $0$. 
Additionally, a comparative analysis between $\alpha(t, \bx)$ and $- \boldsymbol{\omega} \cdot \nabla_{\bx} u(t,\bx)$ is conducted, as illustrated in Figure \ref{fig:compare_control_GL}. The results demonstrate the consistency between these two functions.
We also test our method on $d=4,8,16,32$. The comparison between $\alpha(t, \bx;\theta)$ and $- \boldsymbol{\omega} \cdot \nabla_{\bx} u(t,\bx)$ is shown in Figure  \ref{fig:compare_control_GL_high}. Here $(t,\bx)$ is randomly sampled from 1000 control dynamics.

\subsection*{Systemic Risk Mean Field Control}
In practical applications, there are scenarios where numerous indistinguishable agents, such as multiple traders engaged in buying and selling stocks within financial markets, create a complex, multi-dimensional problem. However, when these traders share similar risk preferences, analyzing the behavior of a single representative trader can suffice to understand the dynamics of the entire group. For example, for a problem with $n$ agents, the control can be modeled as $\ctrl^i(\bx_i,\mu;\theta)$, where $\mu$ is the empirical measure of the $\{\bx_i\}_{i=1}^n$ and $\theta$ are parameters in the neural network. For further details on the network construction used in this setup, refer to Appendix \ref{sec:network}.


Consider the systemic risk mean field control problem, detailed in Appendix \ref{sec:simulation_datail}. The control policy is initially trained using a delta distribution centered on $x_0$ and $n = 100$ and then tested against different values of $n=50,100,200,400,800$. Furthermore, the value function is evaluated by taking the expectation of controlled dynamics starting from different initial distributions $\mu_0$, including Gaussian random variable $x_0= \mathcal{N}(0,0.1)$,  mixture of two Gaussian random variables $x_0 = p(-k+\theta y) +(1-P)(k+\theta z) $ with $P$ a Bernoulli random variable with parameter $\frac{1}{2}$, $k=\frac{\sqrt{3}}{10}$, $\theta=0.1$, $y,z\sim N(0,1)$ and mixture of three Gaussian random variables: $x_0 = [-k_{\lfloor 3U \rfloor =0} + k_{\lfloor 3U \rfloor =1}] + \theta y$ with $k = 0.3$, $\theta=0.07$, $y\sim N(0,1)$. The corresponding value functions for each scenario are shown in Table \ref{tab:value_function_gaussian}, \ref{tab:value_function_two_gaussian}, \ref{tab:value_function_three_gaussian}, respectively.  Our method demonstrates robust generalization across these diverse conditions, in contrast to value-function-based approaches where the control strategy is tied to the specific value function. Since value functions are highly sensitive to initial conditions, traditional methods require retraining for each new initial scenario, limiting their ability to generalize effectively.
\begin{table}[t]
\caption{The value function $u(t,\mu)$ evaluated by Adam-CBO method with $\mu=\mathcal{N}(0,1)$ in the mean filed control problem.}
\label{tab:value_function_gaussian}
\vskip 0.15in
\begin{center}
\begin{small}
\begin{sc}
\begin{tabular}{lccccr}
\toprule
time &  $n=50$ &  $n=100$ & $n=200$ & $n=400$ & Exact\\  
\midrule
0.0 & 0.607 & 0.614 & 0.618 & 0.619  & 0.616 \\ 
0.1 & 0.553 & 0.559 & 0.563 & 0.564  & 0.561 \\ 
0.2 & 0.498 & 0.504 & 0.507 & 0.508  & 0.506 \\ 
0.3 & 0.442 & 0.447 & 0.450 & 0.451  & 0.449 \\ 
0.4 & 0.384 & 0.388 & 0.391 & 0.391  & 0.390 \\ 
0.5 & 0.323 & 0.326 & 0.329 & 0.329  & 0.329 \\ 
0.6 & 0.258 & 0.260 & 0.262 & 0.263  & 0.262 \\ 
0.7 & 0.187 & 0.188 & 0.190 & 0.190  & 0.190 \\ 
0.8 & 0.106 & 0.107 & 0.108 & 0.108  & 0.108 \\ 
0.9 & 0.010 & 0.010 & 0.010 & 0.010  & 0.010 \\ 
\bottomrule
\end{tabular}
\end{sc}
\end{small}
\end{center}
\vskip -0.1in
\end{table}

\begin{table}[t]
\caption{The value function $u(t,\mu)$  with $\mu$ being a mixture of two Gaussian random variables in the mean filed control problem.}
\label{tab:value_function_two_gaussian}
\vskip 0.15in
\begin{center}
\begin{small}
\begin{sc}
\begin{tabular}{lccccr}
\toprule
time &  $n=50$ &  $n=100$ & $n=200$ & $n=400$  & Exact\\  
\midrule
0.0 & 0.621 & 0.628 & 0.633 & 0.634 & 0.630 \\ 
0.1 & 0.567 & 0.574 & 0.578 & 0.579 & 0.576 \\ 
0.2 & 0.513 & 0.518 & 0.522 & 0.523 & 0.521 \\ 
0.3 & 0.457 & 0.462 & 0.465 & 0.466 & 0.465 \\ 
0.4 & 0.399 & 0.404 & 0.407 & 0.408 & 0.407 \\ 
0.5 & 0.339 & 0.343 & 0.346 & 0.346 & 0.346 \\ 
0.6 & 0.276 & 0.279 & 0.281 & 0.281 & 0.281 \\ 
0.7 & 0.207 & 0.209 & 0.211 & 0.211 & 0.211 \\ 
0.8 & 0.129 & 0.131 & 0.132 & 0.132 & 0.132 \\ 
0.9 & 0.039 & 0.040 & 0.040 & 0.040 & 0.040 \\ 
\bottomrule
\end{tabular}
\end{sc}
\end{small}
\end{center}
\vskip -0.1in
\end{table}

\begin{table}[t]
    \caption{The value function $u(t,\mu)$ with $\mu$ being a mixture of three Gaussian random variables in the mean filed control problem.}
    \label{tab:value_function_three_gaussian}
\vskip 0.15in
\begin{center}
\begin{small}
\begin{sc}
\begin{tabular}{lccccr}
\toprule
time &  $n=50$ &  $n=100$ & $n=200$ & $n=400$ & Exact\\  
\midrule
0.0 & 0.633 & 0.640 & 0.645 & 0.646  & 0.642 \\ 
0.1 & 0.579 & 0.586 & 0.590 & 0.591  & 0.588 \\ 
0.2 & 0.524 & 0.531 & 0.535 & 0.536  & 0.534 \\ 
0.3 & 0.469 & 0.475 & 0.478 & 0.47   & 0.478 \\ 
0.4 & 0.412 & 0.417 & 0.421 & 0.421  & 0.420 \\ 
0.5 & 0.353 & 0.357 & 0.360 & 0.361  & 0.360 \\ 
0.6 & 0.290 & 0.294 & 0.297 & 0.297  & 0.297 \\ 
0.7 & 0.223 & 0.226 & 0.228 & 0.228  & 0.228 \\ 
0.8 & 0.148 & 0.151 & 0.152 & 0.152  & 0.152 \\ 
0.9 & 0.063 & 0.064 & 0.064 & 0.065  & 0.065 \\ 
\bottomrule
\end{tabular}
\end{sc}
\end{small}
\end{center}
\vskip -0.1in
\end{table}

\section{Conclusion}
In this work, we present a framework for solving high-dimensional stochastic optimal control problems. Compared with the existing method, our method is gradient-free, which eliminates the high variance in the Monte Carlo estimation of the policy gradient. Also, our method does not depend on solving the high-dimensional Hamiltonian-Jacobi-Bellman equation or on any mesh discretization in the state and action space. These enable us to get rid of the curse of dimensionality and use this method in high-dimensional problems. Theoretically, we show that, under some assumptions, the M-CBO method can converge to the optimal control. In the future, we are interested in applying our method to mean-field game problems and control problems with partial information and constraints~\cite{Ganapathi2024Confidence, Hong2024primal, Qiao2024Near, sun2024Constrained, Wang2024Probabilistic}.






\section*{Impact Statement}


This paper presents work whose goal is to advance the field of 
Machine Learning. There are many potential societal consequences 
of our work, none which we feel must be specifically highlighted here.

\appendix
\section{Well-posedness of the M-CBO method}
\label{app:well-posedness_m_cbo_fix_N}
In this section, we prove that the dynamics of the M-CBO method are well-posed. For an arbitrary but fixed $N$, we begin by studying the existence of a unique process $\l (\boldsymbol{\Theta}_t^{(N)},\boldsymbol{\Omega}_t^{(N)} \r )= \l (\Theta^{(1,N)},\cdots,\Theta^{(N,N)},\Omega^{(1,N)},\cdots, \Omega^{(N,N)} \r )$ that satisfies the M-CBO scheme \eqref{equ:CBO_SDE}
\begin{equation}
\label{equ:M_CBO_fixed_N}
\begin{aligned}
        \intd \boldsymbol{\Theta}\pNt &=    \mathbf{F}_{N,\Theta}\lg\bTheta\pNt,\bOmega\pNt\r )\intd t + \sigma(t)  \intd \mathbf{W}\pNtt,
        \\
     \intd 
 \boldsymbol{\Omega}\pNt &= 
 \mathbf{F}_{N,\Omega}\lg\bTheta^{(N)}_t,\bOmega\pNt \r )\intd t +  \sqrt{m} \sigma(t)   \intd \mathbf{W}\pNot,
\end{aligned}
\end{equation}
where $\mathbf W\pNtt, \mathbf W\pNot$ is the standard Wiener process in $\mathbb{R}^{N D}$, and 
\begin{equation*}
    \begin{aligned}
        \mathbf{F}_{N,\Theta}(\bTheta,\bOmega) &= \l ( F^1_{N,\Theta} (\bTheta,\bOmega),\cdots F^N_{N,\Theta}(\bTheta,\bOmega)\r )  \in \mathbb{R}^{ND},
        \\ \text{ }
        \mathbf{F}_{N,\Omega}(\bTheta,\bOmega) &= \l ( F^1_{N,\Omega} (\bTheta,\bOmega),\cdots F^N_{N,\Omega}(\bTheta,\bOmega)\r )  \in \mathbb{R}^{ND},
        \\
        F^i_{N,\Theta}(\bTheta,\bOmega) 
         &=   \Omega^i  
        - \gamma_1\frac{\sum_{j\neq i}(\Theta^i-\Theta^j)w_\beta(\Theta^j)}{\sum_j w_\beta(\Theta^j)},\\
            F^i_{N,\Omega}(\bTheta,\bOmega) &=  - m\frac{\sum_{j\neq i}(\Theta^i-\Theta^j)w_\beta(\Theta^j)}{\sum_j w_\beta(\Theta^j)}  + \gamma_2\Omega^i.
    \end{aligned}
\end{equation*}

Under the Assumption \ref{assum:Lip_J}, we can easily deduce that $F^i_{N,\Theta}$ and $F^i_{N,\Omega}$ are locally Lipschitz continuous and have linear growth. Consequently, $\l ( {\bF_{N,\Theta}},{\bF_{N,\Omega}}\r ) $ is locally Lipschitz continuous and has linear growth.  More precisely, we have the following lemma. 
\begin{lemma}\label{lem:linear}
    Let $N\in \mathbb N$, $\beta,R>0$ be arbitrary. Then for any $(\bTheta,\bOmega),(\hat\bTheta,\hat\bOmega)\in \mbRd\times \mbRd$ with $\|\bTheta\|+\|\bOmega\|,\|\hat\bTheta\|+\|\hat\bOmega\|\leq R$ and all $i=1,\cdots,N$, it holds
    \begin{equation*}
        \begin{aligned}
            \|F^i_{N,\Theta}(\bTheta,\bOmega)-F^i_{N,\Theta}(\hat\bTheta,\hat\bOmega)\|
            & \leq
            \gamma_1\|\Theta^i-\hat\Theta^i\|+
            \|\Omega^i-\hat\Omega^i\|  
            +\gamma_1\l ( 1 + 2\frac{c_R}{N}\sqrt{N \|\hat{\Theta}^i\|^2 +\|\hat\bTheta\|^2}\r ) \|\bTheta-\hat\bTheta\|,
            \\
            \|F^i_{N,\Omega}(\bTheta,\bOmega)-F^i_{N,\Omega}(\hat\bTheta,\hat\bOmega)\| & \leq \|\Theta^i-\hat\Theta^i\|+
            \gamma_2\|\Omega^i-\hat\Omega^i\|  
            +m\l ( 1 + 2\frac{c_R}{N}\sqrt{N \|\hat{\Theta}^i\|^2 +\|\hat\bTheta\|^2}\r ) \|\bTheta-\hat\bTheta\|,\\
            \|F^i_{N,\Theta}(\bTheta,\bOmega)\| &\leq \gamma_1\|\Theta^i\|+\|\Omega^i\|+\gamma_1\|\bTheta\|,\\
            \|F^i_{N,\Omega}(\bTheta,\bOmega)\| &\leq m\|\Theta^i\|+\gamma_2\|\Omega^i\|+m\|\bTheta\|,
        \end{aligned}
    \end{equation*}
    where $c_R = \alpha \|\nabla \mathcal J \|_{L^\infty(B_{\theta,R}(0))} \exp\l (\beta\|\mathcal J - \underline{\mathcal J }\|_{L^\infty(B_{\theta,R}(0))} \rg$.
\end{lemma}
\begin{proof}
    From Lemma 2.1 \cite{carrillo2018analytical}, we have
    \begin{equation*}
    \begin{aligned}
        \l \|\frac{\sum_{j\neq i}(\Theta^i-\Theta^j)w_\beta(\Theta^j)}{\sum_j w_\beta(\Theta^j)}-\frac{\sum_{j\neq i}(\hat\Theta^i-\hat\Theta^j)w_\beta(\hat\Theta^j)}{\sum_j w_\beta(\hat\Theta^j)}\r  \| 
        &\leq \|\Theta^i-\hat\Theta^i\|   +\l ( 1 + 2\frac{c_R}{N}\sqrt{N \|\hat{\Theta}^i\|^2 +\|\hat\bTheta|^2}\r ) \|\bTheta-\hat\bTheta\|,\\
        \l \|\frac{\sum_{j\neq i}(\Theta^i-\Theta^j)w_\beta(\Theta^j)}{\sum_j w_\beta(\Theta^j)}\r  \| 
        &\leq \|\Theta^i\| +\|\bTheta\|.
    \end{aligned}
    \end{equation*}
    By the triangle inequality, the required estimation is proved.
\end{proof}
Based on Lemma \ref{lem:linear}, we may invoke standard existence results of strong solutions for Equation \eqref{equ:CBO_SDE}.
\begin{proof}[Proof of Theorem \ref{thm:well-posedness_m_cbo_fix_N}]
    We make use of the standard result on the existence of a unique strong solution here. 
    To this end, we show the existence $b_N>0$, such that 
    \begin{equation*}
        \bTheta \cdot \mathbf{F}_{N,\Theta}(\bTheta,\bOmega) 
        + \bOmega\cdot \mathbf{F}_{N,\Omega}(\bTheta,\bOmega) + N (m+1)D\sigma(t)^2 \leq b_N(\|\bTheta\|^2 +\|\bOmega\|^2 + 1 ).
    \end{equation*}
    Notice that 
    \begin{equation*}
    \begin{aligned}
        - \Theta^i\frac{\sum_{j\neq i}(\Theta^i-\Theta^j)w_\beta(\Theta^j)}{\sum_j w_\beta(\Theta^j)} &\leq -\|\Theta^i\|^2 + \|\Theta^i\|\|\bTheta\|, \\
        - \Omega^i\frac{\sum_{j\neq i}(\Theta^i-\Theta^j)w_\beta(\Theta^j)}{\sum_j w_\beta(\Theta^j)} &\leq \|\Omega^i\|\|\Theta^i\| + \|\Omega^i\|\|\bTheta\|, \\
    \end{aligned}
    \end{equation*}
we have the following inequalities
\begin{equation*}
\begin{aligned}
    \Theta^i  F^i_{N,\Theta}(\bTheta,\bOmega)  & = \Theta^i\Omega^i  
        - \gamma_1\Theta^i\frac{\sum_{j\neq i}(\Theta^i-\Theta^j)w_\beta(\Theta^j)}{\sum_j w_\beta(\Theta^j)} \\
    &\leq  \frac{1}{2}\|\Theta^i\|^2 + \frac{1}{2}\|\Omega^i\|^2 -\gamma_1 \|\Theta^i\|^2 + \gamma_1\|\Theta^i\|\|\bTheta\|\\
    & \leq \frac{1}{2}\|\Theta^i\|^2 + \frac{1}{2}\|\Omega^i\|^2 + \gamma_1\|\Theta^i\|\|\bTheta\|\\
    &\leq \l ( \frac{1}{2}+\frac{\gamma_1}{2}\r ) \|\Theta^i\|^2 + \frac{\gamma_1}{2}\|\bTheta\|^2  + \frac{1}{2}\|\Omega^i\|^2,
\end{aligned}
\end{equation*}
and 
\begin{equation*}
\begin{aligned}
    \Omega^i  F^i_{N,\Omega}(\bTheta,\bOmega)  & =  - m\Omega^i \frac{\sum_{j\neq i}(\Theta^i-\Theta^j)w_\beta(\Theta^j)}{\sum_j w_\beta(\Theta^j)}  + \gamma_2\|\Omega^i\|^2   \\
    & \leq m \|\Omega^i\|\|\Theta^i\|  + m \|\Omega^i\|\|\bTheta\| + \gamma_2\|\Omega^i\|^2 \\
    & \leq (m+\gamma_2) \|\Omega^i\|^2  + \frac{m}{2} \|\Theta^i\|^2  + \frac{m}{2}\|\bTheta\|^2.
\end{aligned}
\end{equation*}
Therefore, we conclude that 
\begin{equation}
    \begin{aligned}
     & \bTheta \cdot \mathbf{F}_{N,\Theta}(\bTheta,\bOmega) 
        + \bOmega\cdot \mathbf{F}_{N,\Omega}(\bTheta,\bOmega) + (m+1)D\sigma(t)^2     \\
    &\leq
        N(m+1)D\sigma(t)^2+ 
        \sum_{i=1}^N \l (  
        \Theta^i \cdot F^i_{N,\Theta}(\bTheta,\bOmega) + \Omega^i \cdot F^i_{N,\Omega}(\bTheta,\bOmega)\r ) \\
    & \leq  N(m+1)D\sigma(t)^2+ 
    \sum_{i=1}^N \l ( \l ( \frac{1}{2}+\frac{\gamma_1}{2}\r ) \|\Theta^i\|^2 + \frac{\gamma_1}{2}\|\bTheta\|^2  + \frac{1}{2}\|\Omega^i\|^2 + (m+\gamma_2) \|\Omega^i\|^2  + \frac{m}{2} \|\Theta^i\|^2  + \frac{m}{2}\|\bTheta\|^2 \r ) \\
    & \leq  N(m+1)D\sigma(t)^2 + \l (  \frac{1+\gamma_1+m}{2} + \frac{\gamma_1+ N }{2}  \r ) \|\bTheta\|^2   + \l ( m+\gamma_2+\frac{1}{2}\r ) \|\bOmega\|^2\\
    & \leq b_N(\|\bTheta\|^2+\|\bOmega\|^2+1),
    \end{aligned}
\end{equation}
where $b_N = \max\{ N(m+1)D\bar\sigma^2,  \frac{1+\gamma_1+m}{2} + \frac{\gamma_1+ N }{2}  , m+\gamma_2+\frac{1}{2}  \} >0$. Then we apply Theorem 3.1 in \cite{durrett2018stochastic} to finish the existence and uniqueness proof.
\end{proof}
\section{Well-posedness of the Mean Field Equations}
\label{app:well-posedness_mean_field}
\begin{definition}
    We say $\rho_t \in \MC{[0,T],\PP{4}}$ is a weak solution to the Fokker-Planck equation \eqref{equ:FP_M_CBO}  with initial condition $\rho_0$, if $\forall \phi\in \mcC2d$, we have 
    \begin{equation}
    \begin{aligned}
        \frac{\intd}{\intd t } \int \phi(\theta,\omega) \intd \rho_t =& 
        \int \langle\omega-\gamma_1(\theta-\CM [\mu_t]),\nabla_\theta \phi\rangle \intd 
 \rho_t 
        \\
        &- 
        \int \langle m(\theta - \CM [\mu_t]) +\gamma_2 \omega,\nabla_\omega \phi\rangle\intd  \rho_t \\
        & +\frac{m \sigma(t)^2}{2} \int \Delta_\omega \phi\intd \rho_t +\frac{\sigma(t)^2}{2} \int \Delta_\theta \phi\intd \rho_t,
    \end{aligned}
    \end{equation}
    and $\lim_{t\to \infty}\rho_t = \rho_0$ in a pointwise sense.
\end{definition}
To prove Theorem \ref{thm:well_possness_mean_field}, we start with the following lemma. 
\begin{lemma}\label{lem:wasserstwin}
    If $\mathcal J $ satisfies Assumption  \ref{assum:Lip_J} and $\rho,\hat\rho \in \PP{2}$ with  
    \[
\int \|\theta\|^4 +\|\omega\|^4 \intd \rho, \int \|\hat\theta\|^4 +\|\hat\omega\|^4\intd \hat \rho \leq K,
    \]
    then the following stability estimate holds
    \[
    |\CM [\mu] - \CM[\hat\mu]| \leq c_0 W_2(\rho,\hat\rho)
    \]
    for a constant $c_0>0$ depending on $\beta,L_J$ and $K$, where $\mu(\theta) = \int_{\mbRd} \rho (\theta,\intd \omega),\hat\mu(\hat \theta) = \int_{\mbRd} \hat \rho (\hat \theta,\intd \hat \omega)$.
\end{lemma}
\begin{proof}
    By Lemma 3.2 in \cite{carrillo2018analytical}, we have 
    $|\CM[\mu] - \CM[\hat\mu]| 
    \leq c_0 W_2(\mu,\hat\mu)
    =c_0  \inf \mathbb  E_{(\mu,\hat \mu)}[\|\theta-\hat \theta\|^2]$. Therefore, we have  
    $\inf E_{(\mu,\hat \mu)}[\|\theta-\hat \theta\|^2] \leq  \inf \mathbb  E_{(\rho,\hat\rho)}[\|\theta-\hat \theta\|^2] +  \inf \mathbb  E_{(\rho,\hat\rho)}[\|\omega-\hat \omega\|^2] \leq \inf \mathbb E_{(\rho,\hat\rho)}[\|\theta-\hat \theta\|^2 + \|\omega-\hat \omega\|^2] = W_2(\rho,\hat\rho)$. We prove the boundedness here.
\end{proof}
To prove the existence and uniqueness, we recall the Leray-Schauder fixed point theorem (Theorem 11.3 in \cite{gilbarg2001elliptic}).
\begin{theorem}\label{thm:compact}
  Let $T$ be a compact mapping of a Banach space $\mathcal B$ into itself, and suppose there exists a constant $M$ such that $\|x\|_{\mathcal B}\leq M$ for all $x\in \mathcal B$ and $\eta\in (0,1)$ satisfying $x=\eta T x $. Then T has a fixed point.
\end{theorem}
\begin{proof}[Proof of Theorem \ref{thm:well_possness_mean_field}]

\textbf{Step 1 (Construction of map $T$)}

Let us fix $u_t\in \MC{[0,T]}$. By Theorem 6.2.2 in \cite{Arnold1976stochastic}, there is a unique solution to 
\begin{equation*}\label{equ:map_sde}
\begin{aligned}
    \intd \theta_t  &=   \omega_t \intd t - \gamma_1(\theta_t - u_t) \intd t  +  \sigma(t)\intd W_{\theta,t}, \\
    \intd \omega_t &=  -m (\theta_t - u_t) \intd t  -  \gamma_2 \omega_t   \intd t +  \sqrt{m}  \sigma(t)\intd W_{\omega,t},
\end{aligned}
\end{equation*}
where $(\theta_0,\omega_0)\sim \rho_0$. We use $\rho_t$ to denote the corresponding law of the unique solution. Using $\rho_t$, one can compute $\mu_t(\theta) = \int \rho_t(\theta,  \intd \omega) $ and $\CM [\mu_t]$, which is uniquely determined by $u_t$ and is in $\MC{[0,T]}$. Thus, one can construct a map from $\MC{[0,T]}$ to $\MC{[0,T]}$, which maps $u_t$ to $\CM [\mu_t]$.

\textbf{Step 2 (Compactness)}
    First, by referencing Chapter 7 in \cite{Arnold1976stochastic}, we obtain the inequality for the solution $\theta_t ,\omega_t$ to Equation \eqref{equ:map_sde}:
    \begin{equation*}
        \mathbb E [\|\theta_t\| + \|\omega_t\|]^4 \leq \l ( 1+\mathbb E [\|\theta_0\| + \|\omega_0\|]^4\r ) \exp(ct),
    \end{equation*}
    where $c>0$. Thus one can deduce $\mathbb E [\|\theta_t\|^4 + \|\omega_t\|^4] \lesssim 1 $ and $\mathbb E [\|\theta_t\|^2 + \|\omega_t\|^2] \lesssim 1$.

    By Lemma \ref{lem:wasserstwin}, we have $\|\mathcal M_\beta (\mu_t)- \mathcal M_\beta (\mu_s)\|\leq c_0 W_2(\rho_t,\rho_s)$. For $W_2(\rho_t,\rho_s)$, it holds that 
    $W_2(\rho_t,\rho_s)\leq E[\|\theta_t-\theta_s\|+\|\omega_t-\omega_s\|]\leq\sqrt{2 E[\|\theta_t-\theta_s\|^2+2\|\omega_t-\omega_s\|^2]}$. Further, we can deduce 
    \begin{equation*}
    \begin{aligned}
                \theta_t-\theta_s &=  \int_s^t \omega_\tau -\gamma_1 (\theta_\tau-u_\tau) \intd \tau + \int_s^t \sigma(\tau)\intd W_{\theta,\tau},\\
                \omega_t-\omega_s &=  \int_s^t  -m (\theta_\tau-u_\tau) -\gamma_2 \omega_\tau \intd \tau +  \sqrt{m} \int_s^t \sigma(\tau)\intd W_{\omega,\tau}.\\
    \end{aligned}
    \end{equation*}
    Thus 
    \begin{equation*}
    \begin{aligned}
        \mathbb E [\|\theta_t-\theta_s \|^2+\|\omega_t-\omega_s \|^2]
        \lesssim & \mathbb E \l [\l \| \int_s^t \omega_\tau \intd \tau \r  \|^2\r  ]
        +
        \mathbb E \l [\l \| \int_s^t 
    (\theta_\tau-u_\tau) \intd \tau \r  \|^2\r  ]\\
    & +  \mathbb E \l [\l \| \int_s^t   \sigma(\tau)\intd W_{\theta,\tau} \r  \|^2\r  ]
  + m\mathbb E \l [\l \| \int_s^t   \sigma(\tau)\intd W_{\omega,\tau} \r  \|^2\r  ].
    \end{aligned}
    \end{equation*}
    Let us proceed to bound the four terms on the right-hand side individually. Consider the first term, where we establish \begin{equation*}
        \begin{aligned}
            \mathbb E \l [\l \| \int_s^t \omega_\tau \intd \tau \r  \|^2\r  ] 
            & \leq 
            \mathbb E \l [\l ( \int_s^t \l \|  \omega_\tau  \r  \| \intd \tau \r ) ^2\r  ] \\
            & \leq |t-s|
             \mathbb E \l [\int_s^t \l \|  \omega_\tau  \r  \|^2  \intd \tau \r  ]  \lesssim |t-s|.
        \end{aligned}
    \end{equation*}
    The first inequality in this sequence is derived from the Cauchy-Schwarz inequality, followed by an application of Jensen's inequality for the second inequality. The final inequality is attributed to the boundedness property of the solution.
Similarly, for the second term, we have
\begin{equation*}
\begin{aligned}
        \mathbb E \l [\l \| \int_s^t 
    (\theta_\tau-u_\tau) \intd \tau \r  \|^2\r  ]  
    & \leq 
    \mathbb E \l [\l ( \int_s^t 
    \l \|  \theta_\tau-u_\tau\r  \| \intd \tau\r ) ^2 \r  ]  \\
    & \leq |t-s| \mathbb E \l [\int_s^t 
    \l \|  \theta_\tau-u_\tau\r  \|^2 \intd \tau \r  ] \\
    & \lesssim  |t-s| \l ( \mathbb E \l [\int_s^t 
    \l \|  \theta_\tau \r  \|^2 \intd \tau  
 +\int_s^t\l \| u_\tau\r  \|^2 \intd \tau \r  ]\r )   \lesssim |t-s|.
\end{aligned}
\end{equation*}
For the third and fourth terms, we use the It\^o Isometry,
\begin{equation*}
    \begin{aligned}
    \mathbb E \l [\l \| \int_s^t   \sigma(\tau)\intd W_{\omega,\tau} \r  \|^2\r  ] 
         = \mathbb E \l [\l \| \int_s^t   \sigma(\tau)\intd W_{\theta,\tau} \r  \|^2\r  ] = \mathbb E \l [ 
 \int_s^t   \sigma(\tau)^2 \intd \tau  \r  ] \leq \bar\sigma^2 |t-s|.
    \end{aligned}
\end{equation*}
Finially, we combine the inequiality to deduce that $\|\CM[\mu_t] -\CM[\mu_s]\|\lesssim |t-s|^{1/2}$, which implies that $\mathcal M_\beta(\mu_t)\in \mathcal C ^{0,1/2}[0,T]$. Thus $T$ is compact.

\textbf{Step 3 (Existence)}
We make use of Theorem \ref{thm:compact}. Take $u_t=\eta T u_t$ for $\eta \in [0,1]$. We now try to prove $\|u_t\|_\infty\leq q$ for some finite $q>0$.
First, one has
\begin{equation*}
    \|u_t\|^2 = \eta^2 \|\mathcal M_\beta (\mu_t)\|^2 \leq \eta^2 \exp(\beta (\overline{J}-\underline{J})) \int \|\theta\|^2 \intd\rho_t \leq  \eta^2 \exp(\beta (\overline{J}-\underline{J})) \int |\theta\|^2 + m^{-1}\|\omega\|^2\intd\rho_t.
\end{equation*}
Then we try to prove the boundedness of $\int \|\theta\|^2 + m^{-1}\|\omega\|^2 \intd \rho_t$. Since $\rho_t$ is a weak solution of the Fokker-Planck equation, one has
\begin{equation*}
\begin{aligned}
        \frac{\intd }{\intd  t }\int \l ( \|\theta\|^2 + m^{-1}\|\omega\|^2\r )  \intd \rho_t 
        & =  
        \int \omega \cdot \theta -\gamma_1(\theta-u_t)\cdot \theta  - (\theta-u_t)\cdot \omega -\gamma m^{-1}\omega\cdot\omega \intd \rho_t\\
        & = \int -\gamma_1(\theta-u_t)\cdot \theta  +u_t\cdot \omega -\gamma m^{-1}\omega\cdot\omega \intd \rho_t
\end{aligned}
    \end{equation*}
Since 
\begin{equation*}
    \int \theta \cdot u_t \intd \rho_t \lesssim \int\|\theta\|^2 + \|u_t\|^2  \intd \rho_t \lesssim \int\|\theta\|^2 \intd \rho_t +  \int(\|\theta\|^2 +m^{-1}\|\omega\|^2) \intd \rho_t,
\end{equation*}
and 
\begin{equation*}
    \int \omega \cdot u_t \intd \rho_t \lesssim \int\|\omega\|^2 + \|u_t\|^2  \intd \rho_t \lesssim \int\|\omega\|^2 \intd \rho_t +  \int(\|\theta\|^2 +m^{-1}\|\omega\|^2) \intd \rho_t,
\end{equation*}
we can deduce that 
\begin{equation*}
    \frac{\intd}{\intd t} \int \l ( \|\theta\|^2 + m^{-1}\|\omega\|^2\r )  \intd \rho_t \lesssim \l ( \|\theta\|^2 + m^{-1}\|\omega\|^2\r )  \intd \rho_t.
\end{equation*}
Applying Gronwall's inequality yields that $\int \l ( \|\theta\|^2 + m^{-1}\|\omega\|^2\r )  \intd \rho_t$ is bounded and the above inequality is independent of $u_t$. Thus we have shown that $\|u_t\|_\infty$ is bounded by a uniform constant $q$. Theorem \ref{thm:compact} then gives the existence.

\textbf{Step 4(Uniqueness):}
Suppose we are given two fixed points of $T$: $u_t$ and $\hat u _t $ with $\|u\|_{\infty},\|\hat u\|_\infty\leq q $ and $\sup _{t\in [0,T]}\int \|\theta\|^4 +\|\omega\|^4 \intd \rho_t ,\sup _{t\in [0,T]}\int \|\hat\theta\|^4 +\|\hat\omega\|^4 \intd \hat\rho_t\leq K$ and their corresponding process $(\Theta,\Omega),(\hat\Theta,\hat\Omega)$ satisfying respectively. Then take the difference $\delta \Theta := \Theta - \hat\Theta$ and $\delta \Omega := \Omega-\hat{\Omega}$.
One has 
\begin{equation*}
\begin{aligned}
    \delta \Theta_t  &= \delta \Theta_0 + \int _0 ^ t \delta \Omega_\tau \intd \tau - \gamma_1 \int _0 ^ t 
    \delta \Theta_\tau \intd \tau   + \gamma_1 \int _0^t (u_\tau-\hat u_\tau) \intd \tau ,\\
    \delta \Omega_t &=   \delta \Omega_0  -\gamma_2 \int _0 ^ t \delta \Omega_\tau \intd \tau - m \int _0 ^ t 
    \delta \Theta_\tau \intd \tau   + m \int _0^t (u_\tau-\hat u_\tau) \intd \tau .
\end{aligned}
\end{equation*}
Thus 
\begin{equation*}
\begin{aligned}
    \mathbb E [\|\delta \Theta_t \|^2 + \|\delta \Omega_t\|^2] 
    \lesssim& 
    \mathbb E [\|\delta \Theta_0 \|^2 + \|\delta \Omega_0\|^2] + \mathbb E\l [\l ( \int_0^t \l \|\delta\Omega_\tau\r  \| \intd \tau \r ) ^2\r  ] + 
    \mathbb E\l [\l ( \int_0^t \l \|\delta\Theta_\tau\r  \| \intd \tau \r ) ^2\r  ]\\
    & +    \mathbb   E\l [\l ( \int_0^t \l \|u_\tau -\hat u_\tau\r  \| \intd \tau \r ) ^2\r  ] .
\end{aligned}
\end{equation*}
For the $ \mathbb E\l [\l ( \int_0^t \l \|u_\tau -\hat u_\tau\r  \| \intd \tau \r ) ^2\r  ] $, we have that 
\begin{equation*}
    \begin{aligned}
        \mathbb E\l [\l ( \int_0^t \l \|u_\tau -\hat u_\tau\r  \| \intd \tau \r ) ^2\r  ] 
        &=  \mathbb E\l [\l ( \int_0^t \l \|\CM [\mu_\tau] -\CM[\hat\mu_\tau]\r  \| \intd \tau \r ) ^2\r  ]
   \\
   &\leq t \mathbb E\l [\int_0^t \l \|\CM [\mu_\tau] -\CM[\hat\mu_\tau]\r  \|^2 \intd \tau \r  ].
    \end{aligned}
\end{equation*}
Thus we have 
\begin{equation*}
\begin{aligned}
    \mathbb E [\|\delta \Theta_t \|^2 + \|\Omega_t\|^2] 
    \lesssim& 
    \mathbb E [\|\delta \Theta_0 \|^2 + \|\delta \Omega_0\|^2] + \mathbb E\l [\l ( \int_0^t \l \|\delta\Omega_\tau\r  \| \intd \tau \r ) ^2\r  ] + 
   \mathbb E\l [\l ( \int_0^t \l \|\delta\Theta_\tau\r  \| \intd \tau \r ) ^2\r  ]\\
    & +     \mathbb E\l [\int_0^t \l \|\CM[\mu_\tau] -\CM[\hat\mu_\tau]\r  \|^2 \intd \tau \r  ].
\end{aligned}
\end{equation*}
Notice that by Lemma \ref{lem:wasserstwin}, $\|\CM[\mu_\tau] -\CM [\hat \mu_\tau] \|\lesssim W_2(\rho_\tau,\hat\rho_\tau)\leq \sqrt{\mathbb E \l [\|\delta \Theta_\tau\|^2+\|\delta \Omega_\tau\|^2\r  ]}$. So we can deduce 
\begin{equation*}
    \mathbb E \l [\|\delta\Theta_\tau\|^2 +\|\delta\Omega_\tau\|^2\r  ]
    \lesssim 
    \mathbb E \l [\|\delta\Theta_0\|^2 +\|\delta\Omega_0\|^2\r  ] +  \mathbb E\l [\int_0^t \l \|\delta\Omega_\tau\r  \|^2 + \l \|\delta\Theta_\tau\r  \|^2\intd \tau \r  ].
\end{equation*}
By the Gronwall'sinequality with the fact that $\mathbb E[\|\delta \Theta_0\|^2 + \|\delta\Omega_0\|^2] = 0$ gives that uniqueness result.
\end{proof}
\section{Mean Field Limit}\label{app:mean_filed}
In this section, we prove the connection between the solution to Equation \eqref{equ:CBO_SDE} and the solution of the Fokker-Planck equation \eqref{equ:FP_M_CBO}. We begin with the following boundedness result. 
 \begin{lemma}\label{equ:boundedness}
    Let $\mathcal J$ satisfy Assumption \ref{assum:Lip_J} and $\rho_0\in \PP{4}$. For any $N
\geq 2$, assume that \( \{(\Theta_t^{(i,N)},\Omega_t^{(i,N)})_{t\in [0,T]}\}_{i=1}^N\) is the unique solution to Equation \eqref{equ:M_CBO_fixed_N}  with $ \rho_0^{\otimes N}$-distributed initial data \( \{(\Theta_0^{(i,N)},\Omega_0^{(i,N)})\}_{i=1}^N\). Then there exists a constant $K>0$ independent of $N$ such that 
\begin{equation*}
    \begin{aligned}
        \sup_{i=1\cdots N }\l \{\sup_{t\in [0,T]}\mathbb E \l [\|\Theta_t^{(i,N)}\|^2 + \|\Omega_t^{(i,N)}\|^2\r  ]\r  \} &\leq K, \\
        \sup_{i=1\cdots N }\l \{\sup_{t\in [0,T]}\mathbb E \l [\|\Theta_t^{(i,N)}\|^4 + \|\Omega_t^{(i,N)}\|^4\r  ]\r  \} &\leq K, \\
        \sup_{i=1\cdots N }\l \{\sup_{t\in [0,T]}\mathbb E \l [\|\CM [\hat \mu^N_t]\|^2 \r  ]\r  \} &\leq K, \\
        \sup_{i=1\cdots N }\l \{\sup_{t\in [0,T]}\mathbb E \l [\|\CM [\hat \mu^N_t]\|^4 \r  ]\r  \} &\leq K, \\
    \end{aligned}
\end{equation*}
where $\hat \mu^N_t = \frac{1}{N} \sum_{i=1}^N \delta_{\Theta_t^{(i,N)}} $ is the empirical measure.
\end{lemma}
\begin{proof}
    For each $i$, we have
    \begin{equation*}
        \begin{aligned}
    \intd \Theta _t ^{(i,N)} =&    - \Omega _t ^{(i,N)} \intd t  -  \gamma_1 ( \Theta _t ^{(i,N)}  - \CM [\hat \mu_t ^N]) \intd t + \sigma(t)\intd W^i_{\theta,t},\\
    \intd \Omega _t ^{(i,N)} =&  - m ( \Theta _t ^{(i,N)}  - \CM [\hat \mu_t ^N ]) \intd t -\gamma_2  \Omega _t ^{(i,N)}\intd t +  \sqrt{m} \sigma(t) \intd W^i_{\omega,t}.
        \end{aligned}
    \end{equation*}
    Now we pick $p=1$ or $p=2$. Then 
    \begin{equation*}
        \begin{aligned}
            \mathbb E \l [\l \|\Omega _t ^{(i,N)} \r \|^{2p} \r  ]  +  \mathbb E \l [ \l \|\Theta _t ^{(i,N)} \r \|^{2p} \r  ] \lesssim & \mathbb E \l [ \l \|\Theta _0  ^{(i,N)}\r \|^{2p}  \r  ] + \mathbb E \l [\l \|\Omega _0  ^{(i,N)} \r\|^{2p}  \r  ] \\
            & + \mathbb E \l [\int_0^t \l \|\Theta _\tau  ^{(i,N)} \r \| \intd \tau  \r  ]^{2p} +\mathbb E \l [\int_0^t \l  \|\Omega _\tau  ^{(i,N)}\r\| \intd \tau  \r  ]^{2p} 
            \\
            &  + \mathbb E \l [\int_0^t \l \|\CM [\hat \mu_\tau  ^N] \r \| \intd \tau  \r  ]^{2p}+  \mathbb E \l [\int_0^t \sigma(\tau ) \intd W_{\theta,\tau}  \r  ]^{2p}.
        \end{aligned}
    \end{equation*}
    Now apply the Cauchy's inequality,
    \[
             \mathbb E \l [\int_0^t \|\Theta _\tau  ^{(i,N)}\| \intd \tau  \r  ]^{2p}  \leq  t^p \cdot \mathbb E \l [\int_0^t \|\Theta _\tau  ^{(i,N)}\|^2 \intd \tau  \r  ]^{p} ,
        \quad 
             \mathbb E \l [\int_0^t \|\Omega _\tau  ^{(i,N)}\| \intd \tau  \r  ]^{2p}  \leq  t^p \cdot \mathbb E \l [\int_0^t \|\Omega _\tau  ^{(i,N)}\|^2 \intd \tau  \r  ]^{p} ,
    \]
    and 
    \[
    \mathbb E \l [\int_0^t \|\CM [\hat \mu_\tau ^N ]\| \intd \tau\r  ]^{2p}  \leq t^p  \mathbb E \l [\int_0^t \|\CM [\hat \mu_\tau  ^N]\| ^2\intd \tau\r  ]^{p} .
    \]
    Also, by It\^o Isometry,
    \[
    \mathbb E \l [\int_0^t \sigma(\tau ) \intd W_{\theta,\tau}  \r  ]^{2p} =   \mathbb E \l [\int_0^t \sigma(\tau )^2  \intd \tau  \r  ]^{p}.
    \]
    Thus 
    \begin{equation*}
        \begin{aligned}
    \mathbb E \l [ \l \|\Omega _t ^{(i,N)}\r\|^{2p} + \l \|\Theta _t ^{(i,N)} \r \|^{2p} \r  ] \lesssim & \mathbb E \l [\l \|\Omega _0 ^{(i,N)} \r \|^{2p} \r  ]  +  \mathbb E \l [\l \|\Theta _0 ^{(i,N)}\r \|^{2p} \r  ]  \\
    & + \mathbb E \l [\int_0^t \l \|\Theta _\tau  ^{(i,N)} \r \|^2 \intd \tau  \r  ]^{p} + \mathbb E \l [\int_0^t \l\|\Omega _\tau  ^{(i,N)}\r\|^2 \intd \tau  \r  ]^{p} \\
    & +  \mathbb E \l [\int_0^t \l \|\CM [\hat \mu_\tau ]\r\| ^2\intd \tau\r  ]^{p}   + 1.
        \end{aligned}
    \end{equation*}
    Further, by H\"older inequality,
    \[
    \mathbb E \l [\int_0^t \l \|\Theta _\tau  ^{(i,N)} \r \|^2 \intd \tau  \r  ]^{p} \leq  \mathbb E \l [\int_0^t \l \|\Theta _\tau  ^{(i,N)}\r \|^{2p} \intd \tau  \r  ],
    \quad 
    \mathbb E \l [\int_0^t \|\Omega _\tau  ^{(i,N)}\|^2 \intd \tau  \r  ]^{p} \leq  \mathbb E \l [\int_0^t \|\Omega _\tau  ^{(i,N)}\|^{2p} \intd \tau  \r  ]
    \] and 
    \[
    \mathbb E \l [\int_0^t \l \|M_\beta [\hat \mu_t ^N ] \r \| ^2\intd \tau\r  ]^{p} \leq \mathbb E \l [\int_0^t \l \|M_\beta [\hat \mu_t  ^N] \r | ^{2p}\intd \tau\r  ].
    \]
    So we can deduce 
    \begin{equation*}
        \begin{aligned}
    \mathbb E \l [\l \|\Omega _t ^{(i,N)}  \r \|^{2p} + \l \|\Theta _t ^{(i,N)} \r \|^{2p} \r  ] \lesssim & \mathbb E \l [ \l \|\Omega _0 ^{(i,N)}\r \|^{2p} \r  ]  +  \mathbb E \l [\l \|\Theta _0 ^{(i,N)}\r \|^{2p} \r  ]  \\
    & + \mathbb E \l [\int_0^t \l \|\Theta _\tau  ^{(i,N)}\r \|^{2p} \intd \tau \r  ]  + \mathbb E \l [\int_0^t \l \|\Omega _\tau  ^{(i,N)}\r\|^{2p} \intd \tau  \r  ] \\
    & +  \mathbb E \l [\int_0^t \|\CM [\hat \mu_\tau ^N ]\| ^{2p}\intd \tau\r  ]   + 1.
        \end{aligned}
    \end{equation*}
    Thus 
    \begin{equation*}
        \begin{aligned}
            \mathbb E \l[ \int \l ( \|\theta\|^{2p} +\|\omega\|^{2p}  \r ) \intd \hat \rho_t^N \r] \lesssim & \mathbb E \l[ \int \l (\|\theta\|^{2p} +\|\omega\|^{2p} \r) \intd \hat \rho_0^N \r]\\
            & + \int _0 ^ t\l (  \mathbb E \l [ \int \l(\|\theta\|^{2p} +\|\omega\|^{2p} \r) \intd \hat \rho_\tau^N\r ] \r)  \intd \tau  \\
            & +  \int _0^t  \mathbb E \l [\|\CM [\hat \mu_\tau ^N ]\| \r  ]^{2p} \intd \tau + 1.
        \end{aligned}
    \end{equation*}
    It follows from Lemma 3.1 in \cite{carrillo2018analytical}, we have
    \begin{equation*}
        \int \|\theta\|^2 \frac{w_\beta(\theta)}{\|w_\beta\|_{L^1(\hat\rho^N_\tau)}} \intd \hat\rho_\tau^N \leq b_1 + b_2 \int \|\theta\|^2   \intd \hat\rho _\tau^N\leq b_1 + b_2 \int (\|\theta\|^2+\|\omega\|^2)   \intd \hat\rho _\tau^N.
    \end{equation*}
    Then we can calculate 
    \begin{equation*}
        \begin{aligned}
            \|M_\beta [\hat \mu_t ^N]\| ^{2p}  &= \l \|\int\theta\cdot \frac{w_\beta(\theta)}{\|w_\beta\|_{L^1(\hat\rho^N_\tau)}} \intd \hat\rho_\tau^N  \r  \| ^{2p}\\
             \leq &\l ( \int \|\theta\| \frac{w_\beta(\theta)}{\|w_\beta\|_{L^1(\hat\rho^N_\tau)}} \intd \hat\rho_\tau^N  \r )  ^{2p}\\
             \leq & \l ( \int \|\theta\|^2  \frac{w_\beta(\theta)}{\|w_\beta\|_{L^1(\hat\rho^N_\tau)}} \frac{w_\beta(\theta)}{\|w_\beta\|_{L^1(\hat\rho^N_\tau)}}  \intd \hat\rho_\tau^N  \r )  ^{p}\\
             \lesssim & 
            \l ( \int \|\theta\|^2  \frac{w_\beta(\theta)}{\|w_\beta\|_{L^1(\hat\rho^N_\tau)}} \intd \hat\rho_\tau^N  \r )  ^{p}\\
            \leq&  \l ( b_1 + b_2  \int (\|\theta\|^2+\|\omega\|^2)   \intd \hat\rho _\tau^N\r ) ^{2p}\\
            \lesssim &  1 + \int (\|\theta\|^{2p}+\|\omega\|^{2p})   \intd \hat\rho _\tau^N.
        \end{aligned}
    \end{equation*}
    Combining the above inequality leads to 
    \begin{equation*}
        \begin{aligned}
            \mathbb E \l[\int (\|\theta\|^{2p} +\|\omega\|^{2p} ) \intd \hat \rho_t^N \r]\lesssim & \mathbb E \l [ \int \l (\|\theta\|^{2p} +\|\omega\|^{2p}\r) \intd \hat \rho_0^N \r] \\
            & + \int _0 ^ t\l (  \mathbb E \l[\int \l(\|\theta\|^{2p} +\|\omega\|^{2p} \r)\intd \hat \rho_\tau^N \r]\r )  \intd \tau  +1.
        \end{aligned}
    \end{equation*}
    By applying Gronwall's inequality, it follows that
    $\mathbb E \l[\int\l( \|\theta\|^{2p} +\|\omega\|^{2p} \r)\intd \hat \rho_t^N\r] $ is bounded for $t\in [0,T]$ and the bound does not depend on $N$. Also, we know that
    \[
    \|M_\beta [\hat\mu^N_t \rho_t ]\| ^{2p}  \lesssim  1 + \int (\|\theta\|^{2p}+\|\omega\|^{2p})   \intd \hat\rho _\tau^N,
    \]
    which implies that 
    \[
    \mathbb E\l[ \|M_\beta [\hat \mu_t^N ]\| ^{2p} \r] \lesssim  1 +  \mathbb E\l[ \int(\|\theta\|^{2p}+\|\omega\|^{2p})   \intd \hat\rho _\tau^N\r].
    \]
    So $\mathbb E\l[  \|M_\beta [\hat \mu_t^N ]\| ^{2p}\r]  $ is bounded for $t\in [0,T]$ and the bound does not depend on $N$.
\end{proof}
Now we treat \( \l(\Theta^{(i,N)},\Omega^{(i,N)}\r)\) as a random variable defined on \( (\Omega,\mathcal F, \mathbb P) \) and taking values in  \(\MC{[0,T];\mbRd\times\mbRd}  \). Then $\hat \rho^N= \frac{1}{N}\sum_{i=1}^N \delta_{( \Theta^{(i,N)},\Omega^{(i,N)})}$ is a random measure. Let us denote $\mathcal L (\hat \rho^N):=Law(\rho^N) \in \mathcal P \l(\mathcal P \l(\mathcal C([0,T],\mathcal R^{D} )\times \mathcal C\l([0,T],\mathcal R^{D} \r)\r) \r) $ as a sequence of probability distributions. We can prove that $\l\{\mathcal L (\rho^N )\r\}_{N\geq 2}$ is tight.
Next, we use the Aldous criteria (\cite{bass2011stochastic}, Section 34.3), which could prove the tightness of a sequence of distributions.
\begin{theorem}\label{thm:tightness}
    Let $\mathcal J$  satisfy Assumptions \ref{assum:Lip_J} and $\rho_0 \in \PP{4}$. For any $N\geq 2$, we assume that $\l\{\l(\Theta_t^{(i,N)},\Omega_t^{(i,N)}\r)_{t\in [0,T]}\r\}_{i=1}^N$ is a unique solution to Equation \eqref{equ:M_CBO_fixed_N} with $\rho_0^{\otimes N}$ -distributed initial data $\l\{\l(\Theta_0^{(i,N)},\Omega_0^{(i,N)}\r)\r\}_{i=1}^N$. Then the sequence $\l\{\mathcal L (\hat\rho^N)\r\}_{N\geq 2}$ is tight in $\mathcal P \l(\mathcal P \l(\mathcal C([0,T],\mathcal R^{D} )\times \mathcal C\l([0,T],\mathcal R^{D} \r)\r) \r) $.
\end{theorem}
\begin{proof}
    Because of the exchangeability of the particle system, we only prove the $ \l\{\mathcal L \l( \Theta_t ^{1, N}, \Omega_t^{1, N}\r) \r\}_{N\geq 2}$ is tight. It is sufficient to justify two conditions in Aldous criteria.\\
    \textbf{Condtion 1:} For any $\epsilon>0$, there exist a compact subset $U_\epsilon := \l\{\l(\theta,\omega\r): \|\theta\|^2 +\|\omega\|^2 ]\leq \frac{K}{\epsilon}\r\}$ such that by Markov's inequality 
    \[
    \mathcal L \l( \Theta_t ^{1,N} ,\Omega_t^{1,N}\r) \l((U_\epsilon)^c\r) = \mathbb P \l ( \l\| \Theta_t ^{1,N}\r\|+\l\|\Omega_t^{1,N}\r\| > \frac{K}{\epsilon}\r )   \leq \frac{\epsilon \mathbb E [\| \Theta_t ^{1,N}\|+\|\Omega_t^{1,N}\| ] }{K} \leq \epsilon, \quad \forall N\geq 2,
    \]
    where we have used Lemma \ref{equ:boundedness} in the last inequality. This means that for each $t\in [0,T]$, the sequence $\l\{\mathcal L  \l( \Theta_t ^{1,N} ,\Omega_t^{1,N}\r)\r\}$ is tight.

    \textbf{Condition 2:} 
    {We have to show, for any  $\epsilon,\eta>0$, there exist $\delta_0>0$ and $n_0\in \mathbb N$ such that for all $N\geq n_0$.} Let $\tilde \tau$ be a $\sigma ((\Theta_{s} ^{1,N},\Omega_{s} ^{1,N});s \in [0,T])$-stopping time with discrete values such that $ \tilde \tau + \delta_0 \leq T$, it holds that  
    \[
    \sup _{\delta \in [0,\delta_0]} \mathbb P \l (  \l \| \Theta_{t+\delta} ^{1,N} - \Theta_{t} ^{1,N}\r \| \geq \eta \r )  \leq \epsilon 
    \]
    and
        \[
    \sup _{\delta \in [0,\delta_0]} \mathbb P \l (  \l \| \Omega_{t+\delta} ^{1,N} - \Omega_{t} ^{1,N}\r \| \geq \eta \r )  \leq \epsilon.
    \]
    Recalling Equation \eqref{equ:CBO_SDE}, we have
    \begin{equation*}
        \begin{aligned}
               \Theta_{\tilde \tau +\delta} ^{1,N} - \Theta_{\tilde \tau} ^{1,N} &=  \int_{\tilde \tau }^{\tilde \tau +\delta} \Omega_s ^{1,N} \intd s 
     -  \gamma_1 \int_{\tilde \tau }^{\tilde \tau +\delta}   ( \Theta_s ^{1,N} -\CM [\hat\mu^N_s]) \intd s 
     +  \int_{\tilde \tau }^{\tilde \tau +\delta} \sigma(s)\intd W_{\theta,t}, \\
        \Omega_{\tilde \tau +\delta} ^{1,N} - \Omega_{\tilde \tau} ^{1,N} &= 
     -  m \int_{\tilde \tau }^{\tilde \tau +\delta}   ( \Theta_s ^{1,N} -\CM  [ \hat\mu^N_s]) \intd s 
    - \gamma_2 \int_{\tilde \tau }^{\tilde \tau +\delta} \Omega_s ^{1,N} \intd s 
     +   \sqrt{m} \int_{\tilde \tau }^{\tilde \tau +\delta} \sigma(s)\intd W_{\omega,t}.
        \end{aligned}
    \end{equation*}
    From Theorem 2.1 in \cite{Huang2022mean}, we have 
    \[
    \mathbb E \l [\l | \int_{\tilde \tau }^{\tilde \tau +\delta}   ( \Theta_s ^{1,N} -\CM [\hat\mu^N_s]) \intd s  \r  | ^2 \r  ] \leq 2 T K \delta.
    \]
    Furthermore, we apply It\^o's Isometry 
    \[
    \mathbb E\l [\int_{\tilde \tau }^{\tilde \tau +\delta} \sigma(s)\intd W_{\theta,s} \r  ],\mathbb E\l [\int_{\tilde \tau }^{\tilde \tau +\delta} \sigma(s)\intd W_{\omega,s} \r  ] \leq \l ( \mathbb E\l [\int_{\tilde \tau }^{\tilde \tau +\delta} \sigma(s)^2\intd s  \r  ] \r ) ^{\frac{1}{2}} \leq \overline{\sigma}^2 \delta^{\frac{1}{2}}T ^{\frac{1}{2}}. 
    \]
    Combining the above estimation, one has
    \[
    \mathbb E\l [ \l | \Theta_{\tilde \tau +\delta} ^{1,N} - \Theta_{\tilde \tau} ^{1,N}\r  |^2 +
    \l | \Omega_{\tilde \tau +\delta} ^{1,N} - \Omega_{\tilde \tau} ^{1,N}\r  |^2 \r  ] \lesssim \sqrt{\delta}.
    \] 
    Hence, for any $\epsilon,\eta>0$, there exist $\delta_0>0$ such that for all $N>2$ it holds that 
    \begin{equation*}
        \begin{aligned}
    & \sup _{\delta \in [0,\delta_0]} \mathbb P \l ( \l  \| \Theta_{t+\delta} ^{1,N} - \Theta_{t} ^{1,N} \r \|^2 \geq \eta \r )  
    ,
    \sup _{\delta \in [0,\delta_0]} \mathbb P \l ( \l \| \Omega_{t+\delta} ^{1,N} - \Omega_{t} ^{1,N} \r \|^2 \geq \eta \r )  \\
    \leq &
    \sup _{\delta \in [0,\delta_0]} \mathbb P \l ( \l  \| \Theta_{t+\delta} ^{1,N} - \Theta_{t} ^{1,N} \r \|^2+\| \Omega_{t+\delta} ^{1,N} - \Omega_{t} ^{1,N} \|^2  \geq \eta \r )  \\
    \leq & \sup _{\delta \in [0,\delta_0]} \frac{\mathbb  E\l [ \l | \Theta_{\tilde \tau +\delta} ^{1,N} - \Theta_{\tilde \tau} ^{1,N}\r  |^2 +
    \l | \Omega_{\tilde \tau +\delta} ^{1,N} - \Omega_{\tilde \tau} ^{1,N}\r  |^2 \r  ]}{\eta}\leq \epsilon.
        \end{aligned}
    \end{equation*}
\end{proof}

By Skorokhod's lemma (see \cite{billingsley2013convergence}) and Lemma \ref{lemma:weak_conv_measure}, we may find a common probability space $(\Omega,\mathcal{F},\mathbb P)$ on which the process $\{\hat \rho^N\}_{N\in\mathbb N}$ converges to some process $\rho$ as a random variable valued in $\mathcal {P}\lg\mcCRd\times \mcCRd \rg$ almost surely. In particular, we have that $\forall t\in [0,T]$ and $\phi\in \mathcal C_b(\mbRd\times\mbRd)$,
\begin{equation}
\lim_{N\rightarrow \infty} \l |\la \phi,\hat \rho_t^N-\rho_t\ra \r| + \left|\CM[\hat \mu^N_t]-\CM[\mu_t]\right|= 0\quad \text{a.s.},
\end{equation}
 and we can have a direct result 
 \begin{equation}\label{equ:e_conve}
    \lim_{N\rightarrow \infty} \mathbb E\l[ \l |\la \phi,\hat \rho_t^N-\rho_t\ra \r| + \left|\CM[\hat \mu^N_t]-\CM[\mu_t]\right|\r]= 0.
 \end{equation}
\begin{lemma}\label{lemma:weak_conv_measure}
    \begin{enumerate}
        \item There exist a subsequence of $\{\hat \rho^N\}_{N\geq 2}$ and a random measure $\rho :  \Omega \to \mathcal P \l (\mcCRd\times\mcCRd \r )$ such that $\hat \rho^N \rightharpoonup \rho$ in law as $N\to \infty$
        which is equivalently to say $\mathcal L \l (\hat\rho^N \r )$ converges weakly to $\mathcal L \l ( \rho \rg$ in $\mathcal P \l ( 
         \mathcal P \l (\mcCRd\times \mcCRd \rg\rg$;
         \item For the subsequence in 1, the time marginal $\hat \rho_t^N$ of  $\hat \rho^N$,  as $\mathcal {P}\lg\mbRd\times\mbRd\rg$ valued random measure converges in law to $\rho_t\in \mathcal{P}(\mbRd\times\mbRd)$, the time marginal of $\rho$. Namely $\mathcal{L}(\hat \rho_t^N)$ converges weakly to $\mathcal {L}(\rho_t)$ in $\mathcal {P}\l (\mathcal {P}\lg\mbRd\times \mbRd\rg\rg$.

    \end{enumerate}
\end{lemma}

\begin{definition}
Fix $\phi\ \in \mathcal C^2_c \lg\mbRd \times \mbRd\rg$, define functional $\mathcal F_\phi: \mathcal  P \lg\mcCRd \times \mcCRd \r )\to \mathbb R$:
\begin{equation}
\begin{aligned}
F_{\phi}(\rho_t):=& \la \phi(\theta_t,\omega_t) ,\rho (\intd \theta,\intd \omega )\ra-\la \varphi(\theta_0,\omega_0),\rho(\intd \theta,\intd \omega ) \ra \\
& -\int_0^t\la \l (  \omega_s -\gamma_1 (\theta_s-\CM[ \mu_s])\r )  \cdot\nabla_\theta \phi,\rho(\intd \theta,\intd \omega)\ra \intd s  \\
& - \int_0^t\la \l (  -m  (\theta_s-\CM[\mu_s]) -\gamma_2\omega_s\r )  \cdot\nabla_\omega \phi,\rho(\intd \theta,\intd \omega)\ra \intd s  \\
& - \frac{(m+1)D }{2}\int_0^t \sigma^2(s) 
\intd s 
\end{aligned}
\end{equation}
for all $\rho \in \mathcal  P \l (\mcCRd \times \mcCRd \r) $ and $\theta,\omega \in \mcCRd $, where $\mu(\theta)=\int_{ \mbRd }\rho(\theta,\intd \omega)$.
\end{definition}
\begin{lemma}
    Let \(\mathcal J\) satisfy Assumption \ref{assum:Lip_J} and $\rho_0\in \mcP_4(\mbRd\times \mbRd)$. For all $N\geq 2$, assume that \(
\l\{\l(\Theta_{t}^{(i,N)} ,\Omega_{t} ^{(i,N)}\r)_{t\in[0,T]}  \r\}  _{i=1}^N
    \) is the unique solution to Equation \eqref{equ:M_CBO_fixed_N} with $\rho_0^{\otimes N}$-distributed initial data $\l \{ \l(\Theta_0^{(i,N)}, \Omega_0^{(i,N)}\r)\r\}_{i=1}^N$. There exists a constant $C>0$, such that 
    \[
    \mbbE \l [\l |F_\phi(\hat\rho^N)\r  |^2\r  ]\leq \frac{C}{N},
    \]
    where $\hat\rho^N= \frac{1}{N}\sum_{i=1}^N \delta_{( \Theta^{(i,N)},\Omega^{(i,N)})}$ is the empirical measure.
\end{lemma}
\begin{proof}
    Using the definition of $F_\phi$, one has
    \begin{equation*}
        \begin{aligned}
        	F_{\phi}(\hat\rho^N)  = &  \frac{1}{N} \sum_{i=1}^N \phi(\Theta_t^{(i,N)},\Omega_t^{(i,N)} ) -  \frac{1}{N} \sum_{i=1}^N \phi(\Theta_0^{(i,N)},\Omega_0^{(i,N)} ) \\
& - \frac{1}{N} \sum_{i=1}^N \int_0^t \l (  \Omega_s^{(i,N)}-\gamma_1 \l ( \Theta_0^{(i,N)} -\CM[ \hat\mu^N_s]\r ) \r )  \cdot\nabla_\theta \phi(\Theta_s^{(i,N)},\Omega_s^{(i,N)} ) \intd s \\
         & - \frac{1}{N} \sum_{i=1}^N \int_0^t \l (  -m  \l ( \Theta_s^{(i,N)} -\CM[ \hat\mu^N_s]\r ) -\gamma_2 \Omega_s^{(i,N)}\r )  \cdot\nabla_\omega \phi(\Theta_s^{(i,N)},\Omega_s^{(i,N)} ) \intd s \\
         & - \frac{(m+1)D }{2} \int_0^t \sigma^2(s) \intd s .
        \end{aligned}
    \end{equation*}
    One the other hand, the It\^o-Doeblin formula gives
    \begin{equation*}
        \begin{aligned}
             \phi\l (  \Theta_t^{(i,N)},\Omega_t^{(i,N)} \r )    - \phi \l (\Theta_0^{(i,N)},\Omega_0^{(i,N)} \r )= & 
              \int_0^t \l (  \Omega_s^{(i,N)}-\gamma_1 \l ( \Theta_s^{(i,N)} -\CM[\hat\mu^N_s]\r ) \r )  \cdot\nabla_\theta \phi(\Theta_s^{(i,N)},\Omega_s^{(i,N)} ) \intd s \\
            & +  \int_0^t \l (  -m  \l ( \Theta_s^{(i,N)} -\CM[\hat\mu^N_s]\r ) -\gamma_2 \Omega_s^{(i,N)}\r )  \cdot\nabla_\omega \phi(\Theta_s^{(i,N)},\Omega_s^{(i,N)} ) \intd s \\
            &  + \int_0^t \sigma(s) \intd W^i_{\theta,s}+  \sqrt{m}\int_0^t \sigma(s) \intd W^i_{\omega,s} \\
            & + \frac{(m+1)D }{2} \int_0^t \sigma^2(s) \intd s .
        \end{aligned}
    \end{equation*}
    Then one gets 
    \[
    F_{\phi}(\hat\rho^N)  = \frac{1}{N}\sum_{i=1}^N\lg\int_0^t   \sigma(s) \intd W^i_{\theta,s} +  \sqrt{m} \int_0^t   \sigma(s) \intd W^i_{\omega,s}\rg. 
    \]
    Finally, we can compute
    \begin{equation*}
    \begin{aligned}
    \mbbE \l [ \l|  F_{\phi}(\hat\rho^N) \r| ^2\r ] = & \frac{1}{N^2}\sum_{i=1}^N\mbbE \l [ \l|  \int_0^t   \sigma(s) \intd W^i_{\theta,s} +  \sqrt{m} \int_0^t   \sigma(s) \intd W^i_{\omega,s} \r|\r ]^2 \\
        \leq &  T \frac{\overline{\sigma}^2 (m+1)}{N^2},
    \end{aligned}
    \end{equation*}
    where we use the assumption that the $\sigma(t)$ we choose has an upper bound.
\end{proof}

\begin{proof}[Proof of Theorem \ref{thm:mean_field_limit}]
 Now suppose that we have a convergent subsequence of $ \{\hat{\rho}^N\}_{N\in\mathbb{N}} $, which is denoted by the sequence itself for simplicity and has $ \rho_{t} $ as the limit. We want to prove that $\rho$ is a solution of \eqref{equ:FP_M_CBO}.  For any $\phi\in \mathcal C_c^2\l({\mbRd\times \mbRd}\r)$, using the convergence result in \eqref{equ:e_conve}, we have
	\begin{equation*}
	    \lim_{N\to \infty} \mbbE \l [ \l| \la \phi(\theta,\omega) ,\hat\rho^N_t (\intd \theta,\intd \omega )\ra
     -\la \phi(\theta,\omega) ,\rho_t (\intd \theta,\intd \omega )\ra
     \r| \r]= 0,
	\end{equation*}
 and 
\begin{equation*}
	    \lim_{N\to \infty} \mbbE \l [ \l| \la \phi(\theta,\omega) ,\hat\rho^N_0 (\intd \theta,\intd \omega )\ra
     -\la \phi(\theta,\omega) ,\rho_0 (\intd \theta,\intd \omega )\ra
     \r| \r]= 0.
\end{equation*}
Further, we notice that 
\begin{equation*}
    \begin{aligned}
         & \l\| \int_0^t\la \l (\theta-\CM[\hat\mu^N_s]   \r )\cdot\nabla_\theta \phi,\hat\rho^N_s(\intd \theta,\intd \omega)\ra \intd s - \int_0^t\la \lg\theta-\CM[ \mu_s]    \rg\cdot\nabla_\theta \phi,\rho_s(\intd \theta,\intd \omega)\ra \intd s   \r\|\\
        \leq & \int_0^t\l\|\la\l (\theta-\CM[\hat \mu^N_s]   \r ) \cdot\nabla_\theta \phi,\hat\rho^N_s(\intd \theta,\intd \omega)
        - \rho_s(\intd \theta,\intd \omega)\ra \r\|\intd s  \\
        & +\int_0^t \l\|\la\l (\CM[\mu_s] -\CM[\hat\mu^N_s]  \r ) \cdot\nabla_\theta \phi,\rho_s(\intd \theta,\intd \omega)\ra  \r\|\intd s  \\
        := &   \int_0^t \l\| I_1^N(s)  \r\|\intd s + \int_0^t \l\| I_2^N(s)  \r\|\intd s.
    \end{aligned}
\end{equation*}
For $\l\|I_1^N(s)\r\|$, we have 
\begin{equation}\begin{aligned}
       \mbbE\l[ \|I_1^N(s)\|\r] =&  \mathbb E\l[\l\|\la \theta   \cdot\nabla_\theta \phi,\rho^N_s(\intd \theta,\intd \omega)
        - \rho_s(\intd \theta,\intd \omega)\ra \r\|\r]  + \mathbb E\l[\l\|\la\l (\CM[\mu^N_s]   \r ) \cdot\nabla_\theta \phi,\rho^N_s(\intd \theta,\intd \omega)
        - \rho_s(\intd \theta,\intd \omega)\ra \r\|\r]\\
        \leq &  \mathbb E\l[\l\|\la \theta     \cdot\nabla_\theta \phi,\rho^N_s(\intd \theta,\intd \omega)
        - \rho_s(\intd \theta,\intd \omega)\ra \r\|\r]  + K^{\frac{1}{2}}\mathbb E\l[\l\|\la \nabla_\theta \phi,\rho^N_s(\intd \theta,\intd \omega)
        - \rho_s(\intd \theta,\intd \omega)\ra \r\|\r].
\end{aligned}
\end{equation}
Since $\phi$ has a compact support, applying \eqref{equ:e_conve} leads to \(\lim_{N\to\infty}\mbbE \l[\|I_1^N(s)\|\r]=0\). Moreover, by the uniform boundedness of  \(\mbbE \l[I_1^N(s)\r]=0\) and applying the domained convergence theorem implies 
\[
\lim_{N\to \infty}\int_0^t   \mbbE\l[ \|I_1^N(s)\|\r] \intd s =0.
\]
As for $I_2^N$, we know that 
\[
	 \left\|\la (\CM(\mu_s)-\CM(\hat\mu_s^N))\cdot \nabla\phi(x), \mu_s (dx)\ra\right\|\leq \|\nabla \phi\|_\infty  \|\CM(\mu_s)-\CM(\hat\mu_s^N)\|\,.
\]
Hence by Equation \eqref{equ:e_conve}, we have \(\lim_{N\to\infty}\mbbE \l[ \|I_2^N(s)\|\r]=0\). Again, by the dominated convergence theorem, we have 
\[
\lim_{N\to \infty}\int_0^t   \mbbE\l[ \|I_2^N(s)\|\r] \intd s =0 .
\]
Thus we can get the boundedness
\[
         \lim_{N\to \infty}\mbbE \l[ \l| \int_0^t\la \l (\theta-\CM[\hat\mu^N_s]   \r )\cdot\nabla_\theta \phi,\rho^N_s(\intd \theta,\intd \omega)\ra \intd s - \int_0^t\la \lg\theta-\CM[\hat\mu_s]    \rg\cdot\nabla_\theta \phi,\rho_s(\intd \theta,\intd \omega)\ra \intd s   \r|\r]=0.
\]
Similarly, we can also have
\[
           \lim_{N\to \infty}\mbbE \l[ \l| \int_0^t\la \l (\theta-\CM[\hat\mu^N_s]   \r )\cdot\nabla_\omega \phi,\rho^N_s(\intd \theta,\intd \omega)\ra \intd s - \int_0^t\la \lg\theta-\mathcal M_\beta[\hat\mu_s]    \rg\cdot\nabla_\omega \phi,\rho_s(\intd \theta,\intd \omega)\ra \intd s   \r|\r]=0.
\]
Combining the above results, we get 
\[
    \mbbE \l[ \l|F_\phi(\hat\rho^N )-F_\phi(\rho )\r|\r] = 0,
\]
{which is a direct result of \eqref{equ:e_conve} and the dominated convergence theorem.}
We can deduce 
\begin{equation*}
    \begin{aligned}
\mbbE \l[ \l|F_\phi(\rho )\r|\r] \leq &\mbbE \l[ \l|F_\phi(\hat\rho^N )-F_\phi(\rho )\r|\r]  + \mbbE \l[ \l|F_\phi(\hat\rho^N )\r|\r]\\
& \leq \mbbE \l[ \l|F_\phi(\hat\rho^N )-F_\phi(\rho )\r|\r]  +\sqrt{\frac{C}{N}} \to 0, \quad \text{ as  } \quad N\to \infty.
    \end{aligned}
\end{equation*}
Thus $F_\phi(\rho_t)=0$ almost surely, which implies that $\rho_t$ is a solution to the corresponding Fokker-Planck equation.
Combined with the uniqueness of the solution, proved in Lemma \ref{equ:unique}, we complete the proof.
\end{proof}
\begin{lemma}
    $\forall T>0$, let $f_t\in \MC{[0,T];\mbRd}$ and $\rho_0\in \mathcal P _2(\mbRd\times \mbRd)$. The following linear PDE 
    \begin{equation*}
        \begin{aligned}
              \partial_t \rho_t = -\nabla_\theta\l ( \l ( \omega-\gamma_1(\theta-f_t)\r ) \rho_t\r )  + \nabla_\omega\l ( m (\theta - f_t) +\gamma_2 \omega\rho_t\r )  +\frac{\sigma(t)^2m}{2} \Delta_\omega \rho_t+\frac{\sigma(t)^2 }{2} \Delta_\theta \rho_t
        \end{aligned}
    \end{equation*}
    has a unique solution $\rho_t \in \MC{[0,T];  \mathcal P _2(\mbRd\times \mbRd)}$.
\end{lemma}
\begin{proof}
    The existence is obvious, which can be obtained as the law of the solution to the associated linear SDE. To show  the uniqueness, let us fix $t_0\in [0,T]$ and $\phi\in \mathcal C_c^\infty\l(\mbRd\times\mbRd\r)$, we then can solve the following backward PDE:
\[
\partial_t h_t = -\left( \omega - \gamma_1 (\theta - f_t) \right) \cdot \nabla_\theta h_t + \left( m (\theta - f_t) + \gamma_2 \omega \right) \cdot \nabla_\omega h_t - \frac{\sigma(t)^2}{2} \left( m \Delta_\omega h_t + \Delta_\theta h_t \right),
\]
where $(t,\theta,\omega) \in [0,t_0]\times \mbRd\times \mbRd$ with terminal condition $h_{t_0} =  \phi$. It has a classical solution:
\[
h_t = \mbbE \l[\phi(\Theta_{t_0}^{t,(\theta,\omega)},\Omega_{t_0}^{t,(\theta,\omega)})\r],
\]
where $\l (  \Theta_{s}^{t,(\theta,\omega)},\Omega_{s}^{t,(\theta,\omega)} \rg_{0\leq t\leq s\leq t_0}$ is the strong solution to the following linear SDE:
\[
\begin{aligned}
\intd\Theta_s^{t,(\theta,\omega)} &= \left( \Omega_s^{t,(\theta,\omega)} - \gamma_1 (\Theta_s^{t,(\theta,\omega)} - f_s) \right) \intd s + \sigma(s) \intd W_{\theta,s}, \\
\intd\Omega_s^{t,(\theta,\omega)} &= -\left( m (\Theta_s^{t,(\theta,\omega)} - f_s) + \gamma_2 \Omega_s^{t,(\theta,\omega)} \right) \intd s + \sigma(s)  \sqrt{m} \intd W_{\omega,s}
\end{aligned}
\]
with terminal condition $   \l (\Theta_{t}^{t,(\theta,\omega)},\Omega_{t}^{t,(\theta,\omega)}\r )= (\theta,\omega)$.

Now, suppose \( \rho^1 \) and \( \rho^2 \) are two weak solutions of the PDE with the same initial condition \( \rho^1_0 = \rho^2_0 \). Let \( \delta \rho = \rho^1 - \rho^2 \). To show uniqueness, we need to demonstrate that \( \delta \rho_t = 0 \) for all \( t \in [0, T] \).

Using the backward PDE solution \( h_t \) as a test function, we have:
\[
\langle h_{t_0}, \delta \rho_{t_0} \rangle = \int_0^{t_0} \langle \partial_s h_s, \delta \rho_s \rangle  +  \langle  h_s, \partial_s \delta \rho_s \rangle \intd s  .
\]
By substituting the equation for \( \partial_s h_s \) from the backward PDE and integrating by parts, we get
\[
\langle h_{t_0}, \delta \rho_{t_0} \rangle = 0.
\]
The arbitrariness of \( \psi \in \mathcal{C}_c^\infty(\mathbb{R}^D \times \mathbb{R}^D) \) implies \( \delta \rho_{t_0} = 0 \), and thus \( \rho^1 = \rho^2 \).
\end{proof}

\begin{lemma}\label{equ:unique}
    	Assume that $ \rho^{1},\rho^{2}\in\MC{[0,T];\mathcal{P}_{2}\l(\mbRd\times\mbRd\r) }$ are two weak solutions to Equation \eqref{equ:FP_M_CBO} with the same initial data $ \rho_{0} $. Then it holds that 
	\begin{gather*}
		\sup_{t\in[0,T]}W_{2}\big(\rho^{1}_{t},\rho^{2}_{t}\big)=0,
	\end{gather*}
	where $ W_{2} $ is the 2-Wasserstein distance.
\end{lemma}

\begin{proof}
	Given $ \rho^1 $ and $ \rho^2 $, let us first consider the following two coupled linear SDEs:
	\[
	\begin{aligned}
		\intd\bar\Theta_t^{i} &= \left( \bar\Omega_t^{i} - \gamma_1\left( \bar\Theta_t^{i} - \CM [\mu^i_t] \right) \right)\intd t + \sigma(t) \intd W_{\theta,t}, \\
		\intd\bar\Omega_t^{i} &= -\left( m\left(  \bar\Theta_t^{i} - \CM[\mu^i_t] \right) + \gamma_2 \bar\Omega_t^{i} \right)\intd s + \sigma(t)  \sqrt{m} \intd W_{\omega,t},
	\end{aligned}
	\]
	for \( i=1,2 \), where \( ( \bar\Theta_t^{i},  \bar\Omega_t^{i}) \) are driven by independent Brownian motions \( W_{\theta,s} \) and \( W_{\omega,s} \), with the same initial condition \( (\Theta_0^{i,t}, \Omega_0^{i,t}) \sim \rho_0 \).

	Let \( \tilde{\rho}_t^i \) denote the law of \( (\Theta_t^{i}, \Omega_t^{i}) \) for \( i=1,2 \). By construction, the laws \( \tilde{\rho}_t^i \) satisfy the same Fokker-Planck equation:
	\[
	\begin{aligned}
		\partial_t \tilde{\rho}_t^i &= -\nabla_\theta \cdot \left[ \left( \omega - \gamma_1 (\theta - \CM[\mu^i_t]) \right) \tilde{\rho}_t^i \right] + \nabla_\omega \cdot \left[ \left( m (\theta - \CM[\mu^i_t]) + \gamma_2 \omega \right) \tilde{\rho}_t^i \right] \\
		&\quad + \frac{\sigma(t)^2 m}{2} \Delta_\omega \tilde{\rho}_t^i + \frac{\sigma(t)^2}{2} \Delta_\theta \tilde{\rho}_t^i,
	\end{aligned}
	\]
	in the weak sense, with initial condition \( \tilde{\rho}_0^i = \rho_0 \). Since both \( \rho^i_t \) solve this Fokker-Planck equation and we assumed \( \rho_0^1 = \rho_0^2 = \rho_0 \), by the uniqueness of solutions to this PDE in the last lemma, we have $\tilde \rho_t^i= \rho_t$ for $i=1,2$. As a result,  \( (\Theta_t^{1,t}, \Omega_t^{1,t}) \) and \( (\Theta_t^{2,t}, \Omega_t^{2,t}) \)  both solve Equation \eqref{equ:mean_M_CBO}. By Theorem \ref{thm:well_possness_mean_field}, we have
	\[
	\sup_{t \in [0, T]} \mathbb{E}\left[ \left| \Theta_t^{1,t} - \Theta_t^{2,t} \right|^2 + \left| \Omega_t^{1,t} - \Omega_t^{2,t} \right|^2 \right] = 0.
	\]
	This implies:
	\[
	\sup_{t \in [0, T]} W_2\left( \tilde{\rho}_t^1, \tilde{\rho}_t^2 \right) \leq  \sup_{t \in [0, T]} \mathbb{E}\left[ \left| \bar \Theta_t^{1,t} - \bar \Theta_t^{2,t} \right|^2 + \left|\bar \Omega_t^{1,t} - \bar \Omega_t^{2,t} \right|^2 \right]  = 0.
	\]
\end{proof}

\section{Global Convergence in the Mean Field Law} \label{app:converge}
\begin{lemma}\label{lem:energy_functional}
    Let $E[\rho_t]$ be the energy functional defined in \eqref{equ:energy_functional}. Under Assumption \ref{assum:Lip_J}, \[      \frac{\intd }{\intd t }E[\rho_t]  \leq  -\gamma E[\rho_t] +\lambda \sqrt{E[\rho_t] }\| \CM[\mu_t]-\tilde\theta\| +\frac{\sigma^2(t)D (m+1) }{2},     \] where $\gamma = \min\{\gamma_1,\gamma_2\}$ and $\lambda = \max\{m,\gamma_1\}$ are positive numbers.
\end{lemma}
\begin{proof}
    
From the definition of weak solution of Fokker-Planck equation \eqref{equ:FP_M_CBO}, we define 
$\phi(\theta,\omega) = \frac{1}{2}\|\theta-\tilde\theta\|^2 + \frac{1}{2m}\|\omega\|^2$, then
\begin{equation}
    \begin{aligned}
        \frac{\intd }{\intd t }E[\rho_t]  &=   \frac{\intd }{\intd t }\int \phi(\theta,\omega)  \intd \rho_t \\
        =& \int \langle \omega - \gamma_1(\theta - \CM[\mu_t]) ,\theta-\tilde\theta \rangle \intd \rho_t
        + \int \langle -m (\theta - \CM[\mu_t] )
- \gamma_2\omega ,\frac{1}{m}\omega  \rangle \intd \rho_t +\frac{\sigma^2 D  (m+1)}{2}\\
 = & \int \langle \omega  ,\CM[\mu_t]-\tilde\theta \rangle \intd \rho_t -\gamma_1 \int \langle \theta - \CM[\mu_t] ,\theta-\tilde\theta \rangle \intd \rho_t -\frac{\gamma_2}{m}\int \langle \omega ,\omega \rangle \intd \rho_t +\frac{\sigma^2D (m+1)}{2}\\
 = & \int \langle \omega  ,\CM [\mu_t]-\tilde\theta \rangle \intd \rho_t -\gamma_1 \int \langle \theta - \tilde\theta  ,\theta-\tilde\theta \rangle \intd \rho_t \\
 & - \gamma_1 \int \langle \tilde\theta -  \CM [\mu_t]  ,\theta-\tilde\theta \rangle \intd \rho_t-\frac{\gamma_2}{m}\int \langle \omega ,\omega \rangle \intd \rho_t +\frac{\sigma^2D(m+1)}{2}\\
 \leq&  \|\omega\|\| \CM[\mu_t]-\tilde\theta\| + \gamma_1 \|\theta-\tilde\theta\|\| \CM[\mu_t]-\tilde\theta\| -\gamma_1 \|\theta-\tilde\theta\|^2 -\frac{\gamma_2}{m}\|\omega\|^2  + \frac{\sigma^2D (m+1)}{2}\\
 \leq & -\gamma E[\rho_t] +\lambda \sqrt{E[\rho_t] }\| \CM[\mu_t]-\tilde\theta\| +\frac{\sigma^2D (m+1)}{2},
     \end{aligned}
\end{equation}
where in the last equility we take $\gamma = \min\{\gamma_1,\gamma_2\}$ and $\lambda = \max\{m,\gamma_1\}$.
\end{proof}

Next, we will show that $\|\CM[\mu_t]-\tilde\theta\|$ can be bounded by a suitable scalar multiple of $\sqrt{E[\rho_t]}$, we can apply Gronwall's inequality to bound the energy function.
\begin{lemma}
    Under Assumption~\ref{assum:Lip_J}, $\forall r>0$, we define $J_r :=\sup _{\theta\in B_{\theta,r}(\tilde\theta)} \mathcal J (\theta)$. Then $\forall r\in [0,R_0]$ and $q>0$ such that $q+ J_r - \underline{J} \leq \delta_J $, we have 
    \[
    \|\CM [\mu] - \tilde\theta\|\leq \frac{(q+J_r -\underline{J})}{\eta} + \frac{\exp(-\beta q )}{\rho(B_{\theta,r}(\tilde\theta))}\int\|\theta-\tilde\theta\| \intd \rho(\theta,\omega).
    \]
\end{lemma}
\begin{proof}
    Let $\tilde r = \frac{q+J_r -\underline{J}}{\eta} \geq \frac{J_r-\underline{J}}{\eta} \geq r $, we have
    \begin{equation*} \label{equ:laplace_1}
        \begin{aligned}
             \|\CM [\mu] - \tilde\theta\|  & \leq \int_{B_{\theta,\tilde r} (\tilde\theta)}  \|\theta-\tilde\theta\| \frac{w_\beta (\theta)}{\|w_\beta(\theta)\|_{L^1(\rho)}} \intd \rho + \int_{B^c_{\theta,\tilde r} (\tilde\theta)}  \|\theta-\tilde\theta\| \frac{w_\beta (\theta)}{\|w_\beta(\theta)\|_{L^1(\rho)}}  \intd \rho \\
             & \leq \tilde r  + \int_{B^c_{\theta,\tilde r} (\tilde\theta )}  \|\theta-\tilde\theta\| \frac{w_\beta (\theta)}{\|w_\beta(\theta)\|_{L^1(\rho)}} \intd \rho.
        \end{aligned}
    \end{equation*}  
By Markov's inequality, we have  $\|w_\beta\|_{L^1(\rho)}\geq a \rho (\{(\theta,\omega): \exp(-\beta \mathcal J(\theta)\geq a )\})$. By choosing $a = \exp (- \beta J_r)$, we have
\begin{equation*}
    \begin{aligned}
\|w_\beta\|_{L^1(\rho)} & \geq \exp(- \beta J_r) \rho \l ( \{(\theta,\omega): \exp(-\beta \mathcal J(\theta)\geq \exp(- \beta J_r) )\}\r )  \\
& =   \exp(- \beta J_r) \rho \l ( \{(\theta,\omega): J(\theta)\leq J_r  )\}\r )  \\
&\geq  \exp(- \beta J_r)  \rho (B_{\theta,r} (\tilde\theta) ) ,
    \end{aligned}
\end{equation*}
where the second inequality comes from the definition of $J_r$.
Thus for the second term in \eqref{equ:laplace_1}, we obtain 
\begin{equation*}
    \begin{aligned}
        \int_{B^c_{\theta,\tilde r} (\tilde\theta )}  \|\theta-\tilde\theta\| \frac{w_\beta (\theta)}{\|w_\beta(\theta)\|_{L^1(\rho)}}\intd \rho  &\leq 
    \frac{1}{\exp(- \beta J_r)  \rho (B_{\theta,r} (\tilde\theta) )} \int_{B^c_{\theta,\tilde r} (\tilde\theta )}  \|\theta-\tilde\theta\| w_\beta (\theta) \intd \rho \\
    & \leq  \frac{\exp(-\beta( \inf_{ B^c_{\theta,\tilde r} (\tilde\theta)}J(\theta) -J_r))}{ \rho (B_{\theta,r} (\tilde\theta )} \int_{B^c_{\theta,\tilde r} (\tilde\theta)}  \|\theta-\tilde\theta\|\intd \rho \\ 
    & \leq  \frac{\exp(-\beta( \inf_{ B^c_{\theta,\tilde r} (\tilde\theta )}J(\theta) -J_r))}{ \rho (B_{\theta,r} (\tilde\theta) )} \int  \|\theta-\tilde\theta\|\intd \rho.
    \end{aligned}
\end{equation*}
We also notice
\[
 \inf_{ B^c_{\theta,\tilde r} (\tilde\theta)} J(\theta) -J_r \geq \min\{\delta_J + \underline{J},\eta\tilde r + \underline{J}\} -J_r \geq \eta \tilde r - J_r +\underline J = q,
\]
where the first inequality comes from Assumption \ref{assum:Lip_J} and the second inequality comes from the definition of $\tilde r$ and $q$, $\tilde r = \frac{q + J_r -\underline{J}}{\eta } \leq \frac{\delta_J}{\eta}$. Combining the above inequality and the definition of $\tilde r$, we have 
\begin{equation*}
\begin{aligned}
     \|\CM [\mu] - \tilde\theta\| & \leq \frac{q+J_r-\underline{J}}{\eta} +  \frac{\exp(-\alpha( \inf_{ B^c_{\theta,\tilde r} (\tilde\theta )}J(\theta) -J_r) )}{ \rho (B_{\theta,r} (\tilde\theta) )}  \int  \|\theta-\tilde\theta\|\intd \rho \\ 
     & \leq \frac{q+J_r-\underline{J}}{\eta} +  \frac{\exp(-\beta q)}{ \rho (B_{\theta,r} (\tilde\theta) )}  \int  \|\theta-\tilde\theta\|\intd \rho. 
\end{aligned}
\end{equation*}
\end{proof}

Then we will establish a lower bound for $\rho_t(B_{\theta,r}(\tilde\theta))$. Notice that $  \rho_t(B_{\theta,r}(\tilde\theta)) = \rho_t\l ( \{(\theta,\omega):  \|\theta-\tilde\theta\|^2\leq r^2 \}\r )  \geq \rho_t\l ( \{(\theta,\omega): \|\theta-\tilde\theta\|^2+m^{-1} \|\omega\|^2\leq r^2\}\r )  := \rho_t(B_{r}(\tilde\theta,0)) $. We first define the mollifier $\phi_r (\theta,\omega)$ as follows
\begin{equation*}
    \phi_r (\theta ,\omega ) = \l \{\begin{aligned}
        & \exp\l ( 1-\frac{ r^2}{ r^ 2- (\|\theta-\tilde\theta\|^2 +m^{-1}\|\omega\|^2)}\r )   , & \text{if } \|\theta-\tilde\theta\|^2+m^{-1} \|\omega\|^2\leq   r^2 ,\\
        &0, & \text {else}.
    \end{aligned}\r  .
\end{equation*}
We have $\text{Im}(\phi_r) \in [0,1]$, $ \phi_r \in \mathcal C^{\infty}_c  $.
First, we compute the first-order and second-order derivatives as
\begin{equation*}
    \begin{aligned}
        \nabla_\theta \phi_r &= -2   r^2 \frac{\theta-\tilde\theta}{\l (   r^ 2- (\|\theta-\tilde\theta\|^2 +m^{-1}\|\omega\|^2)\r ) ^2}\phi_r,\\
        \nabla_\omega \phi_r &= -2   r^2 \frac{m^{-1}\omega}{\l (   r^ 2- (\|\theta-\tilde\theta\|^2 +m^{-1}\|\omega\|^2)\r ) ^2}\phi_r,
    \end{aligned}
\end{equation*}
and

\begin{equation*}
    \begin{aligned}
        \Delta_\theta \phi_r =&  -2  r^2 \frac{
        D\l ( r^ 2- (\|\theta-\tilde\theta\|^2 +m^{-1}\|\omega\|^2)\r ) ^2  - 2 \l (   r^ 2- (\|\theta-\tilde\theta\|^2 +m^{-1}\|\omega\|^2)\r )  \l ( -2(\theta - \tilde\theta)\cdot (\theta-\tilde\theta) \r )  }{\l (   r^ 2- (\|\theta-\tilde\theta\|^2 +m^{-1}\|\omega\|^2)\r ) ^4}\phi_r\\
        & +  4   r^4 \frac{\|\theta-\tilde\theta\|^2}{\l (  r^ 2- (\|\theta-\tilde\theta\|^2 +m^{-1}\|\omega\|^2)\r ) ^4}\phi_r
        \\
        = &  2   r^2 \frac{\|\theta-\tilde\theta\|^2(-2 r^2 +4\|\theta-\tilde\theta\|^2 + 4m^{-1}\|\omega\|^2  ) -D\l (   r^ 2- (\|\theta-\tilde\theta\|^2 +m^{-1}\|\omega\|^2)\r ) ^2  }{\l (   r^ 2- (\|\theta-\tilde\theta\|^2 +m^{-1}\|\omega\|^2)\r ) ^4}\phi_r,\\
\Delta_\omega \phi_r =& -2   r^2 \frac{D\l (   r^ 2- (\|\theta-\tilde\theta\|^2 +m^{-1}\|\omega\|^2)\r ) ^2 m^{-1} -2 \l (   r^ 2- (\|\theta-\tilde\theta\|^2 +m^{-1}\|\omega\|^2)\r )  (-m^{-1}\omega)\cdot(m^{-1}\omega)  }{\l (   r^ 2- (\|\theta-\tilde\theta\|^2 +m^{-1}\|\omega\|^2)\r ) ^4}\phi_r \\
         & +  4   r^4 \frac{m^{-2}\omega^2}{\l (   r^ 2- (\|\theta-\tilde\theta\|^2 +m^{-1}\|\omega\|^2)\r ) ^4}\phi_r\\
 = &  2   r^2 \frac{m^{-2}\|\omega\|^2( -2 
  r^ 2 +4 \|\theta-\tilde\theta\|^2 +4m^{-1}\|\omega\|^2 )  - D \l (   r^ 2- (\|\theta-\tilde\theta\|^2 +m^{-1}\|\omega\|^2)\r ) ^2 m^{-1}}{\l (  r^ 2- (\|\theta-\tilde\theta\|^2 +m^{-1}\|\omega\|^2)\r ) ^4}\phi_r .
    \end{aligned}
\end{equation*}
\begin{lemma}
    Let $T>0,r>0$. Choose parameters $\overline{\sigma}\geq\sigma(t)\geq\underline{\sigma}>0$. Assume $\rho\in \mathcal C([0,T],\mathcal P(\mathbb R^{2D}))$ weakly solves the Fokker-Planck equation \eqref{equ:FP_M_CBO} with initial condition $\rho_0$. Then, $\forall t\in [0,T]$, we have 
    \begin{equation*}
        \begin{aligned}
            \rho_t(B_{\theta,r}(\tilde\theta)) \geq \l ( \int \phi_r(\theta,\omega) \intd \rho_0 (\theta,\omega)\r ) \exp(-p t),
        \end{aligned}
    \end{equation*}
    where 
    \begin{equation*}
        p:= \max\l \{\frac{4 \lambda ( \sqrt{k} r  + B )\sqrt{k}}{\l ( 1 - k \r ) ^2r} + \frac{2\overline{\sigma}^2(k+D)}{(1- k )^4r^2} ,   \frac{8(B+r)^2 \lambda ^2 }{(2k-1)\underline{\sigma}^2}\r  \},
    \end{equation*}
    for any $B>0$ with $\sup_{t\in[0,T]} \|\CM[\mu_t] - \tilde\theta\|\leq B $ and for any $k\in (\frac{1 }{2},1 )$ satisfying 
    $ (-1+2 k )k \geq 2 D  (1 -k)^2$ and $\lambda = \max\{\gamma_1,m\}$.
\end{lemma}
\begin{proof}
By the properties of the mollifier $\phi_r$, we have $\mu_t(B_{\theta,r}(\tilde\theta)) \geq\rho_t(B_{r}(\tilde\theta,0)) \geq \int \phi_r(\theta,\omega)\intd \rho_t(\theta,\omega)$.  Using properties of the weak solution $\rho_t$, we have 
\begin{equation*}
    \begin{aligned}
\frac{\intd }{\intd t } \int \phi_r(\theta,\omega) \intd \rho_t(\theta,\omega) =& \int \l \langle\omega-\gamma_1(\theta-\CM[\mu_t]) , -2   r^2 \frac{\theta-\tilde\theta}{\l (  r^ 2- (\|\theta-\tilde\theta\|^2 + \|\omega\|^2)\r ) ^2}\phi_r\r  \rangle \intd \rho_t\\
& +  \int \l \langle -m(\theta-\CM[\mu_t]) -\gamma_2\omega ,  -2   r^2 \frac{m^{-1}\omega}{\l (   r^ 2- (\|\theta-\tilde\theta\|^2 + \|\omega\|^2)\r ) ^2}\phi_r\r  \rangle \intd \rho_t \\
& +\frac{\sigma^2}{2}\int \l ( m \Delta_\omega\phi_r+ \Delta_\theta \phi_r\r ) \intd \rho_t \\
 = &  \int \l ( \gamma \langle \theta-\CM[\mu_t], \theta-\tilde\theta  \rangle + \l \langle\gamma_2\omega,m^{-1}\omega\r  \rangle + \l \langle \omega, \tilde\theta -  \CM[\mu_t]\r  \rangle \r )  \\
& \frac{2\ r^2} {\l (  r^ 2- (\|\theta-\tilde\theta\|^2 +m^{-1}\|\omega\|^2)\r ) ^2}\phi_r\intd \rho_t 
 +\frac{\sigma^2}{2}\int \l ( m\Delta_\omega\phi_r+ \Delta_\theta\phi_r\r ) \intd \rho_t \\ 
 : = & \int T_1(\theta,\omega)\intd \rho_t+\int T_2(\theta,\omega)\intd \rho_t .
        \end{aligned}
    \end{equation*}
Since $\phi_r$ vanishes outside of $D_r := \{(\theta,\omega): \|\theta-\tilde\theta\|^2 + m^{-1}\|\omega\|^2\leq  r^2\}$, we restrict our attention to the open ball $D_r$.
To obtain the lower bound, we introduce the following subsets 
\begin{equation}
    \begin{aligned}
        K_1 &:= \l \{(\theta,\omega): \|\theta-\tilde\theta\|^2 + m^{-1}\|\omega\|^2 >  k r^2  \r  \},\\
        K_2 &:= \l \{(\theta,\omega): -\l (  \gamma_1 \langle\theta - \CM [\mu_t], \theta-\tilde\theta \rangle  +  \langle\gamma_2\omega,m^{-1}\omega\rangle 
 + \langle\omega , \tilde\theta-\mathcal M _\beta[\rho_t]\rangle \r )  (  r^2-\|\theta-\tilde\theta\| -m^{-1}\|\omega\|^2)^2  
        \r  .\\
        &\l . > \tilde k  \frac{\sigma^2}{2}r^2(\|\theta-\tilde\theta\|^2 +m^{-1}\|\omega\|^2) \r  \},
    \end{aligned}
\end{equation}
where $\tilde k = 2 k -1\in (0,1)$.
We divide the integral region into three domains.\\
\textbf{Domain $\Omega_r\cap K_1^c$}: 
We have $\|\theta-\tilde\theta\|^2 + m^{-1}\|\omega\|^2 \leq k r^2  $ in this domain and we can get 
\begin{equation*}
    \begin{aligned}
      T_1(\theta,\omega) = &\l ( \langle \gamma_1(\theta-\CM[\mu_t], \theta-\tilde\theta\rangle + \langle\gamma_2\omega,m^{-1}\omega\rangle + \langle \omega, \tilde\theta -  \CM[\mu_t]\rangle \r )  \frac{2  r^2} {\l (   r^ 2- (\|\theta-\tilde\theta\|^2 +m^{-1}\|\omega\|^2)\r ) ^2}\phi_r\\
      \geq & -\l ( \gamma_1 \| \theta-\CM[\mu_t] \|
      \|\theta-\tilde\theta\| + \|\omega\|\|\tilde \theta- \CM [\mu_t]\|\r ) \frac{2  r^2} {\l (  r^ 2- (\|\theta-\tilde\theta\|^2 +m^{-1}\|\omega\|^2)\r ) ^2}\phi_r\\
      \geq &  - \l ( \gamma_1  \|\theta-\tilde\theta\|  +  \|\omega\|\r ) \frac{2 ( \sqrt{k} r  + B ) r^2} {\l (  r^ 2- (\|\theta-\tilde\theta\|^2 +m^{-1}\|\omega\|^2)\r ) ^2}\phi_r\\
      \geq &  - \l ( \gamma_1  \|\theta-\tilde\theta\|  +  \|\omega\|\r ) \frac{2 ( \sqrt{k} r + B ) r^2} {\l (   r^ 2- k r^2\r ) ^2}\phi_r\\
      \geq & -  \lambda \l (  \|\theta-\tilde\theta\|  +  m^{-1}\|\omega\|\r ) \frac{2 ( \sqrt{k} r + B ) } {\l ( 1 - k \r ) ^2r^2}\phi_r\\
      \geq & - 2\lambda\sqrt{k}r  \frac{2 ( \sqrt{k} r  + B ) } {\l ( 1 - k \r ) ^2r^2}\phi_r  =   -  \frac{4 \lambda ( \sqrt{k} r  + B )\sqrt{k}}{\l ( 1 - k \r ) ^2r}\phi_r  := -p_1 \phi_r,
    \end{aligned}
\end{equation*}
where the first inequality comes from Cauchy-Schwarz inequality and the positiveness of $\|\omega\|^2$, the second inequality comes from the boundedness of $ \|\theta- \CM [\mu_t]\|\leq \|\theta-\tilde\theta\| + \|\tilde\theta - \CM [\mu_t]\|\leq \sqrt{k} r 
 + B $, the third inequality comes from the defintion of domain $K_1^c$, and the fourth inequality comes from the definition of $\lambda$.

For $T_2$, we have 
\begin{equation}
    \begin{aligned}
      T_2(\theta,\omega) =& \sigma^2   r^2 \frac{ (\|\theta-\tilde\theta\|^2  + m^{-1}\|\omega\|^2)( -2 
  r^ 2 +4 \|\theta-\tilde\theta\|^2 +4m^{-1}\|\omega\|^2 )}{\l ( r^ 2- (\|\theta-\tilde\theta\|^2 +m^{-1}\|\omega\|^2)\r ) ^4}\phi_r \\
  &- \sigma^2   r^2\frac{ 2D}{\l ( r^ 2- (\|\theta-\tilde\theta\|^2 +m^{-1}\|\omega\|^2)\r ) ^2}\phi_r \\
  \geq & -\sigma^2   r^2 \frac{ 2 r^ 2(\|\theta-\tilde\theta\|^2  + m^{-1}\|\omega\|^2) }{\l ( r^ 2- (\|\theta-\tilde\theta\|^2 +m^{-1}\|\omega\|^2)\r ) ^4}\phi_r \\
  & - \sigma^2   r^2\frac{ 2D }{\l ( r^ 2- (\|\theta-\tilde\theta\|^2 +m^{-1}\|\omega\|^2)\r ) ^2}\phi_r \\
  \geq & -\sigma^2   r^2 \frac{ 2 k r^ 4 }{\l ( r^ 2- k r^2\r ) ^4}\phi_r - \sigma^2   r^2\frac{ 2D}{\l ( r^ 2- kr^2\r ) ^2}\phi_r \\
  = & -\sigma^2    \frac{ 2 k }{(1- k )^4r^2}\phi_r - \sigma^2   \frac{ 2D}{\l ( 1- k\r ) ^2r^2}\phi_r\\
  \geq & - \frac{2\overline{\sigma}^2(k+D)}{(1- k )^4r^2}\phi_r := -p_2 \phi_r,
    \end{aligned}
\end{equation}
where the first inequlaity uses the positiveness of $\|\theta-\tilde\theta\|^2$ and $\|\omega\|^2$, the second inequality uses the properties of $K_1^c$, and the third inequality use $1-k\in (0,\frac{1}{2})$.\\
\textbf{Domain  $\Omega_r\cap K_1\cap K_2^c$} 

In this domian, we have $\|\theta-\tilde\theta\|^2 + m^{-1}\|\omega\|^2 > k r^2  $ and 
\begin{equation}\label{equ:K2c}
    -\l (  \gamma_1 \langle\theta - \mathcal M_\beta [\rho_t], \theta-\tilde\theta \rangle + \langle\gamma_2\omega,m^{-1}\omega\rangle  + \langle\omega , \tilde\theta-\mathcal M _\beta[\rho_t]\rangle \r )  (  r^2-\|\theta-\tilde\theta\| -m^{-1}\|\omega\|^2)^2  \leq \tilde k  \frac{\sigma^2}{2}r^2(\|\theta-\tilde\theta\|^2 +m^{-1}\|\omega\|^2).
\end{equation}

Our goal is to show $T_1(\theta,\omega)+T_2(\theta,\omega)\geq 0$ in this subset. We first compute 
\begin{equation*}
     \begin{aligned}
      &\frac{T_1(\theta,\omega)+T_2(\theta,\omega)}{2 r^2 \phi_r} \l (  r^ 2- (\|\theta-\tilde\theta\|^2 +m^{-1}\|\omega\|^2)\r ) ^4\\
      & =  \l ( \langle \gamma(\theta-\CM[\mu_t], \theta-\tilde\theta\rangle + \langle\gamma_2\omega,m^{-1}\omega\rangle + \langle \omega, \tilde\theta -  \CM[\mu_t]\rangle \r )  \l (   r^ 2- (\|\theta-\tilde\theta\|^2 +m^{-1}\|\omega\|^2)\r ) ^2 \\
      & +\sigma^2   m^{-1}\|\omega\|^2(-  r^2 +2\|\theta-\tilde\theta\|^2 + 2m^{-1}\|\omega\|^2  ) \\
      & -\frac{\sigma^2 D}{2} \l (  r^ 2- (\|\theta-\tilde\theta\|^2 +m^{-1}\|\omega\|^2)\r ) ^2 \\
      & + \sigma^2 \|\theta-\tilde\theta\|^2(- r^2 +2\|\theta-\tilde\theta\|^2 + 2m^{-1}\|\omega\|^2  ) \\
      & - \frac{\sigma^2 D}{2} \l (   r^ 2- (\|\theta-\tilde\theta\|^2 +m^{-1}\|\omega\|^2)\r ) ^2 .
      \end{aligned}
\end{equation*}
To prove the positiveness, we need to prove
\begin{equation*}
    \begin{aligned}
    &-\l ( \langle \gamma(\theta-\CM[\mu_t], \theta-\tilde\theta\rangle + \langle\gamma_2\omega,m^{-1}\omega\rangle + \langle \omega, \tilde\theta -  \CM[\mu_t]\rangle \r )  \l (   r^ 2- (\|\theta -\tilde\theta \|^2 +m^{-1}\|\omega\|^2)\r ) ^2 \\
      & +\sigma^2 D  \l (   r^ 2- (\|\theta-\tilde\theta\|^2 +m^{-1}\|\omega\|^2)\r ) ^2 \\
      & \leq \sigma^2    (m^{-1} \|\omega\|^2+\|\theta-\tilde\theta\|^2)(-  r^2 +2\|\theta-\tilde\theta\|^2 + 2m^{-1}\|\omega\|^2  )   .  \end{aligned}
\end{equation*}
For the first term, we have 
\begin{equation*}
\begin{aligned}
       & -\l ( \langle \gamma(\theta-\CM[\mu_t], \theta-\tilde\theta\rangle + \langle\gamma_2\omega,m^{-1}\omega\rangle + \langle \omega, \tilde\theta -  \CM[\mu_t]\rangle \r )  \l (   r^ 2- (\|\theta-\tilde\theta\|^2 +m^{-1}\|\omega\|^2)\r ) ^2 \\
    \leq &  \tilde k \frac{\sigma^2}{2} r^2(\|\theta-\tilde\theta\|^2 +m^{-1}\|\omega\|^2]) \\
     = & (2k-1 )\frac{\sigma^2}{2} r^2(\|\theta-\tilde\theta\|^2 +m^{-1}\|\omega\|^2) \\
    \leq &  ( 2 \|\theta-\tilde\theta\|^2+ 2 m^{-1}\|\omega\|^2  -   r^2)\frac{\sigma^2}{2} (\|\theta-\tilde\theta\|^2 +m^{-1}\|\omega\|^2 ),
\end{aligned}
\end{equation*}
where the first inequality comes from the positiveness of the norm and the second inequality comes from Equation \eqref{equ:K2c}. By the definition $\tilde k = 2k-1$, we have the equality in the fourth line and the last inequality comes from $kr^2\leq \|\theta-\tilde\theta\|^2+ m^{-1}\|\omega\|^2 $. 
For the second term, we have 
\begin{equation*}
\begin{aligned}
       &  \sigma^2 D  \l (   r^ 2- (\|\theta-\tilde\theta\|^2 +m^{-1}\|\omega\|^2)\r ) ^2  \\
    &\leq   \sigma^2  D   \l (  1 - k \r ) ^2 r^4    \\ 
    & \leq \frac{\sigma^2}{2}  (2k-1  )r^2 kr^2\\
    & \leq \frac{\sigma^2}{2} ( 2 \|\theta-\tilde\theta\|^2+ 2 m^{-1}\|\omega\|^2  - r^2)(\|\theta-\tilde\theta\|^2 +m^{-1}\|\omega\|^2 ),
\end{aligned}
\end{equation*}
where in the second inequality we use $(-1+2 k )k \geq 2D (1 -k)^2$.
Hence, we have the positiveness of $T_1(\theta,\omega) +T_2(\theta,\omega)$.

\textbf{Domain  $\Omega_r\cap K_1\cap K_2$} 

In this subset, we have $\|\theta-\tilde\theta\|^2 + m^{-1}\|\omega\|^2 > k r^2  $  and 
\begin{equation}
    \begin{aligned}\label{equ:k2}
- \l ( \gamma_1\langle\theta-\CM [\mu_t], \theta  - \tilde\theta \rangle +\langle\gamma_2\omega,m^{-1}\omega\rangle + \langle\omega , \tilde\theta-\mathcal M _\beta[\rho_t]\rangle \r )  ( r^2-\|\theta-\tilde\theta\| -m^{-1}\|\omega\|^2)^2  
         > \tilde k  \frac{\sigma^2}{2}r^2(\|\theta-\tilde\theta\|^2 +m^{-1}\|\omega\|^2).
    \end{aligned}
\end{equation}
In this subset, we have
\begin{equation} \label{equ:k2_ineq1}
    \begin{aligned}
        & \gamma_1   \langle \theta-\CM[\mu_t], \theta-\tilde\theta\rangle + \langle\gamma_2\omega,m^{-1}\omega\rangle + \langle \omega, \tilde\theta -  \CM[\mu_t]\rangle   \\
         \geq  & \gamma_1 \langle  \theta-\CM[\mu_t], \theta-\tilde\theta\rangle   + \langle \omega, \tilde\theta -  \CM[\mu_t]\rangle  \\
         \geq &  -  \gamma_1\| \theta-\CM[\mu_t]\|\| \theta-\tilde\theta\|  - \|\omega\| \|\tilde\theta -  \CM[\mu_t]\|    \\
         \geq  & - \lambda  (\|\tilde \theta-\CM[\mu_t]\|+\|\theta-\tilde \theta\|)( \| \theta-\tilde\theta\|  + m^{-1}\|\omega\| ),
    \end{aligned}
\end{equation}
where the last inequality comes from the definition of $\lambda$.
From the inequality \eqref{equ:k2}, we have
\begin{equation} \label{equ:k2_ineq2}
    \begin{aligned}
       \frac{( \|\theta-\tilde\theta\|  +m^{-1}\|\omega\|)^2 }{( r^2-\|\theta-\tilde\theta\|^2 -m^{-1}\|\omega\|^2)^2  } 
       & \leq  2 \frac{\|\theta-\tilde\theta\|^2  + m^{-1}\|\omega\|^2}{(  r^2-\|\theta-\tilde\theta\|^2 -m^{-1}\|\omega\|^2)^2  } \\
       & <  -\frac{4}{\tilde k\sigma^2 r^2} \l ( \gamma_1 \langle\theta-\mathcal M_\beta [\rho_t], \theta-\tilde\theta \rangle +\langle\gamma_2\omega,m^{-1}\omega\rangle  + \langle\omega , \tilde\theta-\mathcal M _\beta[\rho_t]\rangle \r )  .
    \end{aligned}
\end{equation}
Then we are ready to  prove
\begin{equation*}
    \begin{aligned}
        & \frac{ \langle \gamma_1 (\theta-\CM[\mu]), \theta-\tilde\theta\rangle + \langle\gamma_2\omega,m^{-1}\omega\rangle + \langle \omega, \tilde\theta -  \CM[\mu_t]\rangle   } {\l (   r^ 2- (\|\theta-\tilde\theta\|^2 +m^{-1}\|\omega\|^2)\r ) ^2} \\
        & \geq - \lambda \frac{ (\|\tilde \theta-\CM[\mu_t]\|+\|\tilde \theta -\theta \| )( \| \theta-\tilde\theta\|  + m^{-1}\|\omega\| ) } {\l (   r^ 2- (\|\theta-\tilde\theta\|^2 +m^{-1}\|\omega\|^2)\r ) ^2} \\
        & \geq \lambda \frac{4}{\tilde k\sigma^2 r^2}  \frac{\l ( \gamma_1 \langle\theta-\mathcal M_\beta [\rho_t], \theta-\tilde\theta \rangle +\langle\gamma_2\omega,m^{-1}\omega\rangle  + \langle\omega , \tilde\theta-\mathcal M _\beta[\rho_t]\rangle \r )  (\|\tilde \theta-\CM[\mu_t]\|+\|\tilde \theta -\theta \| ) }{ \| \theta-\tilde\theta\|  + m^{-1}\|\omega\| }\\
        & \geq - \lambda^2 \frac{4}{\tilde k\sigma^2 r^2}  \frac{(\|\tilde \theta-\CM[\mu_t]\|+\|\tilde \theta -\theta \| )^2 \l ( \|\theta-\tilde\theta \|  + m^{-1}\|\omega \| \r )  }{ \| \theta-\tilde\theta\|  + m^{-1}\|\omega\| }\\
        & \geq - \lambda^2 \frac{4}{\tilde k\sigma^2 r^2}  (\|\tilde \theta-\CM[\mu_t]\|+\|\tilde \theta -\theta \| )^2\\
        & \geq - \lambda^2 \frac{4(B+r)^2}{\tilde k\sigma^2 r^2} a,
    \end{aligned}
\end{equation*}
where the first and third inequalities are derived from the inequality \eqref{equ:k2_ineq1}, and the second one is a consequence of \eqref{equ:k2_ineq2}. Utilizing Cauchy–Schwarz inequality and the definition specified in $\lambda$, we have the third and fourth inequalities.
Given this we have 
\begin{equation*}
    \begin{aligned}
        T_1(\theta,\omega) =&  \frac{ \langle \gamma_1 \theta-\CM[\mu_t], \theta-\tilde\theta\rangle + \langle\gamma_2\omega,m^{-1}\omega\rangle + \langle \omega, \tilde\theta -  \CM[\mu_t]\rangle   } {\l (  r^ 2- (\|\theta-\tilde\theta\|^2 +m^{-1}\|\omega\|^2)\r ) ^2} 2  r^2 \phi_r \\
         & \geq - \frac{8(B+r)^2\lambda ^2 }{\tilde k \underline{\sigma}^2} \phi_r  = - \frac{8(B+r)^2 \lambda ^2 }{(2k-1)\underline{\sigma}^2} \phi_r := -p_3 \phi_r.
    \end{aligned}
\end{equation*}
For  $T_2$, it is positive whenever
\[
(\|\theta-\tilde\theta\|^2+m^{-1}\|\omega\|^2)(-2  r^2 +4\|\theta-\tilde\theta\| + 4m^{-1}\|\omega\|^2  )\geq 2 D \l (   r^ 2- (\|\theta-\tilde\theta\|^2 +m^{-1}\|\omega\|^2)\r ) ^2 , 
\]
we have 
\begin{equation*}
    \begin{aligned}
        &(\|\theta-\tilde\theta\|^2+m^{-1}\|\omega\|^2)(-2  r^2 +4\|\theta-\tilde\theta\| + 4m^{-1}\|\omega\|^2  ) \\
        \geq &   (\|\theta-\tilde\theta\|^2+ m^{-1}\|\omega\|^2) (-2  r^2 + 4 k r^2 )\\
        \geq & (\|\theta-\tilde\theta\|^2+ m^{-1}\|\omega\|^2)(-1 +  2k  )2r^2 \\
        \geq & (\|\theta-\tilde\theta\|^2+ m^{-1}\|\omega\|^2) 2D \frac{(1 -k)^2}{k} 2r^2\\
        \geq &   r^2 2D(1  -k)^2 2r^2 \\
         =   & 2D ( r^2 - k r^2)^2 \\
        \geq &  2D ( r^2 -  (\|\theta-\tilde\theta\|^2 +m^{-1}\|\omega\|^2) )^2 .
    \end{aligned}
\end{equation*}
This is safiesifed for all $ \|\theta-\tilde\theta\|^2 + m^{-1}\|\omega\|^2 \geq k r^2  $.

\textbf{Concluding the proof:} Using the evolution of $\phi_r$, we now get 
\begin{equation*}
    \begin{aligned}
        & \frac{\intd }{\intd t } \int \phi_r(\theta,\omega) \intd \rho_t(\theta,\omega) =  \int_{K_1\cap K_2^c\cap\Omega_r}  T_1(\theta,\omega)+T_2(\theta,\omega) \intd \rho_t(\theta,\omega) \\
         &+ \int_{K_1\cap K_2\cap\Omega_r}  T_1(\theta,\omega)+T_2(\theta,\omega) \intd \rho_t(\theta,\omega) +\int_{K_1^c\cap\Omega_r}  T_1(\theta,\omega)+T_2(\theta,\omega) \intd \rho_t(\theta,\omega)\\
        \geq&    -\max\{p_1+p_2,p_3\} \int\phi_r(\theta,\omega)\intd \rho_t(\theta,\omega)  = -p \int\phi_r(\theta,\omega)\intd \rho_t(\theta,\omega).
    \end{aligned}
\end{equation*}
\end{proof}


\begin{proof}[Proof of Theorem \ref{thm:converge}]
    We choose parameters $\beta$ such that 
    \begin{equation*}
        \beta > \beta_0 := \frac{1}{q_\epsilon} \l ( \log\l ( \frac{4\sqrt{2E[\rho_0]}}{\mathrm c(\tau,\lambda) \sqrt{\epsilon}} + \frac{p}{(1-\tau)\lambda}\log\l ( \frac{E[\rho_0]}{\epsilon}\r )  - \log\rho_0(B_{\frac{r_\epsilon}{2}} (\tilde\theta,0))\r )  \r ) ,
    \end{equation*} 
    where we introduce 
    \[\mathrm c(\tau,\lambda) =  \frac{\tau\gamma}{\lambda}, \quad 
 q_\epsilon =  \frac{1}{2}\min\{\frac{c(\tau,\lambda)\sqrt{\epsilon}\eta}{2}, \delta_J \},\text{ and } r_\epsilon = \max_{x\in [0,R_0]} \{\max_{(\theta,\omega)\in B_{s}(\tilde\theta,0)} J(\theta) \leq q_\epsilon+ \underline{J} \}, \]  
    and define the time horizon $T_\beta \geq 0$, which may depend on $\beta$, by 
    \begin{equation*}
        T_\beta = \sup\{t\geq 0 :  E[\mu_{t'}] > \epsilon \text{ and  } \|\mathcal M_\beta [\mu_{t'} ] -\tilde\theta\| < C(t') \text{ for all }t' \in [0,t] \}
    \end{equation*}
    with $C(t) = \mathrm c(\tau,\lambda ) \sqrt{E(\rho_t)}$.
    First we want to prove $T_\beta>0$, which follows from the continunity of the mappings $t\to E[\rho_t]$ and $t\to \|\mathcal M_\beta [\mu_t]-\tilde\theta\|$ since $E[\rho_0]>0$ and $\|\mathcal M_\beta [\mu_{0} ] -\tilde\theta\| < C(0)$. While the former holds by assumption, the latter follows by 
    \begin{equation*}
    \begin{aligned}
        \|\mathcal M_\beta [\mu_{0} ] -\tilde\theta\| &\leq \frac{(q_\epsilon+J_{r_\epsilon} -\underline{J})}{\eta} + \frac{\exp(-\beta q_\epsilon )}{\rho(B_{\theta,r_\epsilon}(\tilde\theta))}\int\|\theta-\tilde\theta\| \intd \rho_0(\theta,\omega)\\
        &\leq \frac{(q_\epsilon+J_{r_\epsilon} -\underline{J})}{\eta} + \frac{\exp(-\beta q_\epsilon )}{\rho(B_{r_\epsilon}(\tilde\theta,0))}\int\|\theta-\tilde\theta\| \intd \rho_0(\theta,\omega)\\
         & \leq \frac{\mathrm c(\tau,\lambda)\sqrt{\epsilon}}{2}   + \frac{\exp(-\beta q_\epsilon )}{\rho(B_{r_\epsilon}(\tilde\theta,0))}\sqrt{2E[\rho_0]}\\
         &  \leq \mathrm c(\tau,\lambda) \sqrt{\epsilon}  \leq \mathrm c(\tau,\lambda) \sqrt{E[\rho_0]} = C(0), 
    \end{aligned}
    \end{equation*}
    where we use the definition of $\beta$ in the first inequality of the last line.
    Recall the Lemma \ref{lem:energy_functional}, up to time $T_\beta$ 
 \begin{equation*}
     \begin{aligned}
         \frac{\intd }{\intd t }E[\rho_t]  &\leq  -\gamma E[\rho_t] +\lambda \sqrt{E[\rho_t] }\| \CM[\mu_t]-\tilde\theta\| +\frac{\sigma^2(t)D (m+1) }{2}\\
         &\leq -\l ( 1-\tau\r ) \gamma E[\rho_t]+\frac{\sigma^2(t)D (m+1) }{2}.
     \end{aligned}
 \end{equation*} 
 Thus we have 
 \begin{equation*}
 \begin{aligned}
     \frac{\intd}{\intd t} \l ( \exp\l (    (1-\tau )\gamma t  \r ) E[\rho_t] \r ) = &  (1-\tau )\gamma\l ( \exp\l (    (1-\tau )\gamma t  \r ) E[\rho_t] \r )   + \exp\l (    (1-\tau )\gamma t  \r ) \frac{\intd}{\intd t} E[\rho_t]  \\
     \leq  & \exp\l (    (1-\tau )\gamma t  \r ) \frac{\sigma^2(t)D(m+1) }{2}.
 \end{aligned}
 \end{equation*}
 Therefore we have 
 \begin{equation*}
     \begin{aligned}
          \l ( \exp\l (    (1-\tau )\gamma t  \r ) E[\rho_t] \r )  -E[\rho_0]  \leq& \int_0^t \exp\l ( (1-\tau ) \lambda s\r )  \sigma^2(s)\intd s  \\
          =&  \frac{\sigma_1^2(1-\exp \l ( (-2\sigma_2 + \lambda(1-\tau))t\r )  }{
  2 \sigma_2  - \lambda(1 - \tau)}.
     \end{aligned}
 \end{equation*}  
We can get the boundedness for $E[\rho_t]$, for $2\sigma_2 -\lambda(1-\tau)< 0$ by the chosen of $\tau$ and $\lambda$, then we have 
 \begin{equation*}
     \begin{aligned}
         E[\rho_t] \leq \exp\l ( - (1- \tau  )t  \lambda \r ) E[\rho_0] .
     \end{aligned}
 \end{equation*} 
 Accordingly, we note that $E(\rho_t)$ is decreasing in $t$, which implies the decay of the function $C(t)$ as well. Hence, recalling the definition of $T_\beta$, we may bound $\max _{t\in [0,T_\beta]}\| \mathcal M_\beta [\rho_{t'} ] -\tilde\theta  \|\leq \max _{t\in [0,T_\beta]} C(t)\leq C(0) $. We now conclude by showing $\min _{t\in [0,T_\beta]} E(\rho_t) \leq \epsilon$ with $T_\beta\leq T^*$. For this, we distinguish the following three cases.\\
 \textbf{Case $T_\beta\geq T^*$}: If $T_\beta\geq T^*$, we can use the definition of $T^* = \frac{1}{(1-\tau)\lambda} \log(\frac{E[\rho_0]}{\epsilon})$ and the time evolution bound of $E[\rho_t]$ to conclude that $E[\rho_{T^*}]\leq \epsilon$. Hence, by definition of $T_\beta$, we find $E[\rho_{T_\beta}]\leq \epsilon$ and $T_\beta  = T^*$.\\
 \textbf{Case $T_\beta< T^*$ and $E[\rho_{T_\beta}] \leq \epsilon$:} Nothing need to discussed in this case.\\
 \textbf{Case $T_\beta<T^*$ and$E[\rho_{T_\beta}] > \epsilon$:} We shall prove that this case will never occur.
 \begin{equation*}
     \begin{aligned}
         \|\mathcal M _\beta[\mu_{T_\beta}]-\tilde\theta \| & \leq \frac{(q_\epsilon+J_{r_\epsilon} -\underline{J})}{\eta} + \frac{\exp(-\beta q_\epsilon )}{\rho(B_{\theta,r_\epsilon}(\tilde\theta))}\int\|\theta-\tilde\theta\| \intd \rho_{T_\beta}(\theta,\omega)\\
         & < \frac{\mathrm c(\tau,\lambda)\sqrt{E[\rho_{T_\beta}]}}{2} + \frac{\exp(-\beta q_\epsilon )}{\rho(B_{\theta,r_\epsilon}(\tilde\theta))}\sqrt{E[\mu_{T_\beta}]}.
     \end{aligned}
 \end{equation*}
    \textcolor{red}{Since, we have $ \max _{t\in [0,T_\beta]}\| \mathcal M_\beta [\mu_{t'} ] -\tilde\theta  \| = B = C(0)$ guarantees that there exist a $p>0$ with 
    \begin{equation}
        \rho_{T_\beta}(B_{\theta,r_\epsilon}(\tilde\theta)) \geq \l ( \int \phi_{r_\epsilon}(\theta,\omega) \intd \rho_0 (\theta,\omega)\r ) \exp(-p T_{\beta}) \geq \frac{1}{2} \rho_0 \l ( B_{\frac{r_\epsilon}{2}}(\tilde\theta,0)) \r ) \exp(-p T^*),
    \end{equation}
    where we used $(\tilde\theta,0) \in supp (\rho_0)$ for bounding the initial mass $\rho_0$ and the fact that $\phi_r$ is bounded from below on $B_{\frac{r_\epsilon}{2}}(\tilde\theta,0)) $  by $1/2$.} With this, we can conclude that 
    \begin{equation*}
     \begin{aligned}
         \|\mathcal M _\beta[\mu_{T_\beta}]-\tilde\theta \|   & < \frac{\mathrm c(\tau,\lambda)\sqrt{E[\rho_{T_\beta}]}}{2} + \frac{2\exp(- \beta q_\epsilon )}{\rho(B_{\frac{r_\epsilon}{2}}(\tilde\theta,0))\exp(-p T^*)}\sqrt{E[\rho_{T_\beta}]}\\
         & \leq c(\tau,\lambda) \sqrt{E[\rho_{T_\beta}]} = C(T_\beta),
     \end{aligned}
 \end{equation*}
 where the first inequality in the last line holds by the choice of $\beta$. This establishes the desired contradiction, against the consequence of the continuity of the mappings $t\to E[\rho_t]$ and $t\to \|\mathcal M _\beta[\mu_t] -\tilde\theta\|$.
\end{proof}

\section{Simulation Details}
\label{sec:simulation_datail}
\subsection{Linear-quadratic-Gaussian Control Problem}
We begin by considering a classical LQG control problem, where the state dynamics is governed by: \begin{equation*}
    \intd \bx_t = 2\ctrl_t \intd t + \sqrt{2} \intd  W_t,
\end{equation*}
incorporating $t \in [0, T]$ and $\bx_0 = x$. The cost functional is given by 
$J({\ctrl_t}) = \mathbb{E}\l [
\int^ T_ 0 \|\ctrl_r\|^2 dt + g(X_T )\r  ]$. Here, the state process $\bx_t$ is a $d$-dimensional vector, while the action process $\ctrl_t$ is a $d$-dimensional vector-valued function. The value function $u$ can be defined as
\begin{equation*}
\label{equ:value}
    u(t,\bx) = \inf_{\ctrl}  \mathbb{E}\l [\int_t^T f(t,\bx_t,\ctrl_t,t)\intd t + g(\bx_T)\l |\bx_t= \bx \r  .\r  ].
\end{equation*}
By solving the Hamilton–Jacobi–Bellman equation for $u$, one can derive an explicit solution with the terminal condition \( u(T, x) = g(x) \), given by
\[
u(t, \bx) = -\ln \l (   \mathbb{E} \l [ \exp \l (  - g \l (  \bx + \sqrt{2} \mathbf W_{T-t} \r )  \r )  \r  ] \r ).
\]
\subsection{Ginzburg-Landau Model}
In this model, superconducting electrons are described by a ``macroscopic” wavefunction, $\varphi(z)$, with Landau free energy $\frac{\mu}{4}\int_0^1 |1-\varphi(z)^2|^2 \intd z$. In order to add fluctuations (local variations in the wavefunction) to this model, Ginzburg suggested adding a term
proportional to $|\nabla_z \varphi(z)|^2$, which can be interpolated as the kinetic energy term in quantum
mechanics or the lowest order fluctuation term allowed by the symmetry of the order parameter.
Adding this term to the
free energy, we have the Ginzburg-Landau theory in zero field,
\begin{equation*}
    U[\varphi] = \frac{\lambda}{2}\int_0^1|\nabla_z \varphi(z)|^2 \intd z  + \frac{\mu}{4}\int_0^1 |1-\varphi(z)^2|^2 \intd z.
\end{equation*}
Upon discretizing the space into $d+1$ points, the potential is defined as
\begin{equation*}
    U(\varphi) = U(x_1,\cdots,x_d) := \l [
    \frac{\lambda}{2} \sum_{i=1}^{d+1} \l ( \frac{x_i-x_{i-1}}{h}\r ) ^2 + \frac{\mu}{4}\sum_{i=1}^d (1-x_i^2)^2
    \r  ]h,
\end{equation*}
where $x_i = \varphi(\frac{i}{d+1})$ for $i=0,\cdots,d+1$ and $x_0 = x_{d+1} = 0$. 
The dynamics is given by 
\begin{equation*}
    \intd \bx_t = \mathbf{b}(\bx_t,\alpha_t) \intd t + \sqrt{2} \intd \mathbf W_t,
\end{equation*}
where the drift term is defined as
$$
b(\bx,a) = -\nabla_\bx U(\bx) + 2 \alpha \boldsymbol\omega.
$$
Here the potential is defined as
\begin{equation*}
    U(\varphi) = U(x_1,\cdots,x_d) := \l [
    \frac{\lambda}{2} \sum_{i=1}^{d+1} \l ( \frac{x_i-x_{i-1}}{h}\r ) ^2 + \frac{\mu}{4}\sum_{i=1}^d (1-x_i^2)^2
    \r  ]h,
\end{equation*}
and $\alpha$ is a scalar-valued function, represents the strength of the external field and the vector $\boldsymbol\omega$ is a $d$-dimensional vector represents the domian of the external field applied, the $i\mhyphen$th element takes the value of $1$ if the condition $\frac{i}{d+1} \in [0.25,0.6]$ is satisfied, and $0$ under other circumstances.
The cost functional is defined 
\begin{equation*}
    J[\alpha] = \mathbb{E}\l [\int_0^T \frac{1}{d}\|\bx_t\|^2 + \|\alpha\| \intd t +  \frac{10}{d}\|\bx_T\|^2
 \r  ].
\end{equation*}

\subsection{Systemic Risk Mean Field Control}

We describe this problem as a network of $N$ banks, where $x_i$ denotes the logarithm of the cash reserves of the  $i$-th bank. The following model introduces borrowing and lending between banks, given by:
\begin{equation*}
    \intd \bx_t  =  [\kappa(\bar{\bx}_t- \bx_t) + \ctrl_t] \intd t + \sigma \intd \mathbf W_t,
\end{equation*}
where $\bar{\bx}_t  = \frac{1}{n} \sum_{i=1}^n \bx^i_t$ represents the average logarithm of the cash reserves across all banks. The control of the representative bank, i.e., the amount lent or borrowed at time $t$ is denoted by $\ctrl_t$. Based on the Almgren-Chriss linear price impact model, the running cost and terminal cost are given by: 
\begin{equation*}
    f(x,\bar{x},\alpha) = \frac{1}{2}\alpha^2 - q \alpha (\bar{x}-x) + \frac{\eta}{2}(\bar{x}-x)^2, g(x,\bar{x}) = \frac{c}{2}(\bar{x}-x)^2,
\end{equation*}
where $\eta$ and $c$ balance the individual bank's behavior with the average behavior of the other banks. $q$ weights the contribution of the components and helps to determine the sign of the control (i.e., whether to borrow or lend).
Specifically, if the logarithmic cash reserve of an individual bank is smaller than the empirical mean, the bank will seek to borrow, choosing $\ctrl_t > 0$, and vice versa. 
We test the performance of our method with parameters $c=2$, $k=0.6$, and $\eta=2$.

\section{Neural Network Structure}\label{sec:network}
In this subsection, we briefly illustrate the network structure used to model the action function. For the LQG and Ginzburg-Landau model problems, we employ traditional fully connected neural networks with a depth of 5 layers and a width of $5d$, where $d$ represents the dimension of the problem.
In the mean-field control problem, we use the cylindrical type mean field neural network structure proposed in \cite{pham2024mean}, where the control $\alpha(t,\bx)$ is parameterized as 
\begin{equation*}
    \ctrl^i(t,\bx;\theta)  =   \Psi\l [x^i, \frac{1}{n}\sum_{i=1}^n \psi(x_i;\theta_2);\theta_1\r  ],
\end{equation*}
where $\theta = \{\theta_1,\theta_2\}$. The advantage of this type of network is its extendibility, i.e, as demonstrated in the numerical results, control policies trained on a small number of agents $N$ can be effectively applied to problems with different values of 
$n$.



\bibliographystyle{plainnat}  
\bibliography{main}

\begin{thebibliography}{10}

\bibitem{stengel1986stochastic}
Robert~F Stengel.
\newblock {\em Stochastic optimal control: theory and application}.
\newblock John Wiley \& Sons, Inc., 1986.

\bibitem{fleming2012deterministic}
Wendell~H Fleming and Raymond~W Rishel.
\newblock {\em Deterministic and stochastic optimal control}, volume~1.
\newblock Springer Science \& Business Media, 2012.

\bibitem{pham2009continuous}
Huy{\^e}n Pham.
\newblock {\em Continuous-time stochastic control and optimization with
  financial applications}, volume~61.
\newblock Springer Science \& Business Media, 2009.

\bibitem{fleming2004stochastic}
Wendell~H Fleming and Jerome~L Stein.
\newblock Stochastic optimal control, international finance and debt.
\newblock {\em Journal of Banking \& Finance}, 28(5):979--996, 2004.

\bibitem{carmona2003pricing}
Ren{\'e} Carmona and Valdo Durrleman.
\newblock Pricing and hedging spread options.
\newblock {\em Siam Review}, 45(4):627--685, 2003.

\bibitem{cousin2011mean}
Areski Cousin, St{\'e}phane Cr{\'e}pey, Olivier Gu{\'e}ant, David Hobson,
  Monique Jeanblanc, Jean-Michel Lasry, Jean-Paul Laurent, Pierre-Louis Lions,
  Peter Tankov, Olivier Gu{\'e}ant, et~al.
\newblock Mean field games and applications.
\newblock {\em Paris-Princeton lectures on mathematical finance 2010}, pages
  205--266, 2011.

\bibitem{lachapelle2016efficiency}
Aim{\'e} Lachapelle, Jean-Michel Lasry, Charles-Albert Lehalle, and
  Pierre-Louis Lions.
\newblock Efficiency of the price formation process in presence of high
  frequency participants: a mean field game analysis.
\newblock {\em Mathematics and Financial Economics}, 10:223--262, 2016.

\bibitem{cardaliaguet2018mean}
Pierre Cardaliaguet and Charles-Albert Lehalle.
\newblock Mean field game of controls and an application to trade crowding.
\newblock {\em Mathematics and Financial Economics}, 12:335--363, 2018.

\bibitem{gueant2009mean}
Olivier Gu{\'e}ant.
\newblock Mean field games and applications to economics, 2009.

\bibitem{gomes2015economic}
Diogo~A Gomes, Levon Nurbekyan, and Edgard~A Pimentel.
\newblock Economic models and mean-field games theory.
\newblock {\em Publicaoes Matematicas, IMPA, Rio, Brazil}, 2015.

\bibitem{gueant2010mean}
O~Gu{\'e}ant, JM~Lasry, and PL~Lions.
\newblock Mean field games and applications. paris-princeton lectures on
  mathematical finance.
\newblock {\em Lect. Notes Math}, 2011:205--266, 2010.

\bibitem{achdou2014partial}
Yves Achdou, Francisco~J Buera, Jean-Michel Lasry, Pierre-Louis Lions, and
  Benjamin Moll.
\newblock Partial differential equation models in macroeconomics.
\newblock {\em Philosophical Transactions of the Royal Society A: Mathematical,
  Physical and Engineering Sciences}, 372(2028):20130397, 2014.

\bibitem{achdou2022income}
Yves Achdou, Jiequn Han, Jean-Michel Lasry, Pierre-Louis Lions, and Benjamin
  Moll.
\newblock Income and wealth distribution in macroeconomics: A continuous-time
  approach.
\newblock {\em The review of economic studies}, 89(1):45--86, 2022.

\bibitem{welch2019describing}
PM~Welch, K{\O}~Rasmussen, and Cynthia~F Welch.
\newblock Describing nonequilibrium soft matter with mean field game theory.
\newblock {\em The Journal of Chemical Physics}, 150(17), 2019.

\bibitem{holdijk2024stochastic}
Lars Holdijk, Yuanqi Du, Ferry Hooft, Priyank Jaini, Berend Ensing, and Max
  Welling.
\newblock Stochastic optimal control for collective variable free sampling of
  molecular transition paths.
\newblock {\em Advances in Neural Information Processing Systems}, 36, 2024.

\bibitem{lachapelle2011mean}
Aim{\'e} Lachapelle and Marie-Therese Wolfram.
\newblock On a mean field game approach modeling congestion and aversion in
  pedestrian crowds.
\newblock {\em Transportation research part B: methodological},
  45(10):1572--1589, 2011.

\bibitem{aurell2018mean}
Alexander Aurell and Boualem Djehiche.
\newblock Mean-field type modeling of nonlocal crowd aversion in pedestrian
  crowd dynamics.
\newblock {\em SIAM Journal on Control and Optimization}, 56(1):434--455, 2018.

\bibitem{achdou2019mean}
Yves Achdou and Jean-Michel Lasry.
\newblock Mean field games for modeling crowd motion.
\newblock {\em Contributions to partial differential equations and
  applications}, pages 17--42, 2019.

\bibitem{hu2024recent}
Ruimeng Hu and Mathieu Lauri\`ere.
\newblock Recent developments in machine learning methods for stochastic
  control and games.
\newblock {\em Numerical Algebra, Control and Optimization}, 14(3):435--525,
  2024.

\bibitem{richardson2006numerical}
Steven Richardson and Song Wang.
\newblock Numerical solution of hamilton--jacobi--bellman equations by an
  exponentially fitted finite volume method.
\newblock {\em Optimization}, 55(1-2):121--140, 2006.

\bibitem{wang2003numerical}
Song Wang, Les~S Jennings, and Kok~Lay Teo.
\newblock Numerical solution of hamilton-jacobi-bellman equations by an upwind
  finite volume method.
\newblock {\em Journal of Global Optimization}, 27:177--192, 2003.

\bibitem{beard1997galerkin}
Randal~W Beard, George~N Saridis, and John~T Wen.
\newblock Galerkin approximations of the generalized {Hamilton-Jacobi-Bellman}
  equation.
\newblock {\em Automatica}, 33(12):2159--2177, 1997.

\bibitem{bea1998successive}
Randal~W Beard.
\newblock Successive galerkin approximation algorithms for nonlinear optimal
  and robust control.
\newblock {\em International Journal of Control}, 71(5):717--743, 1998.

\bibitem{forsyth2007numerical}
Peter~A Forsyth and George Labahn.
\newblock Numerical methods for controlled hamilton-jacobi-bellman pdes in
  finance.
\newblock {\em Journal of Computational Finance}, 11(2):1, 2007.

\bibitem{Li2024value}
Pengyi Li, Jianye Hao, Hongyao Tang, Yan Zheng, and Fazl Barez.
\newblock Value-evolutionary-based reinforcement learning.
\newblock In {\em Proceedings of the 41st International Conference on Machine
  Learning}, pages 27875--27889, 2024.

\bibitem{Lien2024Enhancing}
Yun-Hsuan Lien, Ping-Chun Hsieh, Tzu-Mao Li, and Yu-Shuen Wang.
\newblock Enhancing value function estimation through first-order state-action
  dynamics in offline reinforcement learning.
\newblock In {\em Proceedings of the 41st International Conference on Machine
  Learning}, pages 29782--29794, 2024.

\bibitem{Obando2024value}
Johan~Samir Obando~Ceron, Aaron Courville, and Pablo~Samuel Castro.
\newblock In value-based deep reinforcement learning, a pruned network is a
  good network.
\newblock In {\em Proceedings of the 41st International Conference on Machine
  Learning}, pages 38495--38519, 2024.

\bibitem{zhang2024model}
Renhao Zhang, Haotian Fu, Yilin Miao, and George Konidaris.
\newblock Model-based reinforcement learning for parameterized action spaces.
\newblock In {\em Proceedings of the 41st International Conference on Machine
  Learning}, pages 58935--58954, 2024.

\bibitem{mou2024bellman}
Wenlong Mou and Yuhua Zhu.
\newblock On bellman equations for continuous-time policy evaluation i:
  discretization and approximation.
\newblock {\em arXiv preprint arXiv:2407.05966}, 2024.

\bibitem{e2017deep}
Weinan E, Jiequn Han, and Arnulf Jentzen.
\newblock Deep learning-based numerical methods for high-dimensional parabolic
  partial differential equations and backward stochastic differential
  equations.
\newblock {\em Communications in mathematics and statistics}, 5(4):349--380,
  2017.

\bibitem{han2018solving}
Jiequn Han, Arnulf Jentzen, and Weinan E.
\newblock Solving high-dimensional partial differential equations using deep
  learning.
\newblock {\em Proceedings of the National Academy of Sciences},
  115(34):8505--8510, 2018.

\bibitem{nusken2021solving}
Nikolas N{\"u}sken and Lorenz Richter.
\newblock Solving high-dimensional hamilton--jacobi--bellman pdes using neural
  networks: perspectives from the theory of controlled diffusions and measures
  on path space.
\newblock {\em Partial differential equations and applications}, 2(4):48, 2021.

\bibitem{pham2021neural}
Huy{\^e}n Pham, Xavier Warin, and Maximilien Germain.
\newblock Neural networks-based backward scheme for fully nonlinear pdes.
\newblock {\em SN Partial Differential Equations and Applications}, 2(1):16,
  2021.

\bibitem{lyu2023construction}
Liyao Lyu and Huan Lei.
\newblock Construction of coarse-grained molecular dynamics with many-body
  non-markovian memory.
\newblock {\em Physical Review Letters}, 131(17):177301, 2023.

\bibitem{lu2019nonparametric}
Fei Lu, Ming Zhong, Sui Tang, and Mauro Maggioni.
\newblock Nonparametric inference of interaction laws in systems of agents from
  trajectory data.
\newblock {\em Proceedings of the National Academy of Sciences},
  116(29):14424--14433, 2019.

\bibitem{agrawalpolicy}
Shubhada Agrawal, LA~Prashanth, and Siva~Theja Maguluri.
\newblock Policy evaluation for variance in average reward reinforcement
  learning.
\newblock In {\em Forty-first International Conference on Machine Learning},
  2024.

\bibitem{chenaccelerated}
Yen-Ju Chen, Nai-Chieh Huang, Ching-pei Lee, and Ping-Chun Hsieh.
\newblock Accelerated policy gradient: On the convergence rates of the nesterov
  momentum for reinforcement learning.
\newblock In {\em Forty-first International Conference on Machine Learning},
  2024.

\bibitem{chenE32024}
Dingyang Chen and Qi~Zhang.
\newblock ${\rm e}(3)$-equivariant actor-critic methods for cooperative
  multi-agent reinforcement learning.
\newblock In {\em Proceedings of the 41st International Conference on Machine
  Learning}, 2024.

\bibitem{Daisafe2024}
Juntao Dai, Yaodong Yang, Qian Zheng, and Gang Pan.
\newblock Safe reinforcement learning using finite-horizon gradient-based
  estimation.
\newblock In {\em Proceedings of the 41st International Conference on Machine
  Learning}, 2024.

\bibitem{Hisaki2024Average}
Yukinari Hisaki and Isao Ono.
\newblock {RVI}-{SAC}: Average reward off-policy deep reinforcement learning.
\newblock In {\em Proceedings of the 41st International Conference on Machine
  Learning}, pages 18352--18373, 2024.

\bibitem{Hong2024model}
Mao Hong, Zhengling Qi, and Yanxun Xu.
\newblock Model-based reinforcement learning for confounded {POMDP}s.
\newblock In {\em Proceedings of the 41st International Conference on Machine
  Learning}, pages 18668--18710, 2024.

\bibitem{Hu2024Q}
Shengchao Hu, Ziqing Fan, Chaoqin Huang, Li~Shen, Ya~Zhang, Yanfeng Wang, and
  Dacheng Tao.
\newblock Q-value regularized transformer for offline reinforcement learning.
\newblock In {\em Proceedings of the 41st International Conference on Machine
  Learning}, pages 19165--19181, 2024.

\bibitem{Park2024Max}
Giseung Park, Woohyeon Byeon, Seongmin Kim, Elad Havakuk, Amir Leshem, and
  Youngchul Sung.
\newblock The max-min formulation of multi-objective reinforcement learning:
  From theory to a model-free algorithm.
\newblock In {\em Proceedings of the 41st International Conference on Machine
  Learning}, pages 39616--39642, 2024.

\bibitem{tang2024solving}
Xun Tang, Leah Collis, and Lexing Ying.
\newblock Solving high-dimensional {Kolmogorov} backward equations with
  operator-valued functional hierarchical tensor for markov operators.
\newblock {\em arXiv preprint arXiv:2404.08823}, 2024.

\bibitem{mnih2015human}
Volodymyr Mnih, Koray Kavukcuoglu, David Silver, Andrei~A Rusu, Joel Veness,
  Marc~G Bellemare, Alex Graves, Martin Riedmiller, Andreas~K Fidjeland, Georg
  Ostrovski, et~al.
\newblock Human-level control through deep reinforcement learning.
\newblock {\em nature}, 518(7540):529--533, 2015.

\bibitem{heess2017emergence}
Nicolas Heess, Dhruva Tb, Srinivasan Sriram, Jay Lemmon, Josh Merel, Greg
  Wayne, Yuval Tassa, Tom Erez, Ziyu Wang, SM~Eslami, et~al.
\newblock Emergence of locomotion behaviours in rich environments.
\newblock {\em arXiv preprint arXiv:1707.02286}, 2017.

\bibitem{schulman2015high}
John Schulman, Philipp Moritz, Sergey Levine, Michael Jordan, and Pieter
  Abbeel.
\newblock High-dimensional continuous control using generalized advantage
  estimation.
\newblock {\em arXiv preprint arXiv:1506.02438}, 2015.

\bibitem{schulman2017proximal}
John Schulman, Filip Wolski, Prafulla Dhariwal, Alec Radford, and Oleg Klimov.
\newblock Proximal policy optimization algorithms.
\newblock {\em arXiv preprint arXiv:1707.06347}, 2017.

\bibitem{schulman2015trust}
John Schulman, Sergey Levine, Pieter Abbeel, Michael Jordan, and Philipp
  Moritz.
\newblock Trust region policy optimization.
\newblock In Francis Bach and David Blei, editors, {\em Proceedings of the 32nd
  International Conference on Machine Learning}, volume~37 of {\em Proceedings
  of Machine Learning Research}, pages 1889--1897, Lille, France, 07--09 Jul
  2015.

\bibitem{silver2014deterministic}
David Silver, Guy Lever, Nicolas Heess, Thomas Degris, Daan Wierstra, and
  Martin Riedmiller.
\newblock Deterministic policy gradient algorithms.
\newblock In {\em International conference on machine learning}, pages
  387--395. Pmlr, 2014.

\bibitem{lillicrap2015continuous}
TP~Lillicrap.
\newblock Continuous control with deep reinforcement learning.
\newblock {\em arXiv preprint arXiv:1509.02971}, 2015.

\bibitem{haarnoja2018soft}
Tuomas Haarnoja, Aurick Zhou, Pieter Abbeel, and Sergey Levine.
\newblock Soft actor-critic: Off-policy maximum entropy deep reinforcement
  learning with a stochastic actor.
\newblock In {\em International conference on machine learning}, pages
  1861--1870. PMLR, 2018.

\bibitem{haarnoja2018learning}
Tuomas Haarnoja, Sehoon Ha, Aurick Zhou, Jie Tan, George Tucker, and Sergey
  Levine.
\newblock Learning to walk via deep reinforcement learning.
\newblock {\em arXiv preprint arXiv:1812.11103}, 2018.

\bibitem{jia2022actor-critic}
Yanwei Jia and Xun~Yu Zhou.
\newblock Policy gradient and actor-critic learning in continuous time and
  space: Theory and algorithms.
\newblock {\em Journal of Machine Learning Research}, 23(275):1--50, 2022.

\bibitem{jia2022temporal-difference}
Yanwei Jia and Xun~Yu Zhou.
\newblock Policy evaluation and temporal-difference learning in continuous time
  and space: A martingale approach.
\newblock {\em Journal of Machine Learning Research}, 23(154):1--55, 2022.

\bibitem{gu2021mean}
Haotian Gu, Xin Guo, Xiaoli Wei, and Renyuan Xu.
\newblock Mean-field controls with q-learning for cooperative marl: convergence
  and complexity analysis.
\newblock {\em SIAM Journal on Mathematics of Data Science}, 3(4):1168--1196,
  2021.

\bibitem{carmona2023model}
Ren{\'e} Carmona, Mathieu Lauri{\`e}re, and Zongjun Tan.
\newblock Model-free mean-field reinforcement learning: mean-field mdp and
  mean-field q-learning.
\newblock {\em The Annals of Applied Probability}, 33(6B):5334--5381, 2023.

\bibitem{hua2024simulation}
Mengjian Hua, Matthieu Lauri\`ere, and Eric Vanden-Eijnden.
\newblock A simulation-free deep learning approach to stochastic optimal
  control.
\newblock {\em arXiv preprint arXiv:2410.05163}, 2024.

\bibitem{chen2022consensus}
Jingrun Chen, Shi Jin, and Liyao Lyu.
\newblock A consensus-based global optimization method with adaptive momentum
  estimation.
\newblock {\em Communications in Computational Physics}, 31(4):1296--1316,
  2022.

\bibitem{fornasier2024consensus2}
Massimo Fornasier, Timo Klock, and Konstantin Riedl.
\newblock Consensus-based optimization methods converge globally.
\newblock {\em SIAM Journal on Optimization}, 34(3):2973--3004, 2024.

\bibitem{Ganapathi2024Confidence}
Sriram Ganapathi~Subramanian, Guiliang Liu, Mohammed Elmahgiubi, Kasra Rezaee,
  and Pascal Poupart.
\newblock Confidence aware inverse constrained reinforcement learning.
\newblock In {\em Proceedings of the 41st International Conference on Machine
  Learning}, pages 14491--14512, 2024.

\bibitem{Hong2024primal}
Kihyuk Hong and Ambuj Tewari.
\newblock A primal-dual algorithm for offline constrained reinforcement
  learning with linear {MDP}s.
\newblock In {\em Proceedings of the 41st International Conference on Machine
  Learning}, pages 18711--18737, 2024.

\bibitem{Qiao2024Near}
Dan Qiao and Yu-Xiang Wang.
\newblock Near-optimal reinforcement learning with self-play under adaptivity
  constraints.
\newblock In {\em Proceedings of the 41st International Conference on Machine
  Learning}, pages 41430--41455, 2024.

\bibitem{sun2024Constrained}
Zhongchang Sun, Sihong He, Fei Miao, and Shaofeng Zou.
\newblock Constrained reinforcement learning under model mismatch.
\newblock In {\em Proceedings of the 41st International Conference on Machine
  Learning}, pages 47017--47032, 2024.

\bibitem{Wang2024Probabilistic}
Yanran Wang, Qiuchen Qian, and David Boyle.
\newblock Probabilistic constrained reinforcement learning with formal
  interpretability.
\newblock In {\em Proceedings of the 41st International Conference on Machine
  Learning}, pages 51303--51327, 2024.

\bibitem{carrillo2018analytical}
Jos\'e~A Carrillo, Young-Pil Choi, Claudia Totzeck, and Oliver Tse.
\newblock An analytical framework for consensus-based global optimization
  method.
\newblock {\em Mathematical Models and Methods in Applied Sciences},
  28(06):1037--1066, 2018.

\bibitem{durrett2018stochastic}
Richard Durrett.
\newblock {\em Stochastic calculus: a practical introduction}.
\newblock CRC press, 2018.

\bibitem{gilbarg2001elliptic}
David Gilbarg, Neil~S Trudinger, David Gilbarg, and NS~Trudinger.
\newblock {\em Elliptic partial differential equations of second order}, volume
  224.
\newblock Springer, 2001.

\bibitem{Arnold1976stochastic}
Ludwig Arnold.
\newblock Stochastic differential equations: Theory and applications, 1976.

\bibitem{bass2011stochastic}
Richard~F Bass.
\newblock {\em Stochastic processes}, volume~33.
\newblock Cambridge University Press, 2011.

\bibitem{Huang2022mean}
Hui Huang and Jinniao Qiu.
\newblock On the mean-field limit for the consensus-based optimization.
\newblock {\em Mathematical Methods in the Applied Sciences},
  45(12):7814--7831, 2022.

\bibitem{billingsley2013convergence}
Patrick Billingsley.
\newblock {\em Convergence of probability measures}.
\newblock John Wiley \& Sons, 2013.

\bibitem{pham2024mean}
Huy{\^e}n Pham and Xavier Warin.
\newblock Mean-field neural networks-based algorithms for mckean-vlasov control
  problems.
\newblock {\em Journal of Machine Learning}, 3(2):176--214, 2024.

\end{thebibliography}
\end{document}